\documentclass[11pt]{amsart}

\usepackage{amsmath,amssymb,amsfonts,amsthm,mathtools,mathrsfs}
\usepackage{hyperref}
\usepackage[hmargin=3cm,vmargin=3.2cm]{geometry}
\usepackage{needspace}
\usepackage{xcolor}

\hypersetup{colorlinks=true,
  linkcolor=blue,
  citecolor=blue,
  urlcolor=blue
}

\newtheorem{theorem}{Theorem}[section]
\newtheorem{proposition}[theorem]{Proposition}
\newtheorem{lemma}[theorem]{Lemma}
\newtheorem{corollary}[theorem]{Corollary}
\newtheorem{definition}[theorem]{Definition}
\newtheorem{remark}[theorem]{Remark}
\newtheorem{example}[theorem]{Example}

\theoremstyle{remark}
\newtheorem*{remark*}{Remark}
\theoremstyle{plain}

\numberwithin{equation}{section}

\newcommand{\K}{\mathbb K}

\newcommand{\g}{\mathfrak g}
\newcommand{\h}{\mathfrak h}

\newcommand{\bLie}{\mathfrak b}

\newcommand{\ad}{\operatorname{ad}}
\newcommand{\Der}{\operatorname{Der}}
\newcommand{\Dera}{\operatorname{Der}_{\mathrm a}}

\newcommand{\Imop}{\operatorname{Im}}
\newcommand{\Ker}{\operatorname{Ker}}
\newcommand{\rank}{\operatorname{rank}}
\newcommand{\id}{\operatorname{id}}

\newcommand{\Ann}{\operatorname{Ann}}
\newcommand{\Hom}{\operatorname{Hom}}
\newcommand{\End}{\operatorname{End}}
\newcommand{\Span}{\operatorname{span}}

\title{Poisson and Jacobi structures on quadratic Lie algebras}
\author{Alessandro Arsie and Fabricio Valencia}
\address{}
\date{\today}

\address{Alessandro Arsie - Department of Mathematics and Statistics, The University of Toledo, 2801 W. Bancroft St., Toledo, 43606 Ohio - United States.
	\vspace{0.5cm}
	\newline  
	Fabricio Valencia - Instituto de Matem\'atica, Estat\'istica e Ciência da Computação, Universidade de S\~ao Paulo, Rua do Mat\~ao 1010, Cidade Universit\'aria, S\~ao Paulo, 05508-090 S\~ao Paulo - Brazil.
	\vspace{0.5cm}
	\newline  
	\phantom{xx}
	alessandro.arsie@utoledo.edu, fabricio.valencia@ime.usp.br}

\begin{document}

\begin{abstract}
We study Poisson and Jacobi structures compatible with quadratic structures on finite-dimensional Lie algebras, mostly over fields of characteristic zero.  Building on the theory of generalized double extensions, we give necessary and sufficient reconstruction criteria for the skew operators preserving a marked isotropic ideal, allowing degenerate operators, a nonabelian quadratic quotient, and the full generalized extension data.  Recalling the characterization of quadratic Poisson structures in terms of skew-symmetric solutions of the classical Yang--Baxter equation,  we combine the ambient bracket reconstruction of Kostant and Meyer with prescribed symplectic leaf data to reconstruct Poisson structures whose characteristic subalgebra is nondegenerate for the quadratic form.  We express the standard algebraic Jacobi equations through the invariant form and obtain analogous reconstruction results for nondegenerate characteristic subalgebras. Particular attention is given to the contact and locally conformally symplectic cases. A common cochain framework identifies the remaining scalar obstruction in third cohomology in the Poisson case and in the Jacobi branch whose Euler element belongs to the marked ideal, with coefficients in a one dimensional module in the latter case. Our reconstruction results yield systematic procedures for constructing quadratic Poisson structures in dimension $n+2k$ from those in dimension $n$, and recover the earlier quadratic symplectic extension formulas when the Poisson operator is invertible. Explicit examples illustrate the applicability of the theory, including oscillator algebras, the five-dimensional Lie algebra $A_{5,3}$, and Poisson structures associated with the curvature tensor of a bi-invariant metric.
\end{abstract}

\maketitle
\tableofcontents
\section{Introduction}

A Lie algebra is called \emph{quadratic} if it admits a symmetric, invariant, nondegenerate bilinear form. Classical examples include Lie algebras of compact Lie groups, semisimple Lie algebras, and cotangent Lie algebras of Lie groups. These algebraic objects have been widely studied because of their connections with problems in geometry and mathematical physics. In particular, quadratic Lie algebras provide the algebraic framework for bi-invariant pseudo-Riemannian metrics on Lie groups and play an important role in conformal field theory \cite{FigueroaStanciu1994,FigueroaStanciu1996}. Their study extends beyond the semisimple setting, encompassing a rich class of nonsemisimple Lie algebras. From an algebraic perspective, the ordinary double extension and $T^*$-extension theorems provide systematic approaches to their construction and classification, further developed through generalized double extensions \cite{BajoBenayadiBordemann2007,Bordemann1997,MedinaRevoy1985}. These connections have also motivated the study of additional geometric structures compatible with quadratic Lie algebra structures, including symplectic structures, as investigated in \cite{BajoBenayadiMedina2007,MedinaRevoy1988}.

Motivated by these developments and by the rich geometric framework provided by Poisson and Jacobi manifolds,  we undertake a systematic study of these structures on quadratic Lie algebras, focusing in particular on their behavior under extension procedures. Since Poisson structures encompass possibly degenerate Hamiltonian systems, while Jacobi structures also provide a natural framework for contact and locally conformal symplectic geometry, their study offers a broad perspective on the interaction between Lie theoretic and geometric frameworks.  In the quadratic setting, the invariant bilinear form provides a means of translating geometric compatibility conditions into algebraic ones, opening the way to systematic construction methods. Our main contribution is a reconstruction theory for Poisson and Jacobi structures preserving a specified isotropic ideal $I$, without assuming that the skew operator is invertible. The criteria allow a nonabelian quadratic quotient $I^\perp/I$ and the full generalized extension data. They extend the symplectic constructions of \cite{BajoBenayadiMedina2007} and the symplectic models with an abelian quadratic quotient in \cite{Fischer2019}, and identify obstructions in third cohomology for Poisson structures and for Jacobi structures whose Euler element belongs to the ideal. Whenever the ground field is the real numbers, the constructions and results obtained throughout this work concern left-invariant Poisson and Jacobi structures on the simply connected Lie groups integrating the corresponding quadratic Lie algebras. The case of left invariant Poisson structures on classical non-compact simple Lie groups was studied by Agaoka in \cite{Agaoka2000}.
 
The relevance of quadratic Lie algebras also emerges in recent developments in Poisson geometry and shifted symplectic geometry. In particular, \(2\)-shifted symplectic groups differentiate to quadratic Lie algebras, while, conversely, the van Est theorem ensures that every quadratic structure on the Lie algebra of a simply connected Lie group integrates to a \(2\)-shifted symplectic form, unique up to cohomology \cite{AlvarezBursztynCueca2026,CuecaZhu2023}. These developments highlight the role of quadratic Lie algebras in connecting infinitesimal algebraic data with global geometric structures and further motivate the study of other geometric data compatible with quadratic forms, particularly the Poisson and Jacobi frameworks investigated in this work.

The paper is organized as follows.

\begin{itemize}
\item Section~\ref{sec:double-ext} adapts the generalized double-extension construction introduced by Bajo--Benayadi--Bordemann \cite{BajoBenayadiBordemann2007} from quadratic Lie superalgebras to classical quadratic Lie algebras, providing concise proofs of the main results. The generalized bracket, reconstruction from an isotropic ideal, and inclusion of ordinary and $T^*_\theta$ extensions are already contained in \cite{BajoBenayadiBordemann2007}. Here we make explicit the recognition relative to the marked ideal and the changes of the extension data, skew operator, and Euler element under an adapted Witt splitting, as developed in Subsections~\ref{subsec:gd-recognition} and~\ref{subsec:gd-splittings}. The reconstruction from a prescribed isotropic ideal is proved over any field of characteristic different from $2$, and the ordinary presentation keeps the additional invariant form, which may be degenerate. This framework also provides the algebraic foundation for our subsequent study of Poisson and Jacobi structures on quadratic Lie algebras.

\item Section~\ref{sec:poisson-preliminaries} establishes preliminary results for the study of Poisson and Jacobi structures compatible with quadratic structures on Lie algebras, building on ideas and partial results from \cite{BajoBenayadiMedina2007,BoucettaMedina2010,MedinaRevoy1988}. In particular, we recall the  one-to-one correspondence between quadratic Poisson structures on Lie algebras and skew-symmetric solutions of the classical Yang--Baxter equation. Building on work by Kostant \cite{Kostant1999} and Meyer \cite{Meyer2021}, we explicitly describe the symplectic structure induced on the algebraic characteristic leaf and prove a reconstruction theorem when this leaf is itself quadratic, meaning that the restriction of the ambient quadratic form to the leaf is nondegenerate; see Theorem~\ref{thm:leaf-reconstruction}.

\item Section~\ref{sec:jacobi-preliminaries} expresses the established algebraic Jacobi equations and the description of characteristic subalgebras in \cite{IglesiasMarrero2003} using the invariant quadratic form. This gives explicit identities for the skew operator and its Euler element, which are then used in the reconstruction. We distinguish two important cases, corresponding to Jacobi structures arising from contact and from locally conformally symplectic data. We also establish reconstruction theorems based on the nondegeneracy of the characteristic subalgebra, including results tailored to these two special cases; compare Theorems~\ref{thm:jacobi-nd-reconstruction}, \ref{thm:jacobi-lcs-reconstruction}, and \ref{thm:jacobi-c-reconstruction}.

\item Section~\ref{subsec:shared-cochains} develops a cochain framework that provides a cohomological interpretation of the Poisson case and the inner Jacobi branch for generalized double extensions whose skew operator preserves a marked isotropic ideal. This framework also reveals algebraic relations, such as those established in Lemma~\ref{lem:jacobi-transgression}, that emerge when the Lie algebra structure induced by a solution of the classical Yang--Baxter equation is considered and the solution itself defines a Lie algebra morphism between the underlying Lie algebras.

\item Section~\ref{sec:gd-reconstruction} establishes a general reconstruction theorem for Poisson structures on quadratic Lie algebras preserving a marked isotropic ideal, providing a systematic method for constructing quadratic Poisson data of dimension $n+2k$ from those of dimension $n$; see Theorem~\ref{thm:gd-poisson}. Proposition~\ref{prop:gd-bivector} gives an explicit expression for the reconstructed Poisson bivector, while Theorem~\ref{thm:gd-poisson-obstruction} provides a cohomological interpretation of the main algebraic obstruction to the construction. In particular, the Poisson differential induced by the Schouten bracket preserves $\Lambda^\bullet I$. Under the identification $I=\bLie^*$, its restriction coincides with the Chevalley--Eilenberg differential with scalar coefficients of the induced Lie algebra $\bLie_A$, for $A$ being a solution of the classical Yang--Baxter equation. In the invertible case, these equations recover the constructions of \cite{BajoBenayadiMedina2007}, including their scalar obstruction for symplectic $T^*_\theta$ extensions. Here the obstruction also applies to degenerate operators and full generalized extension data. Corollary~\ref{cor:Sol_Nil_Poisson} gives the reconstruction for any quadratic Poisson algebra with a specified nonzero isotropic ideal preserved by its operator, including nilpotent and solvable algebras with such an ideal.

\item Section~\ref{sec:jacobi-double} establishes the analogous reconstruction for Jacobi structures whose skew operator preserves a marked isotropic ideal, treating all three positions of the Euler element in the isotropic flag. When this element belongs to the ideal, the remaining scalar obstruction lies in third cohomology with coefficients in a one dimensional module, as shown in Theorem~\ref{thm:gd-inner-obstruction}. The corresponding nilpotent and solvable cases are treated when a specified nonzero isotropic ideal is preserved. For an indecomposable algebra, the same reconstruction applies when a chosen minimal isotropic ideal is preserved.

\item Section~\ref{sec:poisson-jacobi-examples} tests our theoretical results on some quadratic Lie algebras of geometric and physical interest, including oscillator algebras and \(A_{5,3}\), the five-dimensional indecomposable real Lie algebra underlying a Wess--Zumino--Witten model related to four-dimensional plane-wave string backgrounds \cite{Kehagias1994}. Oscillator algebras are especially suited to this purpose because of their distinguished role among solvable Lie algebras with irreducible bi-invariant Lorentzian metrics, as well as their connections with the harmonic oscillator and mathematical physics \cite{Medina1985,GadeaOubina1999,Levichev1986,Levichev1990,MullerRicci1988,NappiWitten1993,Streater1967}. The oscillator Poisson classification covers repeated and resonant frequencies, extending the generic setting of \cite{BoucettaMedina2010}. We also examine Poisson structures associated with the curvature tensor of a bi-invariant metric, thereby illustrating another interaction between the algebraic properties of quadratic Lie algebras and their differential-geometric structures.
\end{itemize}

\vspace{.2cm}

\noindent {\bf Acknowledgments:} This work was completed while the first named author was on sabbatical and visiting the Department of Mathematics of the University of São Paulo.  We thank the Department for partial support through the Fapesp grant 2022/16097-2, the pleasant atmosphere and the kind hospitality.  We thank also Alessandro Tomasiello for useful discussions concerning examples of generalized double extensions relevant to Theoretical Physics. 

\section{Quadratic generalized double extension revisited}\label{sec:double-ext}

We begin with the generalized double extension of Bajo--Benayadi--Bordemann \cite{BajoBenayadiBordemann2007}. Their construction is formulated for quadratic Lie superalgebras. Here we give the construction and its proofs for ordinary Lie algebras, following \cite{BajoBenayadiBordemann2007}. We then explain how it includes the ordinary double extension of Medina--Revoy \cite{MedinaRevoy1985} and the $T^*_\theta$-extension of Bordemann \cite{Bordemann1997}. This also makes it possible to describe the role of the splitting in each construction. See Remark~\ref{Rmk1_NewInsights} and Subsections~\ref{subsec:gd-recognition} and~\ref{subsec:gd-splittings}.

Unless otherwise stated, all Lie algebras are finite-dimensional over a field $\K$ of characteristic different from $2$. Statements about Lie groups are made over $\K=\mathbb R$.

\begin{definition}
	\normalfont 
A \emph{quadratic Lie algebra} is a pair $(\g,B)$ consisting of a Lie algebra $\g$ and a symmetric nondegenerate bilinear form $B$ satisfying
	\begin{equation}\label{eq:invariant-B}
		B([X,Y],Z)=B(X,[Y,Z]),\qquad X,Y,Z\in \g.
	\end{equation}
\end{definition}

Quadratic Lie algebras are also called \emph{metric} or \emph{self-dual} Lie algebras \cite{Bordemann1997,MedinaRevoy1985}. Over $\K=\mathbb R$, the algebra $\g$ integrates to a simply connected Lie group $G$ with a bi-invariant pseudo-Riemannian metric. The adjoint action of $G$ preserves $B$, and the flat map $B^\flat\colon\g\longrightarrow\g^*$ induces a $G$-equivariant diffeomorphism between the adjoint and coadjoint orbits of $G$.

\subsection{The generalized bracket}
\label{subsec:double-def} 

The construction below produces a quadratic Lie algebra $(\tilde{\mathfrak g},T)$ of dimension $n+2k$ from a quadratic Lie algebra $(\mathfrak g,B)$ of dimension $n$. In particular, it includes the constructions used to obtain nilpotent and solvable quadratic Lie algebras.

Let $(\tilde{\mathfrak g},T)$ be a quadratic Lie algebra. A subspace $I\subset \tilde{\mathfrak g}$ is \emph{isotropic} if $T$ vanishes on $I\times I$, or equivalently if $I\subseteq I^\perp$. It is \emph{Lagrangian} if $I=I^\perp$. We call an ideal \emph{marked} when it is specified as part of the data.  In general, we will assume that this marked ideal $I$ is not zero. When the zero ideal is allowed, it gives the identity presentation in the reconstruction below. A nontrivial increase in dimension requires a nonzero ideal. An isomorphism must then carry it to the marked ideal in the target. An isotropic ideal determines the flag
\[
I\subseteq I^\perp\subseteq\tilde{\mathfrak g}.
\]
The quotient $I^\perp/I$ carries the form induced by $T$, while $\tilde{\mathfrak g}/I^\perp$ pairs with $I$. These two observations lead to the generalized double extension. The bracket and the reconstruction from an isotropic ideal are the versions for ordinary Lie algebras of the results of Bajo--Benayadi--Bordemann (see \cite{BajoBenayadiBordemann2007}, Theorems~2.3 and~3.1). We give direct proofs over any field of characteristic zero and explain how the splitting and the additional invariant form enter the ordinary construction. We do not assume that $\g$ is abelian.

We first describe the bracket and the conditions it must satisfy. Let $\bLie$ be a Lie algebra and let $(\g,B)$ be a quadratic Lie algebra. Write
\begin{equation}\label{eq:Dera}
	\Dera(\g)=\Der(\g)\cap\mathfrak o(\g,B).
\end{equation}
for the Lie algebra of \emph{skew-symmetric derivations}. Thus $D\in\Dera(\g)$ means that
\[
D[x,w]=[Dx,w]+[x,Dw],\qquad
B(Dx,w)+B(x,Dw)=0
\]
hold for all $x,w\in\g$.  
Choose linear maps and an alternating form
\[
\psi:\bLie\longrightarrow\Dera(\g),\qquad
\alpha:\Lambda^2\bLie\longrightarrow\g,\qquad
\Gamma\in\Lambda^3\bLie^*.
\]
The map $\psi$ is not assumed to be a homomorphism.  We set
$\tilde{\mathfrak g}=\bLie\oplus\g\oplus I$, where $I=\bLie^*$, and define
\begin{equation}\label{eq:double-form}
	T(y+x+f,z+w+k)=B(x,w)+f(z)+k(y).
\end{equation}
At this stage, the direct sum is a sum of vector spaces. For $X=y+x+f$, we call $y$, $x$, and $f$ its $\bLie$-, $\mathfrak{g}$-, and $I$-components, respectively. The projections are linear maps and need not be Lie algebra homomorphisms. The subspace $I$ is isotropic for $T$, and $I^\perp=\g\oplus I$. Define
\begin{equation}\label{eq:phi-def}
	\phi(x,w)(y)=B(\psi_yx,w).
\end{equation}
Since $\psi_y$ is skew, $\phi$ is alternating in $x,w$. The derivation property also makes $\phi$ a two-cocycle of $\g$ with values in $I$, where $\g$ acts trivially. The cocycle condition is $\sum_{\mathrm{cyc}}\phi([x,w],u)=0$. It ensures that adding the central term $\phi$ to the bracket of $\g$ preserves the Jacobi identity. We use the coadjoint convention $(\ad_y^*f)(z)=-f([y,z]_{\bLie})$ for $y,z\in\bLie$ and $f\in\bLie^*$. We can now define the generalized bracket on the summands by
\begin{equation}\label{eq:gd-bracket}
	\begin{aligned}
		[y,z]&=[y,z]_{\bLie}+\alpha(y,z)+\Gamma(y,z,\cdot),\\
		[y,x]&=\psi_yx-B(\alpha(y,\cdot),x),\\
		[x,w]&=[x,w]_{\g}+\phi(x,w),\\
		[y,f]&=\ad_y^*f,\qquad [x,f]=0,\qquad [f,k]=0,
	\end{aligned}
\end{equation}
where, in formulas for the data, a bracket without a subscript on two elements of $\bLie$ or two elements of $\g$ denotes the bracket of that algebra. Bilinearity and skew-symmetry determine all other brackets.

Recall that an extension of Lie algebras is an exact sequence $0\to\mathfrak n\to\tilde{\mathfrak g}\to\bLie\to0$, in which $\mathfrak n$ is an ideal and $\tilde{\mathfrak g}/\mathfrak n\cong\bLie$. When the compatibility conditions below hold, our construction gives two successive extensions
\[
0\longrightarrow I\longrightarrow\g\oplus I=I^{\perp}
\longrightarrow\g\longrightarrow0,\qquad
0\longrightarrow I^{\perp}\longrightarrow\tilde{\mathfrak g}
\longrightarrow\bLie\longrightarrow0.
\]
The first is central, while the second relates the bracket on $\g\oplus I$ to the quadratic form and the action of $\bLie$.

\needspace{5\baselineskip}
The following theorem gives the compatibility conditions. It is the restriction of \cite{BajoBenayadiBordemann2007}, Theorem~2.3, to vectors of even degree. Their action map corresponds to $\psi$, their curvature map to $\alpha$, and their covector term to $\Gamma(y,z,\cdot)$.

\begin{theorem}
	\label{thm:gd-bracket}
	The operations \eqref{eq:gd-bracket} and \eqref{eq:double-form} define
	a quadratic Lie algebra if and only if
	\begin{align}
		[\psi_y,\psi_z]-\psi_{[y,z]}&=\ad_{\alpha(y,z)},
		\label{eq:gd-curvature}\\
		\sum_{\mathrm{cyc}(y,z,w)}
		\bigl(\psi_y\alpha(z,w)-\alpha([y,z],w)\bigr)&=0,
		\label{eq:gd-bianchi}\\
		d\Gamma(y,z,w,t)
		&=B(\alpha(y,z),\alpha(w,t))
		-B(\alpha(y,w),\alpha(z,t))\notag\\
		&\hspace{8mm}+B(\alpha(y,t),\alpha(z,w)).
		\label{eq:gd-gamma}
	\end{align}
	Here $d$ stands for the Chevalley--Eilenberg differential\footnote{ For $\omega\in\Lambda^k\bLie^*$,
		our convention is
		\begin{equation}\label{eq:gd-scalar-differential}
			\begin{aligned}
				(d\omega)(y_0,\ldots,y_k)
				=\sum_{0\leq i<j\leq k}(-1)^{i+j}
				\omega\bigl([y_i,y_j],y_0,\ldots,\widehat y_i,
				\ldots,\widehat y_j,\ldots,y_k\bigr),
			\end{aligned}
		\end{equation}
		where a hat means that the argument is omitted and all other
		arguments keep their order.} of $\bLie$ with
	trivial scalar coefficients.
\end{theorem}

\begin{proof}
The form $T$ is symmetric and nondegenerate by its definition. We first check invariance, and then reduce the Jacobi identity to the three stated conditions. For $X_i=y_i+x_i+f_i$, the bracket gives
	\begin{equation}\label{eq:gd-cartan}
		\begin{aligned}
			T([X_1,X_2],X_3)
			={}&\Gamma(y_1,y_2,y_3)+B([x_1,x_2],x_3)\\
			&+\sum_{\mathrm{cyc}(1,2,3)}
			\bigl(B(\psi_{y_1}x_2,x_3)
			+B(\alpha(y_1,y_2),x_3)
			+f_1([y_2,y_3])\bigr).
		\end{aligned}
	\end{equation}
Invariance of $B$, skewness of $\psi_y$, and alternation of $\alpha$ and $\Gamma$ make every term alternating, with each cyclic sum considered as a whole. Thus $T$ is invariant without using the three compatibility conditions.
	
To check the Jacobi identity, define $J(X,Y,Z)=[[X,Y],Z]+[[Y,Z],X]+[[Z,X],Y]$. By invariance,
	\begin{equation}\label{eq:gd-jacobi-fourform}
		\begin{aligned}
			F(X,Y,Z,W):=T(J(X,Y,Z),W)
			={}&T([X,Y],[Z,W])+T([Y,Z],[X,W])\\
			&+T([Z,X],[Y,W]).
		\end{aligned}
	\end{equation}
Symmetry of $T$ and skew-symmetry of the bracket imply that $F$ is alternating. Since $T$ is nondegenerate, the Jacobi identity is equivalent to $F=0$. We check this on arguments from the three summands. If $\pi_{\bLie}$ denotes projection onto $\bLie$, then $\pi_{\bLie}[X,Y]=[\pi_{\bLie}X,\pi_{\bLie}Y]_{\bLie}$, so
	\[
	F(X,Y,Z,f)
	=f\bigl(J_{\bLie}(\pi_{\bLie}X,\pi_{\bLie}Y,
	\pi_{\bLie}Z)\bigr)=0
	\qquad(f\in I).
	\]
This proves the assertion whenever an argument belongs to $I$. For $x,u,v,r\in\g$ and $y\in\bLie$, we also have
	\[
	\begin{aligned}
		F(x,u,v,r)&=B(J_{\g}(x,u,v),r)=0,\\
		F(y,x,u,v)
		&=B\bigl([\psi_yx,u]+[x,\psi_yu]-\psi_y[x,u],v\bigr)=0,
	\end{aligned}
	\]
by the Jacobi identity in $\g$ and the derivation property of $\psi_y$. For the remaining choices, substitution in \eqref{eq:gd-jacobi-fourform} gives
	\[
	\begin{aligned}
		F(y,z,x,u)
		&=B\bigl((\psi_{[y,z]}+\ad_{\alpha(y,z)}
		-[\psi_y,\psi_z])x,u\bigr),\\
		F(y,z,w,x)
		&=-B\Bigl(\sum_{\mathrm{cyc}(y,z,w)}
		\bigl(\psi_y\alpha(z,w)-\alpha([y,z],w)\bigr),x\Bigr),\\
		F(y,z,w,t)
		&=-d\Gamma(y,z,w,t)+B(\alpha(y,z),\alpha(w,t))\\
		&\qquad-B(\alpha(y,w),\alpha(z,t))
		+B(\alpha(y,t),\alpha(z,w)).
	\end{aligned}
	\]
Here $y,z,w,t\in\bLie$ and $x,u\in\g$. By \eqref{eq:gd-scalar-differential}, the six terms involving $\Gamma$ in the last equality sum to $-d\Gamma$. Nondegeneracy of $B$ shows that these expressions vanish for all arguments exactly when \eqref{eq:gd-curvature}--\eqref{eq:gd-gamma} hold. Multilinearity and alternation then give the equivalence with $F=0$.
\end{proof}

We use the following terminology for this construction.

\begin{definition}
	\normalfont 
The quadratic Lie algebra $(\tilde{\mathfrak g},T)$ defined by \eqref{eq:gd-bracket} and \eqref{eq:double-form}, with data satisfying Theorem~\ref{thm:gd-bracket}, is called a \emph{quadratic generalized double extension} with data $(\psi,\alpha,\Gamma)$.
\end{definition}

The first two compatibility conditions have an interpretation in terms of curvature.

\begin{remark} \label{Rmk1_NewInsights}
	\normalfont 
The names \emph{curvature equation} for \eqref{eq:gd-curvature} and \emph{Bianchi equation} for \eqref{eq:gd-bianchi} can be explained using covariant cochains. Consider the extension of Lie algebras of invariant vector fields associated with a principal bundle. A section that is linear over functions is a horizontal lift, and its failure to preserve the bracket is measured by the curvature. See \cite{AlekseevskyMichorRuppert2004}, Sections~2--4, especially equations~(2.2)--(2.4). For $\beta\in\Hom(\Lambda^p\bLie,\g)$, define the \emph{covariant Chevalley--Eilenberg operator} associated with $\psi$ by
	\[
	\begin{aligned}
		(d_\psi\beta)(y_0,\ldots,y_p)
		={}&\sum_{i=0}^p(-1)^i\psi_{y_i}
		\beta(y_0,\ldots,\widehat y_i,\ldots,y_p)\\
		&+\sum_{i<j}(-1)^{i+j}
		\beta([y_i,y_j],y_0,\ldots,\widehat y_i,\ldots,
		\widehat y_j,\ldots,y_p).
	\end{aligned}
	\]
This operator is defined even when $\psi$ is not a representation. Its formula has the same form as that of the exterior covariant derivative $d_{\nabla}$ of a connection $\nabla$ (see \cite{Lee2009}). The first two conditions of Theorem~\ref{thm:gd-bracket} become
	\[
	R_\psi(y,z):=[\psi_y,\psi_z]-\psi_{[y,z]}
	=\ad_{\alpha(y,z)},
	\qquad d_\psi\alpha=0.
	\]
	Thus $\alpha$ is a lift of the curvature with values in $\g$, and the
	second condition is a \emph{lifted Bianchi equation}.  The universal
	Bianchi identity, together with the curvature equation, gives only $d_\psi\alpha\in\Hom\!\bigl(\Lambda^3\bLie,Z(\g)\bigr)$, where $Z(\g)$ is the center of $\g$.  Requiring $d_\psi\alpha=0$
	is therefore stronger in general, see
	\cite{Neeb2006}, Proposition~I.8 and Lemma~I.9(3). The compatibility condition \eqref{eq:gd-gamma} on scalars is the
	$I$-component of the lifted Bianchi equation for the full extension
	by $\g\oplus I$.  This formulation is contained in
	Bajo--Benayadi--Bordemann
	(see \cite{BajoBenayadiBordemann2007}, equations~(22)--(23) in the proof of Theorem~2.3), specialized to the case in which all vectors have even degree.
\end{remark}

\subsection{Ordinary double extensions and
	\texorpdfstring{$T^*_{\theta}$}{T*-theta}-extensions}
\label{subsec:gd-special-cases}
Ordinary double extensions and $T^*_{\theta}$-extensions precede the generalized construction. Their inclusion as special cases is established for superalgebras in \cite{BajoBenayadiBordemann2007}, Theorem~2.4 and Proposition~2.5 in Section~2.3. We describe these inclusions for ordinary Lie algebras. In the ordinary double extension, we also explain how to normalize the invariant bilinear form $q$ on $\bLie$, which may be degenerate.
The map $\phi$ is not independent of $\psi$.  In fact, $\phi=0$ if and only of $\psi=0$, since $B(\psi_yx,w)=0$ for every $w$ implies $\psi_yx=0$ by nondegeneracy.

\subsubsection{Ordinary double extensions}

Favre and Santharoubane construct extensions by a derivation skew with respect to the quadratic form. Together with orthogonal decompositions, these extensions describe quadratic Lie algebras with nonzero center. They also compare the extensions by isometries that preserve the central line (see \cite{FavreSantharoubane1987}, Propositions~2.3, 2.9 and~2.11). Here the isotropic ideal $I$ may have arbitrary dimension and need not be central. We reconstruct the algebra from $I^\perp/I$ by a generalized double extension and then reconstruct Poisson and Jacobi structures whose operators preserve $I$, as described in Sections~\ref{sec:gd-reconstruction} and~\ref{sec:jacobi-double}.

The ordinary double extension of Medina--Revoy starts with a homomorphism $\psi:\bLie\to\Dera(\g)$ and a symmetric invariant bilinear form $q$ on $\bLie$. The form $q$ need not be nondegenerate (see \cite{MedinaRevoy1985}, Section~0.2). We keep this form so that the inclusion below covers the full construction.

\begin{proposition}
	\label{thm:double-quadratic}\label{rem:double-structure}
Let $\psi:\bLie\to\Dera(\g)$ be a homomorphism, and let $q$ be a symmetric invariant bilinear form on $\bLie$. The ordinary double extension has bracket
	\begin{equation}\label{eq:double-bracket}
		\begin{gathered}
			[y,z]_{\mathrm{ord}}=[y,z]_{\bLie},\qquad
			[y,x]_{\mathrm{ord}}=\psi_yx,\qquad
			[y,f]_{\mathrm{ord}}=\ad_y^*f,\\
			[x,w]_{\mathrm{ord}}=[x,w]_{\g}+\phi(x,w),\qquad
			[x,f]_{\mathrm{ord}}=[f,k]_{\mathrm{ord}}=0
		\end{gathered}
	\end{equation}
	and form
	\begin{equation}\label{eq:gd-ordinary-extra-form}
		T_q(y+x+f,z+w+k)=q(y,z)+B(x,w)+f(z)+k(y).
	\end{equation}
	For $q=0$, these are exactly \eqref{eq:gd-bracket} and
	\eqref{eq:double-form} with $\alpha=\Gamma=0$.
	For arbitrary $q$, the quadratic Lie algebra is isometrically
	isomorphic to the generalized double extension with form $T$ and data
	\begin{equation}\label{eq:gd-ordinary-canonical-data}
		\psi,\qquad \alpha=0,\qquad
		\Gamma_q(y,z,w)=-\tfrac12q([y,z],w).
	\end{equation}
	More precisely, define $q^\flat y=q(y,\cdot)$.  The isomorphism from
	these generalized data to the ordinary extension is
	\begin{equation}\label{eq:gd-ordinary-isometry}
		F_q(y+x+f)=y+x+f-\tfrac12q^\flat y.
	\end{equation}
	Consequently every ordinary double extension, including every choice
	of $q$, occurs in the generalized construction up to isometry.
\end{proposition}

\begin{proof}
For $\alpha=\Gamma=0$, the conditions in Theorem~\ref{thm:gd-bracket} reduce to the homomorphism property of $\psi$. This proves the assertion for $q=0$. For general $q$, the difference between $T_q$ and $T$ is the pullback of the invariant form $q$ along the quotient map onto $\bLie$, so $T_q$ is invariant. The map $F_q$ is invertible, and symmetry of $q$ gives
	\[
	T_q(F_qX,F_qY)
	=T(X,Y)+q(y,z)-\tfrac12q(y,z)-\tfrac12q(z,y)=T(X,Y).
	\]
	This also proves nondegeneracy of $T_q$.  Invariance of $q$ gives
	$\ad_y^*q^\flat z=q^\flat[y,z]$.  Thus
	\[
	\begin{aligned}
		[F_qy,F_qz]_{\mathrm{ord}}&=[y,z]-q^\flat[y,z],\\
		F_q^{-1}[F_qy,F_qz]_{\mathrm{ord}}
		&=[y,z]-\tfrac12q^\flat[y,z].
	\end{aligned}
	\]
	All other bracket types are unchanged, since $F_q$ fixes $\g$
	and $I$.  Thus transport by $F_q$ gives exactly
	\eqref{eq:gd-bracket} with the data
	\eqref{eq:gd-ordinary-canonical-data}.  The form $\Gamma_q$ is
	alternating by symmetry and invariance of $q$, and its cocycle
	condition follows from Theorem~\ref{thm:gd-bracket}.
\end{proof}

In the presentation \eqref{eq:double-bracket}, the summand $\bLie$
is a Lie subalgebra.  It is isotropic for $T_q$ exactly when $q=0$.
The summand $\g$ is a Lie subalgebra exactly when $\phi=0$,
equivalently $\psi=0$.

\subsubsection{\texorpdfstring{$T^*_{\theta}$}{T*-theta}-extensions}

Let $\theta:\Lambda^2\bLie\to\bLie^*$ be a two-cocycle for
the coadjoint representation of $\bLie$.  It is said to be \emph{cyclic} when
\begin{equation}\label{eq:gd-tstar-cyclicity}
	\theta(y,z)(w)=\theta(z,w)(y).
\end{equation}
Bordemann's $T^*_\theta$-extension uses these cocycles with the canonical pairing. Cyclicity and the correspondence with scalar three-cocycles are established in \cite{Bordemann1997}, Lemma~3.1 and Section~3.1.

\begin{proposition}\label{prop:gd-all-tstar}
	Setting $\g=0$ in the generalized double extension gives exactly
	all quadratic $T^*_{\theta}$-extensions with their canonical pairing.
	The bracket and form are
	\begin{equation}\label{eq:gd-special-tstar}
		\begin{aligned}
			[y+f,z+k]_{\theta}
			&=[y,z]_{\bLie}+\theta(y,z)+\ad_y^*k-\ad_z^*f,\\
			T(y+f,z+k)&=f(z)+k(y).
		\end{aligned}
	\end{equation}
	The correspondence between their data is
	\begin{equation}\label{eq:gd-tstar-correspondence}
		\theta(y,z)(w)=\Gamma(y,z,w),\qquad
		d_{\ad^*}\theta=0\quad\Longleftrightarrow\quad d\Gamma=0,
	\end{equation}
	where cyclicity of $\theta$ corresponds to alternation of $\Gamma$.
	The ideal $I=\bLie^*$ is Lagrangian.
\end{proposition}

\begin{proof}
When $\g=0$, we have $\psi=\alpha=\phi=0$. The bracket and form reduce to \eqref{eq:gd-special-tstar}, and the only remaining compatibility condition is $d\Gamma=0$. An alternating map $\theta:\Lambda^2\bLie\to\bLie^*$ is cyclic exactly when $\Gamma(y,z,w):=\theta(y,z)(w)$ is alternating. To compare the cocycle conditions, evaluate the coadjoint differential:
	\[
	\begin{aligned}
		(d_{\ad^*}\theta)(y,z,w)(t)
		={}&-\Gamma(z,w,[y,t])+\Gamma(y,w,[z,t])
		-\Gamma(y,z,[w,t])\\
		&-\Gamma([y,z],w,t)+\Gamma([y,w],z,t)
		-\Gamma([z,w],y,t)\\
		={}&d\Gamma(y,z,w,t).
	\end{aligned}
	\]
	This proves both directions of the cocycle correspondence and
	thus the inclusion of all quadratic $T^*_{\theta}$-extensions.
	Finally, the canonical pairing gives $I^\perp=I$.
\end{proof}

The constructions by double extension and $T^*_{\theta}$-extension overlap at
\begin{equation}\label{eq:gd-special-tstar-zero}
	T^*_0\bLie=\bLie\ltimes_{\ad^*}\bLie^*,
\end{equation}
which is the ordinary double extension with $\g=0$ and
$q=0$. 
\subsection{Adapted Witt splitting and recognition from an isotropic ideal}
\label{subsec:gd-recognition}

We now reverse the construction. Starting with a quadratic Lie algebra $(\mathfrak d,T)$ and a marked isotropic ideal $I$, we recover the data of a generalized double extension. We also determine when the algebra admits a presentation as an ordinary double extension, including the additional invariant form $q$ in \eqref{eq:gd-ordinary-extra-form}.

The reconstruction from an arbitrary isotropic ideal is the version for ordinary Lie algebras of \cite{BajoBenayadiBordemann2007}, Theorem~3.1, which is stated over an algebraically closed field of characteristic zero. Our proof constructs the splitting by the correction $s=s_0-\tfrac12\ell\beta^\flat$ and therefore requires only $\operatorname{char}\K\ne2$. For the ordinary reconstruction, we use the classical argument based on a complement that is a Lie subalgebra (see \cite{MedinaRevoy1985}, Sections~2.1--2.2, and \cite{FigueroaStanciu1996}, Section~3). We apply it to the specified ideal $I$, without assuming minimality, and keep the possibly degenerate form $q$.

An \emph{adapted Witt splitting} of the flag
$I\subseteq I^\perp\subseteq\mathfrak d$ consists of linear sections
$j:\g:=I^\perp/I\to I^\perp$ and
$s:\bLie:=\mathfrak d/I^\perp\to\mathfrak d$ such that
$s(\bLie)$ is isotropic and orthogonal to $j(\g)$.
Together with the inclusion of $I$, they identify
$\mathfrak d$ with $\bLie\oplus\g\oplus\bLie^*$ and its form
with \eqref{eq:double-form}. The sections need not be Lie algebra homomorphisms. The existence of such a splitting is a relative form of Witt decomposition for the prescribed totally isotropic subspace $I$. See \cite{Witt1937} and \cite{ElmanKarpenkoMerkurjev2008}, Lemma 8.10.

\needspace{10\baselineskip}
\begin{theorem}
	\label{thm:gd-recognition}
	Let $(\mathfrak d,T)$ be a quadratic Lie algebra and let $I$ be an isotropic ideal.
	Then
	\begin{equation}\label{eq:gd-isotropic-flag}
		[I,I]=0,\qquad [I,I^\perp]=0.
	\end{equation}
	Define $\bLie:=\mathfrak d/I^\perp$ and $\g :=I^\perp/I.$
	The quotient $\g$ carries the quadratic form induced by $T$, and
	$T$ identifies $I$ with $\bLie^*$.  Adapted Witt splittings exist,
	and every such splitting gives an isometric presentation
	\eqref{eq:gd-bracket} whose data satisfy
	\eqref{eq:gd-curvature}--\eqref{eq:gd-gamma}.
	
A presentation as an ordinary double extension relative to $I$, allowing a symmetric invariant bilinear form $q$ as in \eqref{eq:gd-ordinary-extra-form}, exists if and only if the quotient
	map $\pi:\mathfrak d\to\bLie$ admits a Lie algebra section
	$t:\bLie\to\mathfrak d$.  Equivalently, there is a Lie subalgebra
	$L\subseteq\mathfrak d$ complementary to $I^\perp$.  For the presentation associated with $t$, the
	additional form is
	\[
	q(y,z)=T(t(y),t(z)),\qquad y,z\in\bLie,
	\]
	it is automatically symmetric and invariant, but it does not need to be
	nondegenerate.
	
	Equivalently, an adapted Witt splitting can be chosen for which
	\[
	\alpha=0,\qquad
	\Gamma(y,z,w)=\Gamma_q(y,z,w)
	=-\tfrac12q([y,z]_{\bLie},w)
	\]
for some symmetric invariant bilinear form $q$ on $\bLie$.
	A presentation with $q=0$ exists if and only if the complement $L$
	can be chosen isotropic, or equivalently if an adapted Witt splitting
	can be chosen with $\alpha=\Gamma=0$.
\end{theorem}

\begin{proof}
Let $p:I^\perp\to\g=I^\perp/I$ be the quotient map.
We first describe the quotients determined by $I$, and then construct the splitting. Invariance of $T$ gives $[I,I]=0$, $[I^{\perp},I]=0$, and shows that $I^{\perp}$ is an ideal. The restriction of $T$ to $I^\perp$ has radical $I$, so it induces a nondegenerate invariant form $B$ on $\g$. The pairing between $I$ and $\mathfrak d/I^\perp$ is perfect. The inverse identification $\ell:\bLie^*\to I$ is characterized by
	\[
	T(\ell(f),X)=f(\pi X).
	\]
	
Choose a linear section $j:\g\to I^\perp$. Since
	$T(jx,ju)=B(x,u)$, the subspace $j(\g)$ is nondegenerate.
	Its orthogonal complement $E$ satisfies
	$E\cap I^\perp=I$, and $\pi|_E:E\to\bLie$ is onto.
Choose a linear section $s_0:\bLie\to E$. To make its image isotropic, define
	\[
	\beta(y,z)=T(s_0y,s_0z),\qquad
	\beta^\flat y=\beta(y,\cdot),\qquad
	s(y)=s_0(y)-\tfrac12\ell(\beta^\flat y).
	\]
	Because $I$ is isotropic and $T(\ell(f),s_0y)=f(y)$,
	\[
	T(sy,sz)=\beta(y,z)-\tfrac12\beta(y,z)
	-\tfrac12\beta(z,y)=0.
	\]
Moreover, $s(\bLie)\subset E$. Thus $s,j$ form an adapted Witt splitting, and $\Phi(y+x+f)=s(y)+j(x)+\ell(f)$ is an isometry for \eqref{eq:double-form}.
	
We now recover the extension data. For any adapted splitting, define
	\[
	\begin{aligned}
		\psi_yx&=p([s(y),j(x)]),\\
		\alpha(y,z)&=p([s(y),s(z)]-s([y,z]_{\bLie})),\\
		\Gamma(y,z,w)&=T([s(y),s(z)],s(w)).
	\end{aligned}
	\]
	The adjoint action of $s(y)$ induces a derivation on $I^\perp/I$,
	and invariance of $T$ makes it $B$-skew.  Thus
	$\psi_y\in\Dera(\g)$, while $\alpha$ and $\Gamma$ are alternating.
	The quotient maps give the $\bLie$- and $\g$-components of the
	bracket.  Its $I$-components are determined by pairing with
	$s(\bLie)$, using
	\[
	\begin{aligned}
		T([s(y),j(x)],s(z))&=-B(\alpha(y,z),x),\\
		T([j(x),j(u)],s(y))&=B(\psi_yx,u),\\
		[s(y),\ell(f)]&=\ell(\ad_y^*f).
	\end{aligned}
	\]
	Together with the definition of $\Gamma$ and $[I,I^\perp]=0$,
	these identities give exactly \eqref{eq:gd-bracket} under $\Phi$.
	The compatibility conditions follow from
	Theorem~\ref{thm:gd-bracket}.
	
It remains to characterize ordinary presentations. Such a presentation gives a Lie algebra section through its $\bLie$-summand. Conversely, suppose a Lie algebra section $t$ exists and define $L=t(\bLie)$. The pairing between $I$ and $L$
	shows that $W=I\oplus L$ is nondegenerate, even if
	$q(y,z):=T(t(y),t(z))$ is degenerate.  The complement
	$H=W^\perp$ lies in $I^\perp$, and $p|_H:H\to\g$ is an
	isomorphism.  Let $j$ be its inverse.  Since $W$ is a subalgebra,
	invariance gives $[L,H]\subset H$.  Therefore
	\[
	[t(y),j(x)]=j(\psi_yx),\qquad
	\psi:\bLie\longrightarrow\Dera(\g)
	\]
	is a homomorphism, by Jacobi.  The form $q$ is symmetric and
	invariant because $t$ is a Lie algebra homomorphism.  Pairing with
	$t(\bLie)$ as above now gives
	\[
	[j(x),j(u)]=j([x,u]_{\g})+\ell(\phi(x,u)),\qquad
	[t(y),\ell(f)]=\ell(\ad_y^*f).
	\]
	Consequently $\Phi_q(y+x+f)=t(y)+j(x)+\ell(f)$ is a Lie algebra
	isomorphism for \eqref{eq:double-bracket}, carries $\bLie^*$ onto
	$I$, and pulls $T$ back to $T_q$.  This proves the criterion for a presentation as an ordinary double extension.  A Lie algebra section and a Lie
	subalgebra complement to $I^\perp$ are equivalent by restriction
	of $\pi$.
	
	To obtain adapted coordinates from this ordinary presentation,
	define $q^\flat y=q(y,\cdot)$ and subtract half of this contribution:
	\[
	s(y)=t(y)-\tfrac12\ell(q^\flat y).
	\]
	As before, $s(\bLie)$ is isotropic and orthogonal to $j(\g)$.
	Using $\ad_y^*q^\flat z=q^\flat[y,z]_{\bLie}$, we obtain
	\[
	\begin{aligned}
		[s(y),s(z)]
		&=t([y,z]_{\bLie})-\ell(q^\flat[y,z]_{\bLie})=s([y,z]_{\bLie})-\tfrac12\ell(q^\flat[y,z]_{\bLie}),\\
		[s(y),j(x)]&=j(\psi_yx).
	\end{aligned}
	\]
	Thus $\alpha=0$ and $\Gamma=\Gamma_q$.  Conversely, these data
	give an ordinary presentation by Proposition~\ref{thm:double-quadratic},
	with Lie algebra section $t(y)=s(y)+\tfrac12\ell(q^\flat y)$.
	Finally, $q=0$ precisely when the associated complement $L$ is
	isotropic.  In adapted coordinates, $s(\bLie)$ is a Lie subalgebra
	precisely when $\alpha=\Gamma=0$, proving the last equivalences.
\end{proof}

From now on, an adapted Witt splitting of a quadratic Lie algebra $(\mathfrak d,T)$ with a marked isotropic ideal $I$ identifies it with the generalized double extension $\widetilde{\g}:=\bLie\oplus\g\oplus\bLie^*$. We use the form \eqref{eq:double-form}, the bracket \eqref{eq:gd-bracket}, and the data $\psi$, $\alpha$, and $\Gamma$ determined by that splitting. In particular, we identify $I$ with $\bLie^*$ and $\Lambda^k I$ with $\Lambda^k\bLie^*$.

We record a consequence for nilpotent and solvable algebras. It follows from the preceding reconstruction and the results of \cite{Bordemann1997}.

\begin{corollary}\label{cor:reconstruction_nilpotent_solvable}
	Let $(\mathfrak d,T)$ be a quadratic Lie algebra of even dimension such that either $\mathfrak d$ is nilpotent and $\operatorname{char}\K\ne2$, or $\mathfrak d$ is solvable and $\operatorname{char}\K=0$. Then the quadratic Lie algebra structure on $\mathfrak d$ can be obtained as a generalized double extension.
\end{corollary}

If $(\mathfrak d,T)$ is nonzero, indecomposable, and neither simple nor one-dimensional, then every minimal nonzero proper ideal is isotropic. Indeed, the radical of the restricted form is an ideal and is therefore either zero or the whole minimal ideal. The first alternative would give a nontrivial orthogonal decomposition of $\mathfrak d$. Thus Theorem~\ref{thm:gd-recognition} applies to such an ideal; see \cite{MedinaRevoy1985}.

\subsection{Changes of adapted splitting}
\label{subsec:gd-splittings}

The reconstruction in Subsection~\ref{subsec:gd-recognition} depends on an adapted Witt splitting. Different choices may give different data $\psi,\alpha,\Gamma$ for the same quadratic Lie algebra and the same ideal $I$. To compare these presentations, we describe how the data change with the splitting. We also compute the effect on operators and vectors, since these formulas will be used for the Poisson and Jacobi structures in Sections~\ref{sec:gd-reconstruction} and~\ref{sec:jacobi-double}. A change of splitting may give $\alpha'\neq0$ or $\Gamma'\neq0$ even when $\alpha=\Gamma=0$. We must therefore distinguish the vanishing of these data in a particular splitting from the existence of an ordinary double-extension presentation relative to $I$.

Keep the identification $I=\bLie^*$ and fix the identifications of $I$,
$I^\perp/I$, and $\widetilde{\g}/I^\perp$.  The isometries that we use have
the same matrix as those in \cite{KathOlbrich2006}, Section 3, with the
summands in reverse order, but the description here does not assume that $\g$
is abelian.  For $F\in\Hom(\bLie,\g)$ and
$\omega\in\Lambda^2\bLie^*$, define $(F^\dagger x)(y)=B(Fy,x)$ and $(\omega^\flat y)(z)=\omega(y,z)$.  Except in the coadjoint notation $\ad_y^*$ fixed above, the superscript $*$ on an endomorphism of $\bLie$ denotes its ordinary dual. For $\tau:\bLie\to\g$ and $G\in\Lambda^2\bLie^*$, define
\begin{equation}\label{eq:gd-splitting-map}
	T_{\tau,G}(y+x+f)
	=y+(x+\tau y)
	+f-\tau^\dagger x+G^\flat y-\tfrac12\tau^\dagger\tau y.
\end{equation}
We will also need the effect of a change of splitting on skew operators that preserve $I$. Let $\mathcal S\in\End(\widetilde{\g})$ satisfy $T(\mathcal SX,Y)+T(X,\mathcal SY)=0$ and $\mathcal S(I)\subseteq I$. Relative to $\bLie\oplus\g\oplus\bLie^*$, it has the form
\[
\mathcal S=
\begin{pmatrix}
	A&0&0\\
	V&C&0\\
	\eta^\flat&-V^\dagger&-A^*
\end{pmatrix},
\]
where $A\in\End(\bLie)$, $V\in\Hom(\bLie,\g)$,
$C\in\mathfrak o(\g,B)$, and $\eta\in\Lambda^2\bLie^*$.
This is the case of Lemma~\ref{lem:gd-skew} in which the ideal $I$
is preserved.  The operators induced on $\widetilde{\g}/I^\perp$ and
$I^\perp/I$ are $A$ and $C$, respectively.  Thus $A$ and $C$ do not
depend on the splitting, whereas $V$ and $\eta$ may change.
For a vector $E\in\widetilde{\g}$, write
\[
E=y_E+x_E+f_E,\qquad
y_E\in\bLie,\quad x_E\in\g,\quad f_E\in\bLie^*.
\]
When $(\mathcal S,E)$ is a Jacobi structure, $E$ is the Euler element
in the notation of Section~\ref{sec:jacobi-double}. Primes will denote the data obtained by transport through
$T_{\tau,G}^{-1}$.  More precisely, we define
$[X,Y]'=T_{\tau,G}^{-1}[T_{\tau,G}X,T_{\tau,G}Y]$,
$\mathcal S'=T_{\tau,G}^{-1}\mathcal S T_{\tau,G}$, and
$E'=T_{\tau,G}^{-1}E$.

\needspace{5\baselineskip}
The next proposition characterizes the changes of splitting and gives their effect on the operator and the vector.

\begin{proposition}
	\label{prop:gd-splitting}
	The maps \eqref{eq:gd-splitting-map} are exactly the isometries that fix
	$\bLie^*$ pointwise and induce the identity on
	$I^\perp/I$ and $\widetilde{\g}/I^\perp$.
	For $\mathcal S$ as above, its blocks transform according to
	\begin{equation}\label{eq:gd-splitting-ACV}
		A'=A,\qquad C'=C,\qquad V'=V+C\tau-\tau A,
	\end{equation}
	\begin{equation}\label{eq:gd-splitting-eta}
		\begin{aligned}
			\eta'^\flat={}&\eta^\flat-V^\dagger\tau+\tau^\dagger V
			+\tau^\dagger C\tau-A^*G^\flat-G^\flat A\\
			&+\tfrac12A^*\tau^\dagger\tau
			-\tfrac12\tau^\dagger\tau A,
		\end{aligned}
	\end{equation}
	and the components of $E$ satisfy
	\begin{equation}\label{eq:gd-splitting-E}
		y_E'=y_E,\qquad x_E'=x_E-\tau y_E,\qquad
		f_E'=f_E+\tau^\dagger x_E-G^\flat y_E-\tfrac12\tau^\dagger\tau y_E.
	\end{equation}
	In particular, the conditions $E\in I$,
	$E\in I^\perp\setminus I$, and $E\notin I^\perp$ do not depend on the
	splitting.  
\end{proposition}

\begin{proof}
	Let $P$ be an isometry that fixes $\bLie^*$ pointwise and induces
	the identity on $I^\perp/I$ and $\widetilde{\g}/I^\perp$.
	Since $P$ preserves $I^\perp$, its matrix has the form
	\[
	P=\begin{pmatrix}1&0&0\\ \tau&1&0\\ R&S&1\end{pmatrix}.
	\]
	Pairing $Py$ with $Px$ and with $Pz$, for $y,z\in\bLie$ and
	$x\in\g$, gives
	\[
	S=-\tau^\dagger,\qquad
	(Ry)(z)+(Rz)(y)=-B(\tau y,\tau z).
	\]
	Thus $R=G^\flat-\tfrac12\tau^\dagger\tau$ for a unique
	$G\in\Lambda^2\bLie^*$, so $P=T_{\tau,G}$.  Conversely, these
	identities and \eqref{eq:double-form} show that $T_{\tau,G}$ is an
	isometry with the required action on $\bLie^*$ and the quotients. Multiplication of the matrices gives
	$$T_{\tau,G}T_{\sigma,H}=T_{\tau+\sigma,K},\qquad K^\flat=G^\flat+H^\flat
	+\tfrac12(\sigma^\dagger\tau-\tau^\dagger\sigma).$$
	The form $K$ is alternating by symmetry of $B$.  In particular,
	the identity is $T_{0,0}$ and
	\[
	T_{\tau,G}^{-1}
	=\begin{pmatrix}
		1&0&0\\ -\tau&1&0\\
		-G^\flat-\tfrac12\tau^\dagger\tau&\tau^\dagger&1
	\end{pmatrix}=T_{-\tau,-G}.
	\]
Thus changes of splitting form a group. The composition rule for transporting the bracket, operator, and vector follows from their definitions.
	
	For $\mathcal S'=P^{-1}\mathcal SP$, multiplication gives
	$A'=A$, $C'=C$, and $V'=V+C\tau-\tau A$, proving
	\eqref{eq:gd-splitting-ACV}.  Define
	$\widetilde R=-G^\flat-\tfrac12\tau^\dagger\tau$.
	The block from $\bLie$ to $\bLie^*$ is
	\[
	\widetilde R A+\tau^\dagger(V+C\tau)
	+\eta^\flat-V^\dagger\tau-A^*R.
	\]
	Substituting $R$ and $\widetilde R$ gives
	\eqref{eq:gd-splitting-eta}.  The other blocks follow from preservation of $I$ and skewness
	relative to $T$, both of which hold for $P^{-1}\mathcal SP$.
	In particular, the block from $\g$ to $\bLie^*$ is
	$-(V')^\dagger$ and the restriction to $\bLie^*$ is $-A^*$.	Finally, \eqref{eq:gd-splitting-E} follows by applying the inverse above to $E$.  Since $P$ preserves $I$ and $I^\perp$, the three
	conditions on $E$ are independent of the splitting.
\end{proof}

Before computing the transformation of the bracket, let us recall the
operator $d_\psi$ introduced in Remark~\ref{Rmk1_NewInsights}.
On cochains of degree one it is given by
\[
(d_\psi\tau)(y,z)=\psi_y\tau z-\psi_z\tau y-\tau[y,z].
\]
As in Remark~\ref{Rmk1_NewInsights}, this formula is meaningful even when $\psi$ is not a representation. On alternating cochains with values in $\g$, the operator $d_\psi$ is given by the Chevalley--Eilenberg formula with $\psi_y$ in place of the action on coefficients. However, these cochains need not form a complex. For $x\in\g$, regarded as a cochain of degree zero, we have $(d_\psi x)(y)=\psi_yx$, and therefore
$$(d_\psi^2x)(y,z)=\psi_y\psi_zx-\psi_z\psi_yx-\psi_{[y,z]}x=[\alpha(y,z),x],$$
by \eqref{eq:gd-curvature}. The right-hand side need not vanish, so $d_\psi$ need not be a differential. If $\alpha=0$, equation \eqref{eq:gd-curvature} makes $\psi$ a representation and $d_\psi$ its Chevalley--Eilenberg differential, regardless of $\Gamma$. The same conclusion holds whenever the curvature vanish, even if $\alpha\neq0$. Also, \eqref{eq:gd-bianchi} reads $d_\psi\alpha=0$. The terms involving $\psi$ are $\psi_y\alpha(z,w)-\psi_z\alpha(y,w)+\psi_w\alpha(y,z)$, while those involving the bracket of $\bLie$ are $-\alpha([y,z],w)+\alpha([y,w],z)-\alpha([z,w],y)$.

\needspace{5\baselineskip}
We can now compute the changes in $\psi,\alpha,\Gamma$. In particular, the formula for $\Gamma$ determines the $\bLie^*$-component of $[y,z]'$ for $y,z\in\bLie$.

\begin{proposition}
	\label{prop:gd-bracket-splitting}
	Under the change of splitting \eqref{eq:gd-splitting-map}, the data
	$\psi,\alpha,\Gamma$ satisfy
	\begin{align}
		\psi'_y&=\psi_y+\ad_{\tau y},
		\label{eq:gd-splitting-psi}\\
		\alpha'(y,z)&=\alpha(y,z)+(d_\psi\tau)(y,z)+[\tau y,\tau z],
		\label{eq:gd-splitting-alpha}\\
		\Gamma'(y,z,w)
		&=\Gamma(y,z,w)+dG(y,z,w)
		+\sum_{\mathrm{cyc}}B(\alpha(y,z),\tau w)\notag\\
		&\quad+\tfrac12\sum_{\mathrm{cyc}}
		B((d_\psi\tau)(y,z),\tau w)
		+B([\tau y,\tau z],\tau w).
		\label{eq:gd-splitting-Gamma}
	\end{align}
	Therefore, the data again define a generalized double extension.
	The condition $\alpha=\Gamma=0$ does not imply
	$\alpha'=\Gamma'=0$ after changing the splitting.  Nevertheless,
	the algebra and its ideal $\bLie^*$ are unchanged, so the existence
	of a presentation as an ordinary double extension is preserved.
\end{proposition}

\begin{proof}
	Define $P=T_{\tau,G}$ and
	$\beta_y=G^\flat y-\tfrac12\tau^\dagger\tau y$.  Then $Py=y+\tau y+\beta_y$, $Px=x-\tau^\dagger x$, and $Pf=f$. By \eqref{eq:gd-bracket}, the $\g$-component of $[Py,Px]$ is
	$(\psi_y+\ad_{\tau y})x$, and its $\bLie$-component vanishes.
	Applying $P^{-1}$ therefore gives \eqref{eq:gd-splitting-psi}.
	Similarly, the $\g$-component of $[Py,Pz]$ is
	\[
	N(y,z)=\alpha(y,z)+\psi_y\tau z-\psi_z\tau y+[\tau y,\tau z],
	\]
	while its $\bLie$-component is $[y,z]$.  Thus the $\g$-component
	of $[y,z]'$ is $N(y,z)-\tau[y,z]$, which proves
	\eqref{eq:gd-splitting-alpha}. The $\bLie^*$-components of $[y,x]'$ and $[x,w]'$ are
	$-B(\alpha'(y,\cdot),x)$ and $\phi'(x,w)$, respectively, where
	$\phi'(x,w)(y)=B(\psi'_yx,w)$.  This follows using
	\eqref{eq:gd-bracket}, skewness of $\psi_y$, and invariance of $B$.
	In particular, the $\bLie^*$-component of $[x,w]'$ is
	$\phi(x,w)+\tau^\dagger[x,w]$, whose value at $y$ is
	$B((\psi_y+\ad_{\tau y})x,w)$.  The brackets involving $\bLie^*$
	are unchanged, since $P$ fixes this summand and
	$[\bLie^*,\g\oplus\bLie^*]=0$.
	
	It remains to compute $\Gamma'$.  Since $P$ is an isometry,
	\[
	\Gamma'(y,z,w)=T(P^{-1}[Py,Pz],w)=T([Py,Pz],Pw).
	\]
	Substitution in \eqref{eq:gd-cartan} gives
	\[
	\begin{aligned}
		\Gamma'(y,z,w)
		={}&\Gamma(y,z,w)+B([\tau y,\tau z],\tau w)\\
		&+\sum_{\mathrm{cyc}}B(\alpha(y,z),\tau w)
		+\sum_{\mathrm{cyc}}B(\psi_y\tau z,\tau w)
		+\sum_{\mathrm{cyc}}\beta_y([z,w]).
	\end{aligned}
	\]
	By \eqref{eq:gd-scalar-differential} and symmetry of $B$,
	\[
	\sum_{\mathrm{cyc}}\beta_y([z,w])
	=dG(y,z,w)-\tfrac12\sum_{\mathrm{cyc}}B(\tau[y,z],\tau w).
	\]
	Moreover, skewness of $\psi_y$ gives $\sum_{\mathrm{cyc}}B(\psi_z\tau y,\tau w)
	=-\sum_{\mathrm{cyc}}B(\psi_y\tau z,\tau w)$. Therefore,
	
	$$\sum_{\mathrm{cyc}}B(\psi_y\tau z,\tau w)
	-\tfrac12\sum_{\mathrm{cyc}}B(\tau[y,z],\tau w)=\tfrac12\sum_{\mathrm{cyc}}
	B((d_\psi\tau)(y,z),\tau w),$$
	which proves \eqref{eq:gd-splitting-Gamma}.
	
	Finally, the Jacobi identity and invariance of $T$ are preserved
	by transport through the isometry $P$.  Thus the data satisfy
	Theorem~\ref{thm:gd-bracket}.  Notice that even when $\alpha=\Gamma=0$, the terms $d_\psi\tau$,
	$[\tau y,\tau z]$, and those in the formula for $\Gamma'$ need
	not vanish.  Thus, the image of the section $y\mapsto Py$ need
	not be a Lie subalgebra, even when the image of $y\mapsto y$ has
	this property.  Nevertheless, $P$ is an isomorphism and fixes $\bLie^*$.
	It preserves isotropy and complements of $I^\perp$.  The characterization in
	Theorem~\ref{thm:gd-recognition} therefore shows that the existence
	of a presentation as an ordinary double extension is independent of
	the splitting.
\end{proof}
\enlargethispage{3\baselineskip}

\section{Poisson structures on quadratic Lie algebras}
\label{sec:poisson-preliminaries}

We now use the quadratic form to express the Poisson condition as an equation for an endomorphism. We then describe the symplectic structure on its image and reconstruct the ambient algebra when the quadratic form is nondegenerate on that image. This approach builds on \cite{BajoBenayadiMedina2007,BoucettaMedina2010,MedinaRevoy1988}. As before, the Lie algebras are finite-dimensional over a field $\K$ of characteristic different from $2$, unless otherwise stated, and statements about Lie groups are made over $\K=\mathbb R$.

Let $(\g, B)$ be a quadratic Lie algebra. An endomorphism $A\in \End(\g)$ is called $B$-\emph{skew} if
\begin{equation}\label{eq:B-skew}
	B(AX,Y)+B(X,AY)=0,
	\qquad X,Y\in\g.
\end{equation}
The form $B$ identifies $\g$ with $\g^*$. We write
\[
        B^\flat\colon\g\longrightarrow\g^*,
        \qquad B^\flat(X)=B(X,\cdot),
        \qquad B^\sharp:=(B^\flat)^{-1}.
\]
We also abbreviate $X^\flat:=B^\flat X$ and $\alpha^\sharp:=B^\sharp\alpha$.
If $A$ is $B$-skew, it defines a bivector $r_A\in\Lambda^2\g$ by
\begin{equation}\label{eq:rA-def}
        r_A\bigl(B(X,\cdot),B(Y,\cdot)\bigr)=B(AX,Y),
        \qquad X,Y\in\g.
\end{equation}
Over $\K=\mathbb R$, let $G$ be a simply connected Lie group integrating $\g$. The left-invariant extension of $r_A$ is
\begin{equation}\label{eq:piA-def}
        \pi_A:=r_A^L,
        \qquad
        \pi_A(g)=(L_g)_*r_A.
\end{equation}
The bivector $\pi_A$ is Poisson exactly when $[r_A,r_A]=0$ in the Gerstenhaber algebra $\Lambda^\bullet\g$. Over any field of characteristic zero, we take this tensor equation as the algebraic Poisson condition. It is the classical Yang--Baxter equation for $r_A$ \cite{Drinfeld1983,Drinfeld1987,KosmannSchwarzbach1997}.

We use both the tensor and operator forms of this equation. Neither requires a metric, but $B$ identifies the $B$-skew maps with skew tensors and relates the two formulations. For the operator form, see \cite{Bordemann1990,Kupershmidt1999, MedinaRevoy1988,SemenovTianShansky1983}.

\begin{definition}\label{def:cybe-two-notions}
	\normalfont 
\begin{enumerate}
\item[\textnormal{(i)}] Let $\h$ be a Lie algebra. A \emph{tensor solution} of the classical Yang--Baxter equation on $\h$, also called a skew-symmetric classical $r$-matrix, is a bivector $r\in\Lambda^{2}\h$ satisfying
\begin{equation}\label{eq:cybe-tensor}
        [r,r]=0\qquad\text{in }\Lambda^{3}\h.
\end{equation}
\item[\textnormal{(ii)}] Let $\h$ be a Lie algebra.  An \emph{operator solution} of the classical
Yang--Baxter equation (CYBE\@ for short) on $\h$ is an endomorphism $A\in\End(\h)$ satisfying
\[
        [AX,AY]=A\bigl([AX,Y]+[X,AY]\bigr),\qquad X,Y\in\h.
\]
Equivalently, $A$ is a Rota--Baxter operator of weight zero on the Lie algebra $\h$.
\end{enumerate}
\end{definition}

Condition~\textnormal{(ii)} implies that $\h_A:=\Imop A$ is a Lie subalgebra of $\h$. For a quadratic Lie algebra $(\g,B)$, the identification \eqref{eq:rA-def} also makes \textnormal{(i)} and \textnormal{(ii)} equivalent. We recall the proof to fix the notation and signs used later. The role of $B$ is to relate the two equations, without appearing explicitly in either of them.

\begin{lemma}\label{lem:operator-cybe-equivalence}
Let $(\g,B)$ be a quadratic Lie algebra and let $A\in\End(\g)$ be $B$-skew.
Then $[r_A,r_A]=0$ if and only if $A$ satisfies the classical
Yang--Baxter equation in operator form
\begin{equation}\label{eq:operator-cybe-lemma}
        [AX,AY]
        =A\bigl([AX,Y]+[X,AY]\bigr),
        \qquad X,Y\in\g.
\end{equation}
When $\K=\mathbb R$, these conditions are equivalent to $\pi_A$ being
Poisson on every simply connected Lie group that integrates $\g$.
\end{lemma}

\begin{proof}
Write $X^{\flat}=B(X,\cdot)$. We first identify the map associated with the bivector, and then evaluate its Schouten bracket. If $r_A^{\sharp}:\g^{*}\to\g$ is defined by $\beta(r_A^{\sharp}\alpha)=r_A(\alpha,\beta)$, nondegeneracy of $B$ gives
\begin{equation}\label{eq:rsharp-is-A}
        r_A^{\sharp}(X^{\flat})=AX,
        \qquad X\in\g.
\end{equation}
For the Schouten convention used here, a bivector $r$ satisfies\footnote{For this evaluation formula, see \cite{Vaisman1994}, Section~2, and for the identity for the associated operator, see
	\cite{Kupershmidt1999}, equations~(2.9)--(2.13) in Section~2. For the quadratic Rota--Baxter formulation, see
	\cite{Goncharov2019}, equation~(3) and Proposition~1 in Section~3.}
\begin{equation}\label{eq:schouten-evaluation}
\begin{aligned}
\frac12[r,r](\alpha,\beta,\gamma)
={}&\alpha([r^{\sharp}\beta,r^{\sharp}\gamma])
 +\beta([r^{\sharp}\gamma,r^{\sharp}\alpha])+\gamma([r^{\sharp}\alpha,r^{\sharp}\beta]).
\end{aligned}
\end{equation}
Applying \eqref{eq:schouten-evaluation} to
$X^{\flat},Y^{\flat},Z^{\flat}$ and using
\eqref{eq:rsharp-is-A} gives
\begin{equation}\label{eq:J-from-schouten}
\frac12[r_A,r_A](X^{\flat},Y^{\flat},Z^{\flat})
= B([AX,AY],Z)+B([AY,AZ],X)+B([AZ,AX],Y).
\end{equation}
Set
\begin{equation}\label{eq:J-A}
J_A(X,Y,Z)
:=B([AX,AY],Z)+B([AY,AZ],X)+B([AZ,AX],Y).
\end{equation}
Since $B^\flat:\g\to\g^*$ is an isomorphism,
\begin{equation}\label{eq:poisson-iff-J}
        [r_A,r_A]=0
        \quad\Longleftrightarrow\quad
        J_A=0.
\end{equation}
To compare this scalar expression with the operator equation, define
\begin{equation}\label{eq:M-A}
        M_A(X,Y):=[AX,AY]-A\bigl([AX,Y]+[X,AY]\bigr).
\end{equation}
Invariance of $B$ and the $B$-skewness of $A$ yield
\[
\begin{aligned}
B(M_A(X,Y),Z)
&=B([AX,AY],Z)+B([AX,Y],AZ)+B([X,AY],AZ)\\
&=J_A(X,Y,Z).
\end{aligned}
\]
Since $B$ is nondegenerate, $J_A=0$ is equivalent to $M_A=0$, which is \eqref{eq:operator-cybe-lemma}. Over $\K=\mathbb R$, the Schouten bracket of left-invariant multivector fields is left-invariant and equals the algebraic Schouten bracket at the identity. Thus $[r_A,r_A]=0$ is equivalent to $[\pi_A,\pi_A]=0$ on every integrating Lie group (see \cite{Vaisman1994}, Section~2).

\end{proof}

We therefore use the following definition.

\begin{definition}\label{def:poisson-structure}
	\normalfont 
Let $(\g,B)$ be a quadratic Lie algebra. A \emph{Poisson structure} on $(\g,B)$ is a $B$-skew endomorphism $A\in\End(\g)$ satisfying \eqref{eq:operator-cybe-lemma}. It is \emph{degenerate} if $\Ker A\ne0$ and \emph{nondegenerate} if $A$ is invertible.
\end{definition}

\subsection{The characteristic leaf}

By Lemma~\ref{lem:operator-cybe-equivalence}, $A$ is Poisson if and only if $[r_A,r_A]=0$. Over $\mathbb R$, this is the Poisson condition for $\pi_A$ in \eqref{eq:piA-def}. We now describe the symplectic structure on the image of $A$. It gives the algebraic counterpart of the leaf through the identity.

\begin{definition}\label{def:characteristic-leaf}
	\normalfont 
A \emph{symplectic Lie algebra} is a Lie algebra $\h$ equipped with a nondegenerate alternating form $\Omega\in\Lambda^2\h^*$ that is a Chevalley--Eilenberg $2$-cocycle for the trivial action of $\h$ on $\K$. A \emph{symplectic Lie subalgebra} of $\g$ is a Lie subalgebra $\h\subseteq\g$ with such a form. The restriction of an ambient quadratic form to $\h$ need not be nondegenerate.
\end{definition}

For a $B$-skew map $A$, we set $\h_A:=\Imop A$ and define
\begin{equation}\label{eq:leaf-form}
	\Omega_A(AX,AY):=B(AX,Y),
	\qquad X,Y\in\g.
\end{equation}
When $A$ is a Poisson structure, we call $(\h_A,\Omega_A)$ its \emph{characteristic leaf}. Algebraically, it is the subalgebra $\Imop A$ with the symplectic form $\Omega_A$. Over $\mathbb R$, $\h_A$ integrates to the leaf through the identity, and $\Omega_A$ is the value of its symplectic form there. Proposition~\ref{prop:leaf-cosets} gives the precise statement for Lie groups. We often refer simply to the leaf of $A$.

The next result characterizes the Poisson condition in terms of the leaf.

\begin{proposition}\label{prop:leaf-criterion}
Let $(\g,B)$ be a quadratic Lie algebra and let $A\in\End(\g)$ be $B$-skew.  Then
\eqref{eq:leaf-form} is well defined and nondegenerate on the vector space $\h_A$.  Moreover,
the following conditions are equivalent:
\begin{enumerate}
\item $A$ is a Poisson structure.
\item $\h_A$ is a Lie subalgebra of $\g$ and $\Omega_A$ is a Chevalley--Eilenberg
$2$-cocycle on $\h_A$.
\end{enumerate}
Whenever $\h_A$ is a Lie subalgebra, one has the identity
\begin{equation}\label{eq:leaf-differential}
\begin{aligned}
 (d_{\h_A}\Omega_A)(AX,AY,AZ)
 &=-J_A(X,Y,Z)=-B\bigl(M_A(X,Y),Z\bigr),
\end{aligned}
\end{equation}
for all $X,Y,Z\in\g$,
where $J_A$ and $M_A$ are defined by \eqref{eq:J-A} and \eqref{eq:M-A}.
\end{proposition}

\begin{proof}
If $AY=AY'$, then $A(Y-Y')=0$, and $B$-skewness gives
\[
 B(AX,Y-Y')=-B\bigl(X,A(Y-Y')\bigr)=0.
\]
Thus \eqref{eq:leaf-form} is independent of the representative of its second argument. The same calculation for the first argument also shows that the form is alternating. If $AX\in\h_A$ satisfies $\Omega_A(AX,AY)=0$ for every $AY\in\h_A$, then $B(AX,Y)=0$ for every $Y\in\g$. Nondegeneracy of $B$ gives $AX=0$, so $\Omega_A$ is nondegenerate.

If $A$ is a Poisson structure, then
\[
 [AX,AY]=A\bigl([AX,Y]+[X,AY]\bigr)\in\h_A,
\]
so $\h_A$ is a Lie subalgebra. To prove both directions, suppose now only that $\h_A$ is a subalgebra. With the convention
\[
 (d\Omega)(u,v,w)
 =-\Omega([u,v],w)+\Omega([u,w],v)-\Omega([v,w],u),
\]
formula \eqref{eq:leaf-form} gives
\[(d_{\h_A}\Omega_A)(AX,AY,AZ)
=-B([AX,AY],Z)+B([AX,AZ],Y)
-B([AY,AZ],X)\]
\[=-\underset{X,Y,Z}{\mathrm{cyc}}\,B([AX,AY],Z)
=-J_A(X,Y,Z).\]
By the computation in Lemma~\ref{lem:operator-cybe-equivalence}, the last expression is $-B(M_A(X,Y),Z)$. This proves \eqref{eq:leaf-differential}. If $A$ is Poisson, then $M_A=0$ and therefore $d_{\h_A}\Omega_A=0$. Conversely, if $\h_A$ is a subalgebra and $d_{\h_A}\Omega_A=0$, nondegeneracy of $B$ in \eqref{eq:leaf-differential} gives $M_A=0$. Thus $A$ is Poisson.
\end{proof}

At the group level, this algebraic leaf is the usual symplectic leaf through
the identity.

\begin{proposition}\label{prop:leaf-cosets}
Let $A$ be a Poisson structure on a quadratic Lie algebra $(\g,B)$ over
$\mathbb R$, and let $G$ be a connected Lie group that integrates $\g$.  Let
$H_A$ be the connected immersed subgroup with Lie algebra
$\h_A=\Imop A$.  The symplectic leaves of $\pi_A$ are the left cosets of
$H_A$.  The form on $H_A$ is left-invariant and has value $\Omega_A$ at the
identity.  Left translation gives the forms on the other leaves.
\end{proposition}

\begin{proof}
The characteristic distribution of $\pi_A$ is left-invariant and has fiber
$\Imop r_A^\sharp=\Imop A=\h_A$ by \eqref{eq:rsharp-is-A}.  By
Proposition~\ref{prop:leaf-criterion}, $\h_A$ is a Lie subalgebra.  Its
connected immersed subgroup therefore has left cosets as the maximal
integral manifolds.  The symplectic form induced on this subgroup has value $\Omega_A$ at the identity.  Left translation gives the form on every other leaf.
\end{proof}

Conversely, every symplectic Lie subalgebra determines a unique Poisson structure.

\begin{proposition}\label{prop:leaf-converse}
Let $(\g,B)$ be a quadratic Lie algebra, and let $(\h,\Omega)$ be a symplectic Lie subalgebra of
$\g$. Write $\Omega^\flat\colon\h\to\h^*$ for the map $\Omega^\flat(u)=\Omega(u,\cdot)$. Then
\begin{equation}\label{eq:operator-from-leaf}
 A_{\h,\Omega}X
:=-\bigl(\Omega^\flat\bigr)^{-1}\bigl(B(X,\cdot)|_{\h}\bigr),
 \qquad X\in\g,
\end{equation}
where the right-hand side is regarded as an element of $\h\subseteq\g$,
defines the unique Poisson structure whose characteristic leaf is
$(\h,\Omega)$.  Consequently, $A\mapsto(\h_A,\Omega_A)$ is a bijection
between the Poisson structures on $(\g,B)$ and the symplectic Lie subalgebras of
$\g$.
\end{proposition}

\begin{proof}
Nondegeneracy of $B$ makes the restriction map $X\mapsto B(X,\cdot)|_{\h}$ from $\g$ onto $\h^*$ surjective. Since $\Omega$ is nondegenerate, the operator in \eqref{eq:operator-from-leaf} has image $\h$. Define $A=A_{\h,\Omega}$. Its defining identity is
\begin{equation}\label{eq:operator-from-leaf-identity}
        \Omega(AX,u)=-B(X,u),
        \qquad X\in\g,\quad u\in\h.
\end{equation}
Taking $u=AY$ and then interchanging $X$ and $Y$ shows
\[
 B(AX,Y)=\Omega(AX,AY)=-B(X,AY),
\]
so $A$ is $B$-skew and $\Omega_A=\Omega$.  Since $\h$ is a Lie
subalgebra and $d_{\h}\Omega=0$,
Proposition~\ref{prop:leaf-criterion} shows that $A$ is a Poisson structure.
Finally, \eqref{eq:operator-from-leaf-identity} determines $AX$ uniquely for
each $X$, which proves uniqueness and the bijection.
\end{proof}

When $A$ is invertible, the Poisson equation becomes the derivation identity for $A^{-1}$. We recall this standard observation from \cite{BajoBenayadiMedina2007}.

\begin{proposition}\label{prop:invertible-case}
Let $(\g,B)$ be a quadratic Lie algebra and let $A\in\End(\g)$ be
invertible and $B$-skew.  Then $A$ is a Poisson structure if and only if
$\delta:=A^{-1}$ is a $B$-skew derivation of $\g$.  In this case
\[
        \omega(X,Y):=B(\delta X,Y)
\]
is a Chevalley--Eilenberg $2$-cocycle and is nondegenerate.
\end{proposition}

\begin{proof}
The inverse of an invertible $B$-skew map is $B$-skew. Substituting $X=\delta U$ and $Y=\delta V$ in \eqref{eq:operator-cybe-lemma} gives
\[
        [U,V]=A\bigl([U,\delta V]+[\delta U,V]\bigr),
\]
which is equivalent to
\[
        \delta[U,V]=[\delta U,V]+[U,\delta V].
\]
Thus $A$ satisfies the operator equation if and only if $\delta$ is a
derivation.  The form $\omega(X,Y)=B(\delta X,Y)$ is alternating and
nondegenerate.  For its closedness, see \cite{BajoBenayadiMedina2007},
Lemma~1.1. 
\end{proof}

For this inverse formulation of the classical Yang--Baxter equation, see \cite{BajoBenayadiMedina2007}, Remark~2, and \cite{Kupershmidt1999}, Theorem~2.4. With our convention \eqref{eq:leaf-form}, the form in Proposition~\ref{prop:invertible-case} satisfies $\omega=-\Omega_A$. The sign distinguishes the two common conventions for the form defined by an inverse, but does not affect nondegeneracy or closedness.

\begin{remark}\label{rem:JacobsonDerivation}
	\normalfont
A finite-dimensional Lie algebra over a field of characteristic zero that
admits an invertible derivation is nilpotent \cite{Jacobson1955}.  Therefore,  in characteristic zero a Lie
algebra that is both quadratic and symplectic is nilpotent.
\end{remark}

Our convention can also be expressed in the terminology of Lie bialgebras.

\begin{remark}\label{rem:terminology}
\normalfont 
In the language of Lie bialgebras, $r_A$ is a triangular $r$-matrix.  It
defines the coboundary $\delta(X)=[X,r_A]$ and therefore a triangular coboundary Lie bialgebra.  The multiplicative Poisson tensor corresponding to this Lie bialgebra on a group is $r_A^L-r_A^R$, whereas the tensor
studied in this paper is the left-invariant bivector $\pi_A=r_A^L$.  Two
$r$-matrices that differ by an invariant bivector may define the same
cobracket, although they need not define the same $\pi_A$. Equation \eqref{eq:operator-cybe-lemma} also says that $A$ is a Rota--Baxter
operator of weight zero. Thus our Poisson structures are precisely the $B$-skew Rota--Baxter operators of weight zero. See
\cite{Drinfeld1983,Drinfeld1987,KosmannSchwarzbach1997,Goncharov2019,LangSheng2023}
for these correspondences and \cite{Bordemann1990,Kupershmidt1999} for the
operator formulation.
\end{remark}

Diatta and Medina relate solutions of the CYBE to Lagrangian graphs in the cotangent Lie algebra and study the induced affine structures on dual Lie groups (see \cite{DiattaMedina2004}, Theorems~1 and~2). They characterize completeness of the corresponding connection in terms of the symplectic leaf of $r^L$ through the identity. Our question is different: when $B|_{\Imop A}$ is nondegenerate, we reconstruct the ambient quadratic Lie algebra together with its Poisson structure, as in Theorem~\ref{thm:leaf-reconstruction}.

\subsection{Reconstruction when the leaf is quadratic}
\label{sec:leaf-reconstruction}

Propositions~\ref{prop:leaf-criterion} and~\ref{prop:leaf-converse} recover $A$ once its leaf is known. We now fix a symplectic Lie algebra $(\h,\Omega)$ with a nondegenerate invariant form $B_{\h}$ and seek a quadratic Lie algebra containing it as a nondegenerate subalgebra. There are two steps: reconstruct the ambient bracket, and recover the Poisson structure with leaf $(\h,\Omega)$. The first step goes back to Kostant. We use the version for ordinary Lie algebras of \cite{Meyer2021}, Theorem~3.3. This gives a formulation in exterior algebra, and an extension to colour Lie algebras, of \cite{Kostant1999}, Theorem~1.50. We state the construction in our notation and then add the Poisson structure.

The following linear algebra observation will also be used in the converse direction.

\begin{lemma}\label{lem:image-kernel-orthogonal}
Let $(\g,B)$ be a quadratic vector space and let $A\in\End(\g)$ be $B$-skew.  Then
\begin{equation}\label{eq:image-kernel-orthogonal}
        (\Imop A)^{\perp}=\Ker A.
\end{equation}
If the restriction of $B$ to $\Imop A$ is nondegenerate, then
\begin{equation}\label{eq:image-kernel-orthogonal-splitting}
        \g=\Imop A\mathbin{\perp}\Ker A,
\end{equation}
and the restriction $A|_{\Imop A}$ is an automorphism of the vector space $\Imop A$.
\end{lemma}

%\begin{proof}
%For $X,Y\in\g$, skewness gives $B(AX,Y)=-B(X,AY)$.  Thus
%\[
%Y\in(\Imop A)^\perp
%\quad\Longleftrightarrow\quad
%B(X,AY)=0\ \text{for every }X
%\quad\Longleftrightarrow\quad
%AY=0,
%\]
%which proves \eqref{eq:image-kernel-orthogonal}.  If
%$B|_{\Imop A}$ is nondegenerate, then
%$\Imop A\cap\Ker A=0$.  Rank--nullity gives the orthogonal direct sum
%\eqref{eq:image-kernel-orthogonal-splitting}, and
%$A|_{\Imop A}$ is injective and therefore invertible.
%\end{proof}

Let $(\h,B_{\h})$ be a quadratic Lie algebra, and let $(\mathfrak m,B_{\mathfrak m})$ be a finite-dimensional vector space with a symmetric nondegenerate bilinear form. Suppose that
\begin{equation}\label{eq:orthogonal-representation}
        \rho\colon\h\longrightarrow
        \mathfrak{so}(\mathfrak m,B_{\mathfrak m})
\end{equation}
is an orthogonal representation.  Its moment map is the map
\begin{equation}\label{eq:moment-map-beta}
        \beta_{\rho}\colon\Lambda^2\mathfrak m\longrightarrow\h
\end{equation}
defined by
\begin{equation}\label{eq:moment-map-beta-pairing}
        B_{\h}(\beta_{\rho}(u,v),x)
        =B_{\mathfrak m}(\rho(x)u,v),
        \qquad x\in\h,\quad u,v\in\mathfrak m.
\end{equation}
Nondegeneracy of $B_{\h}$ determines $\beta_{\rho}$ uniquely, and skew-symmetry of $\rho(x)$ makes it alternating. To describe the other component of the bracket on $\mathfrak m$, choose an invariant cubic form
\begin{equation}\label{eq:invariant-cubic-gamma}
        \gamma\in(\Lambda^3\mathfrak m^*)^{\h},
\end{equation}
where invariance means
\begin{equation}\label{eq:gamma-invariance}
\begin{aligned}
0={}&\gamma(\rho(x)u,v,w)
   +\gamma(u,\rho(x)v,w)
   +\gamma(u,v,\rho(x)w)
\end{aligned}
\end{equation}
for all $x\in\h$ and $u,v,w\in\mathfrak m$.  Let
\begin{equation}\label{eq:varphi-from-gamma}
        \varphi_{\gamma}\colon\Lambda^2\mathfrak m\longrightarrow\mathfrak m
\end{equation}
be the unique alternating map determined by
\begin{equation}\label{eq:varphi-from-gamma-pairing}
        B_{\mathfrak m}(\varphi_{\gamma}(u,v),w)=\gamma(u,v,w).
\end{equation}
The invariance of $\gamma$ is equivalent to the equivariance identity
\begin{equation}\label{eq:varphi-equivariance}
        \rho(x)\varphi_{\gamma}(u,v)
        =\varphi_{\gamma}(\rho(x)u,v)
         +\varphi_{\gamma}(u,\rho(x)v).
\end{equation}
Indeed, pairing the difference between the two sides with $w$ and using the
$B_{\mathfrak m}$-skewness of $\rho(x)$ gives the negative of the left-hand side of
\eqref{eq:gamma-invariance}.

The moment map is equivariant as well.

\begin{lemma}\label{lem:moment-map-equivariance}
For the orthogonal representation \eqref{eq:orthogonal-representation}, one has
\begin{equation}\label{eq:moment-map-equivariance}
        [x,\beta_{\rho}(u,v)]
        =\beta_{\rho}(\rho(x)u,v)
         +\beta_{\rho}(u,\rho(x)v)
\end{equation}
for all $x\in\h$ and $u,v\in\mathfrak m$.
\end{lemma}

\begin{proof}
For $y\in\h$, invariance of $B_{\h}$, the representation identity, and
skewness of $\rho(x)$ give
\[
\begin{aligned}
B_{\h}([x,\beta_{\rho}(u,v)],y)
&=-B_{\mathfrak m}(\rho([x,y])u,v)=B_{\mathfrak m}(\rho(y)\rho(x)u,v)
+B_{\mathfrak m}(\rho(y)u,\rho(x)v).
\end{aligned}
\]
The defining identity \eqref{eq:moment-map-beta-pairing} and nondegeneracy of $B_{\h}$ now give \eqref{eq:moment-map-equivariance}.
\end{proof}

For the pair $(\rho,\gamma)$, equip $\h\oplus\mathfrak m$ with the orthogonal sum $B:=B_{\h}\mathbin{\perp}B_{\mathfrak m}$ and define
\begin{equation}\label{eq:kostant-combined-map}
        \Psi_{\rho,\gamma}
:=\beta_{\rho}+\varphi_{\gamma}
        \colon\Lambda^2\mathfrak m\longrightarrow
        \h\oplus\mathfrak m.
\end{equation}
The remaining condition for the Jacobi identity is expressed by a $4$-form. Let $\operatorname{Sh}(2,2)\subseteq S_4$ be the set of $(2,2)$-shuffles. The contracted exterior square of $\Psi_{\rho,\gamma}$, called Kostant's obstruction, is
\begin{equation}\label{eq:kostant-quartic-shuffle}
\mathcal K_{\rho,\gamma}(u_1,u_2,u_3,u_4)
:=\frac{1}{2}\sum_{\sigma\in\operatorname{Sh}(2,2)}
\operatorname{sgn}(\sigma)
\cdot B\bigl(
\Psi_{\rho,\gamma}(u_{\sigma(1)},u_{\sigma(2)}),
\Psi_{\rho,\gamma}(u_{\sigma(3)},u_{\sigma(4)})
\bigr).
\end{equation}
Pairing complementary shuffles and using symmetry of $B$ gives the equivalent
formula with three terms
\begin{equation}\label{eq:kostant-quartic}
\begin{aligned}
\mathcal K_{\rho,\gamma}(u,v,w,z)
={}&B_{\h}(\beta_{\rho}(u,v),\beta_{\rho}(w,z))
   -B_{\h}(\beta_{\rho}(u,w),\beta_{\rho}(v,z))\\
&\quad+B_{\h}(\beta_{\rho}(u,z),\beta_{\rho}(v,w))\\
&\quad+B_{\mathfrak m}(\varphi_{\gamma}(u,v),\varphi_{\gamma}(w,z))
   -B_{\mathfrak m}(\varphi_{\gamma}(u,w),\varphi_{\gamma}(v,z))\\
&\quad+B_{\mathfrak m}(\varphi_{\gamma}(u,z),\varphi_{\gamma}(v,w)).
\end{aligned}
\end{equation}
Thus alternation is built into the definition.  In the notation of Meyer
(see \cite{Meyer2021}, Definition~1.20 and Theorem~3.3), the Kostant norm of $\Psi_{\rho,\gamma}$ is
\[
N_{B}(\Psi_{\rho,\gamma})
:=\Psi_{\rho,\gamma}\wedge_{B}\Psi_{\rho,\gamma}
=2\mathcal K_{\rho,\gamma}.
\]
The factor $2$ occurs because Meyer uses the sum over all $(2,2)$-shuffles, whereas
\eqref{eq:kostant-quartic-shuffle} includes the factor $1/2$.  Therefore the vanishing of
$\mathcal K_{\rho,\gamma}$ is exactly the vanishing of the Kostant invariant.

The following theorem explains how this obstruction determines whether the reconstructed bracket satisfies the Jacobi identity.

\begin{theorem}
\label{thm:kostant-quadratic-reconstruction}
Let $(\h,B_{\h})$ be a quadratic Lie algebra and let
$(\mathfrak m,B_{\mathfrak m})$ be a quadratic vector space.  For an orthogonal
representation $\rho$ as in \eqref{eq:orthogonal-representation} and an invariant cubic form $\gamma$ as in \eqref{eq:invariant-cubic-gamma}, define on
\begin{equation}\label{eq:orthogonal-sum-g-h-m}
        \g:=\h\oplus\mathfrak m,
        \qquad B:=B_{\h}\mathbin{\perp}B_{\mathfrak m},
\end{equation}
the skew bilinear operation
\begin{equation}\label{eq:kostant-reconstructed-bracket}
\begin{aligned}
[x+u,y+v]_{\rho,\gamma}
={}&[x,y]_{\h}+\beta_{\rho}(u,v)+\rho(x)v-\rho(y)u+\varphi_{\gamma}(u,v),
\end{aligned}
\end{equation}
where $x,y\in\h$ and $u,v\in\mathfrak m$. The form $B$ is invariant under \eqref{eq:kostant-reconstructed-bracket}, and this operation is a Lie bracket if and only if
\begin{equation}\label{eq:kostant-obstruction-vanishing}
        \mathcal K_{\rho,\gamma}=0.
\end{equation}
Equivalently, after the Jacobi identities containing an element of $\h$ have been imposed, the
only remaining vector equation is
\begin{equation}\label{eq:kostant-vector-obstruction}
        \sum_{\mathrm{cyc}}
        \left(
        \varphi_{\gamma}\bigl(u,\varphi_{\gamma}(v,w)\bigr)
        -\rho(\beta_{\rho}(v,w))u
        \right)=0,
        \qquad u,v,w\in\mathfrak m.
\end{equation}

Conversely, every Lie bracket on the vector space $\h\oplus\mathfrak m$ which restricts to the
given bracket on $\h$, makes $\h$ a Lie subalgebra, and makes
$B_{\h}\perp B_{\mathfrak m}$ invariant is obtained uniquely from a pair $(\rho,\gamma)$ as above satisfying \eqref{eq:kostant-obstruction-vanishing}.  Thus
\eqref{eq:kostant-reconstructed-bracket} gives a bijection between such quadratic Lie brackets
and pairs $(\rho,\gamma)$ for which the Kostant obstruction vanishes.
\end{theorem}

\begin{proof}
We apply the version for ordinary Lie algebras of \cite{Meyer2021}, Theorem~3.3, which extends \cite{Kostant1999}, Theorem~1.50. In Meyer's notation, take
\[
\mathfrak g_{\mathrm{Meyer}}=\h,
\qquad V=\mathfrak m,
\qquad \mu=\beta_{\rho},
\qquad \phi=\varphi_{\gamma}.
\]
The defining equation for $\mu$ is \eqref{eq:moment-map-beta-pairing}. Meyer's two conditions on $\phi$ require $B_{\mathfrak m}(\phi(u,v),w)$ to be alternating and $\h$-invariant. Thus $\phi=\varphi_{\gamma}$ for a unique $\gamma\in(\Lambda^3\mathfrak m^*)^{\h}$. The bracket in \cite{Meyer2021}, Theorem~3.3, extends the bracket of $\h$ and its action on $\mathfrak m$, while on $\Lambda^2\mathfrak m$ it is $\mu+\phi$. These identifications give \eqref{eq:kostant-reconstructed-bracket}. Since $N_B(\mu+\phi)=2\mathcal K_{\rho,\gamma}$ by the normalization above, Meyer's norm condition is \eqref{eq:kostant-obstruction-vanishing}. The $\mathfrak m$-component of the Jacobi identity on three elements of $\mathfrak m$ is \eqref{eq:kostant-vector-obstruction}.

For the converse, invariance and orthogonality imply
$[\h,\mathfrak m]\subseteq\mathfrak m$.  Thus
$\rho(x)u=[x,u]$ is an orthogonal representation.  For
$u,v\in\mathfrak m$ and $x\in\h$,
\[
B_{\h}(\operatorname{pr}_{\h}[u,v],x)
=B_{\mathfrak m}(\rho(x)u,v),
\]
so the $\h$-component of $[u,v]$ is $\beta_{\rho}(u,v)$. The converse in \cite{Meyer2021}, Theorem~3.3 part~(a), determines the unique admissible map $\phi$ and therefore the unique invariant form $\gamma$.
\end{proof}

We now add the symplectic structure on the prescribed subalgebra. Let $(\h,\Omega)$ be a symplectic Lie algebra, and let $B_{\h}$ be a symmetric nondegenerate invariant form on $\h$. Define $D\in\End(\h)$ by
\begin{equation}\label{eq:D-from-Omega}
        \Omega(x,y)=B_{\h}(Dx,y),
        \qquad x,y\in\h.
\end{equation}
Nondegeneracy and alternation of $\Omega$ make $D$ invertible and $B_{\h}$-skew. The cocycle condition for $\Omega$ is equivalent to the derivation identity for $D$ (see \cite{BajoBenayadiMedina2007}, Lemma~1.1). Here this also follows from Propositions~\ref{prop:leaf-criterion} and~\ref{prop:invertible-case}, applied to the invertible $B_{\h}$-skew map $-D^{-1}$. Equation~\eqref{eq:D-from-Omega} fixes the sign used below.

The next theorem combines this construction with
Proposition~\ref{prop:leaf-converse}.  It reconstructs the ambient algebra
and the Poisson structure from the prescribed leaf.

\begin{theorem}
\label{thm:leaf-reconstruction}
Let $(\h,\Omega)$ be a symplectic Lie algebra, let $B_{\h}$ be a symmetric nondegenerate invariant form on $\h$, and let $(\mathfrak m,B_{\mathfrak m})$ be a quadratic vector space. Let $D$ be the invertible $B_{\h}$-skew derivation defined by \eqref{eq:D-from-Omega}. For every orthogonal representation $\rho$ of $\h$ on $\mathfrak m$ and every invariant cubic form $\gamma\in(\Lambda^3\mathfrak m^*)^{\h}$ satisfying
\begin{equation}\label{eq:kostant-obstruction}
        \mathcal K_{\rho,\gamma}=0,
\end{equation}
let $\g_{\rho,\gamma}$ be the quadratic Lie algebra reconstructed by
\eqref{eq:kostant-reconstructed-bracket}.  Define
\begin{equation}\label{eq:reconstructed-poisson-structure}
        A_{\rho,\gamma}(x+u):=-D^{-1}x,
        \qquad x\in\h,\quad u\in\mathfrak m.
\end{equation}
Then $A_{\rho,\gamma}$ is a Poisson structure on
$(\g_{\rho,\gamma},B_{\h}\perp B_{\mathfrak m})$, and
\begin{equation}\label{eq:reconstructed-image-kernel}
        \Imop A_{\rho,\gamma}=\h,
        \qquad
        \Ker A_{\rho,\gamma}=\mathfrak m.
\end{equation}
Its characteristic leaf is $(\h,\Omega)$ and it is the unique Poisson structure with this leaf on the algebra reconstructed above.

Conversely, let $A$ be a Poisson structure on a quadratic Lie algebra
$(\g,B)$ and assume that
\begin{equation}\label{eq:nondegenerate-image-hypothesis}
        \h:=\Imop A
        \quad\text{is nondegenerate for }B.
\end{equation}
Then
\begin{equation}\label{eq:canonical-leaf-kernel-splitting}
        \g=\h\mathbin{\perp}\mathfrak m,
        \qquad \mathfrak m:=\Ker A,
\end{equation}
$B_{\h}:=B|_{\h}$ makes $\h$ quadratic, while $\Omega_A$ makes it
symplectic.  Relative to
\eqref{eq:canonical-leaf-kernel-splitting}, the ambient bracket is reconstructed uniquely
from the data
\begin{equation}\label{eq:converse-rho-gamma}
        \rho(x)u=[x,u],
        \qquad
        \gamma(u,v,w)
        =B\bigl(\operatorname{pr}_{\mathfrak m}[u,v],w\bigr),
\end{equation}
which satisfy \eqref{eq:kostant-obstruction}, and
\begin{equation}\label{eq:converse-A-D-formula}
        A=-D^{-1}\operatorname{pr}_{\h},
        \qquad
        D=-(A|_{\h})^{-1},
\end{equation}
where $D$ is determined by
$\Omega_A(x,y)=B_{\h}(Dx,y)$.

\end{theorem}

\begin{proof}
By Theorem~\ref{thm:kostant-quadratic-reconstruction}, condition \eqref{eq:kostant-obstruction} makes $(\g_{\rho,\gamma},B_{\h}\perp B_{\mathfrak m})$ a quadratic Lie algebra with $\h$ as a nondegenerate quadratic subalgebra. Define $J:=-D^{-1}$. Since $D$ is invertible and $B_{\h}$-skew, so is $J$. Therefore $A_{\rho,\gamma}=J\oplus0$ is $B_{\h}\perp B_{\mathfrak m}$-skew, with image and kernel given by \eqref{eq:reconstructed-image-kernel}. We need to check that its leaf has the prescribed form, with the correct sign. Since $J$ is invertible, any two elements of $\h$ can be written as $Jx$ and $Jy$. By \eqref{eq:leaf-form},
\[
        \Omega_{A_{\rho,\gamma}}(Jx,Jy)=B_{\h}(Jx,y).
\]
On the other hand, $DJ=-\id_{\h}$ and the $B_{\h}$-skewness of $J$ give
\[
\begin{aligned}
        \Omega(Jx,Jy)
        &=B_{\h}(DJx,Jy)
          =-B_{\h}(x,Jy)=B_{\h}(Jx,y).
\end{aligned}
\]
Thus $\Omega_{A_{\rho,\gamma}}=\Omega$. Since $\h$ is a Lie subalgebra and $\Omega$ is closed, Proposition~\ref{prop:leaf-criterion} makes $A_{\rho,\gamma}$ a Poisson structure, and Proposition~\ref{prop:leaf-converse} gives uniqueness. No further compatibility between $D$ and $\rho$ is needed. The operator vanishes on $\mathfrak m$, while its leaf lies in $\h$.

Conversely, suppose $(\g,B,A)$ satisfies \eqref{eq:nondegenerate-image-hypothesis}. Lemma~\ref{lem:image-kernel-orthogonal} gives $\mathfrak m:=\Ker A=(\Imop A)^\perp=\h^\perp$ and the orthogonal decomposition \eqref{eq:canonical-leaf-kernel-splitting}. By Proposition~\ref{prop:leaf-criterion}, $\h$ is a Lie subalgebra and $\Omega_A$ is a nondegenerate $2$-cocycle on it. The restriction $B_{\h}$ is nondegenerate by assumption and invariant because $\h$ is a subalgebra. Thus $\h$ is both quadratic and symplectic. Applying Theorem~\ref{thm:kostant-quadratic-reconstruction} to \eqref{eq:canonical-leaf-kernel-splitting} recovers the bracket uniquely from \eqref{eq:converse-rho-gamma}. In particular, the resulting pair $(\rho,\gamma)$ satisfies \eqref{eq:kostant-obstruction}.

It remains to identify $A$. By Lemma~\ref{lem:image-kernel-orthogonal}, $J:=A|_{\h}$ is invertible as a linear map on $\h$. Define $D$ by $\Omega_A(p,q)=B_{\h}(Dp,q)$. For $x,y\in\h$, take $p=Jx$ and $q=Jy$. Then
\[
\begin{aligned}
        B_{\h}(DJx,Jy)
        &=\Omega_A(Jx,Jy)
          =B_{\h}(Jx,y)=-B_{\h}(x,Jy).
\end{aligned}
\]
Since $J$ is surjective and $B_{\h}$ is nondegenerate, we obtain $DJ=-\id_{\h}$ and $D=-J^{-1}$, so $J=-D^{-1}$. The operator $A$ vanishes on $\mathfrak m=\Ker A$ and equals $J$ on $\h$, proving \eqref{eq:converse-A-D-formula}. Thus $(\g,B,A)$ determines all the reconstruction data.
\end{proof}
Once the leaf $(\h,\Omega)$ and the restricted form $B_{\h}$ are fixed, the theorem describes every quadratic Poisson structure $A$ for which $B$ is nondegenerate on $\Imop A$. The remaining choices are a quadratic space $(\mathfrak m,B_{\mathfrak m})$ and a pair $(\rho,\gamma)$ satisfying \eqref{eq:kostant-obstruction}. In characteristic zero,  the algebra $\h$ is nilpotent by Remark~\ref{rem:JacobsonDerivation}.

\section{Jacobi structures and characteristic leaves}
\label{sec:jacobi-preliminaries}

We now consider bivectors whose failure to satisfy the Poisson equation is controlled by a vector field. This leads to Jacobi structures. We recall their algebraic equations and the distinction between contact and locally conformally symplectic characteristic subalgebras, following \cite{IglesiasMarrero2003}, Definition~A.1, Examples~A.2, and Proposition~A.4 in Appendix~A. The calculations below express these results using the quadratic form and fix the signs needed for reconstruction. They hold over any field $\K$ of characteristic zero, while statements about manifolds and Lie groups are made over $\mathbb R$. For background, see \cite{Lichnerowicz1978,Kirillov1976,Vaisman1994} and \cite{IglesiasMarrero2001}, Section~2.1.

Let $(\g,B)$ be a quadratic Lie algebra. For a $B$-skew endomorphism $C$, let $r_C$ be the bivector in \eqref{eq:rA-def} and let $M_C$ be the Poisson defect \eqref{eq:M-A}, with $A$ replaced by $C$. Our signs are fixed by \eqref{eq:schouten-evaluation}, \eqref{eq:rsharp-is-A}, and the bracket on functions below.

\begin{definition}\label{exp4:def-jacobi-manifold}
A \emph{Jacobi structure} on a smooth real manifold $M$ is a bivector
field $\pi$ and a vector field $E$ satisfying
\[
 [\pi,\pi]=2E\wedge\pi,\qquad [E,\pi]=0.
\]
\end{definition}

With the Schouten convention used in this paper, its bracket on smooth functions is
\begin{equation}\label{exp4:eq-function-jacobi}
	\{f,g\}=\pi(df,dg)-fE(g)+gE(f).
\end{equation}
This bracket is alternating, satisfies the Jacobi identity, and is a differential operator of order at most one in each argument. Our normalization gives $\{1,f\}=-E(f)$, whereas \cite{IglesiasMarrero2001}, Section~2.1, uses $\{1,f\}=E(f)$. To translate into their notation while keeping the bivector and the bracket on functions fixed, take $\Lambda_{\mathrm{IM}}=\pi$ and $E_{\mathrm{IM}}=-E$. Their Schouten normalization (see \cite{IglesiasMarrero2001}, Section~2.2) satisfies $[\pi,\pi]_{\mathrm{IM}}=-[\pi,\pi]$ on bivectors. The tensor equations therefore agree after this change. We use our normalization for all subsequent calculations with operators and forms on leaves.

Using the invariant metric, we express the algebraic definition in these conventions (see \cite{IglesiasMarrero2003}, Definition~A.1 in Appendix~A).

\begin{definition}\label{def:jacobi-structure}
	\normalfont
A \emph{Jacobi structure} on $(\g,B)$ is a pair $(C,e)$, where $C$ is
$B$-skew and $e\in\g$, satisfying
\begin{equation}\label{eq:operator-jacobi-algebraic}
[r_C,r_C]=2e\wedge r_C,
\qquad [e,r_C]=0.
\end{equation}
\end{definition}
Over $\mathbb R$, these equations mean that the left-invariant fields $\pi_C=r_C^L$ and $E_e=e^L$ on an integrating Lie group $G$ satisfy
\begin{equation}\label{eq:jacobi-def}
	[\pi_C,\pi_C]=2E_e\wedge\pi_C,
	\qquad [E_e,\pi_C]=0.
\end{equation}

The next theorem translates the tensor equations into the operator identities used in the reconstruction.

\begin{theorem}\label{thm:operator-jacobi}
The pair $(C,e)$ is a Jacobi structure if and only if
\begin{align}
M_C(X,Y)&=B(e,X)CY-B(e,Y)CX+B(CX,Y)e,
&&X,Y\in\g,
\label{eq:J1}\\
[\ad_e,C]&=0.
\label{eq:J2}
\end{align}
Equation \eqref{eq:J1} is equivalent to
\begin{equation}\label{eq:J1-scalar}
\underset{X,Y,Z}{\mathrm{cyc}}\,B([CX,CY],Z)
=\underset{X,Y,Z}{\mathrm{cyc}}\,B(e,X)B(CY,Z).
\end{equation}
\end{theorem}

\begin{proof}
Write $X^\flat=B(X,\cdot)$.  Lemma~\ref{lem:operator-cybe-equivalence}
already gives the identity
\begin{equation}\label{eq:J-equals-BM}
\begin{aligned}
 \tfrac12[r_C,r_C](X^\flat,Y^\flat,Z^\flat)
 &=\underset{X,Y,Z}{\mathrm{cyc}}\,B([CX,CY],Z)=B(M_C(X,Y),Z).
\end{aligned}
\end{equation}
On the other hand,
\[
\begin{aligned}
 (e\wedge r_C)(X^\flat,Y^\flat,Z^\flat)
 ={}&B(e,X)B(CY,Z)-B(e,Y)B(CX,Z)+B(e,Z)B(CX,Y).
 \end{aligned}
\]
Nondegeneracy of $B$ makes equality of these expressions equivalent to \eqref{eq:J1}, and alternation gives the cyclic form \eqref{eq:J1-scalar}. For the second equation, the adjoint action on a bivector gives $[e,r](\alpha,\beta)=r(\alpha\circ\ad_e,\beta)+r(\alpha,\beta\circ\ad_e)$. Since $\ad_e$ is $B$-skew,
\[
\begin{aligned}
 [e,r_C](X^\flat,Y^\flat)
 &=-B(C\ad_eX,Y)-B(CX,\ad_eY)=B\bigl((\ad_eC-C\ad_e)X,Y\bigr).
\end{aligned}
\]
Nondegeneracy of $B$ now gives the equivalence between $[e,r_C]=0$ and \eqref{eq:J2}.
\end{proof}

For $e=0$, we recover the Poisson equation $M_C=0$ in Definition~\ref{def:poisson-structure}. The Jacobi equations are also preserved by $(C,e)\mapsto(sC,se)$ for $s\in\K^\times$. Each term in \eqref{eq:J1} gains a factor $s^2$, and $[\ad_{se},sC]=s^2[\ad_e,C]$. When $e\neq0$, the tensor equation determines precisely when $C$ can still be Poisson.

\begin{proposition}\label{prop:poisson-rank2}
Let $(C,e)$ be a Jacobi structure with $e\neq0$.  Then $M_C=0$ if and only
if $e\wedge r_C=0$, equivalently $r_C=e\wedge v$ for some $v\in\g$.
When these equivalent conditions hold, $\rank C\leq2$, and if $C\neq0$, then
$\Imop C=\Span\{e,v\}$.
\end{proposition}

\begin{proof}
By \eqref{eq:operator-jacobi-algebraic} and
Lemma~\ref{lem:operator-cybe-equivalence}, $M_C=0$ is equivalent to
$e\wedge r_C=0$.  For a complement $\g=\K e\oplus W$, write
$r_C=e\wedge v+r_W$ with $r_W\in\Lambda^2W$.
The map $e\wedge\cdot:\Lambda^2W\to\Lambda^3\g$ is injective,
so $e\wedge r_C=0$ exactly when $r_W=0$.  Finally,
\[
 r_C^\sharp\lambda=\lambda(e)v-\lambda(v)e
\]
has image $\Span\{e,v\}$ when $r_C\neq0$, and zero image otherwise.
Since $C=r_C^\sharp B^\flat$, the rank assertions follow.
\end{proof}

Define $\h=\Imop C$. As in the Poisson case, the form
\begin{equation}\label{eq:Omega-C}
\Omega_C(CX,CY):=B(CX,Y),
\end{equation}
is well defined and nondegenerate on $\h$, by the same argument as in Proposition~\ref{prop:leaf-criterion}. In particular, $\dim\h$ is even.

The position of $e$ relative to $\h$ leads to contact or locally conformally symplectic structures. We recall the terminology from \cite{IglesiasMarrero2003}, Examples~A.2 in Appendix~A, and \cite{AngellaBazzoniParton2020}, Introduction.

\begin{definition}\label{exp4:def-contact}
\normalfont
Let $M$ be a smooth manifold. 
\begin{itemize}
\item[(i)] If $M$ has dimension $2n+1$, a one-form $\alpha$ is a
\emph{contact form} if $\alpha\wedge(d\alpha)^n$ is nowhere zero.
Its kernel $\ker\alpha$ is the \emph{contact distribution}. The \emph{Reeb vector field} is the unique vector field $R$ satisfying $\alpha(R)=1$ and $\iota_Rd\alpha=0$.
\item[(ii)] A \emph{locally conformally symplectic form} (LCS\@ form for short) on $M$ consists of a nondegenerate two-form $\Omega$ and a closed one-form $\theta$ satisfying $d\Omega=\theta\wedge\Omega$. The form $\theta$ is called a \emph{Lee form}.
\end{itemize}
\end{definition}

On a Lie algebra, contact and locally conformally symplectic structures are defined by the same equations with the Chevalley--Eilenberg differential. Left translation extends them to structures of the same kind on any real integrating Lie group.

For a Jacobi structure on a manifold $M$, define $X_f=\pi^\sharp(df)-fE$. Then \eqref{exp4:eq-function-jacobi} becomes $\{f,g\}=X_f(g)+gE(f)$. The \emph{characteristic distribution} is the span of these vectors at each point. Since the value and differential of a local function can be prescribed independently at a point, this distribution is
\begin{equation}\label{exp4:eq-characteristic-distribution}
 \mathcal D_p=\Imop\pi_p^\sharp+\mathbb R E_p.
\end{equation}
A \emph{characteristic leaf} is a maximal connected immersed integral manifold of this distribution. For a left-invariant structure, the rank is constant. The next two propositions give the operator form of the standard dichotomy for characteristic subalgebras (see \cite{IglesiasMarrero2003}, Proposition~A.4 in Appendix~A). We include the calculations to fix the signs $d\alpha=-\Omega_C$ on the contact distribution and $\theta=-\iota_e\Omega_C$ in the even-dimensional case.

\begin{proposition}\label{prop:contact-type}
Let $(C,e)$ be a Jacobi structure with $e\notin\Imop C$.  Then
$\mathfrak t=\Imop C\oplus\K e$ is a Lie subalgebra, and the form
$\alpha\in\mathfrak t^*$ defined by $\alpha(e)=1$ and $\alpha|_{\Imop C}=0$ satisfies $d\alpha|_{\Imop C\times\Imop C}=-\Omega_C$ and
$\iota_e d\alpha=0$.  Thus, $\alpha$ is a contact form with Reeb element
$e$.
\end{proposition}

\begin{proof}
Define $\h=\Imop C$.  Equation~\eqref{eq:J1} gives
\[
 [CX,CY]=C\bigl([CX,Y]+[X,CY]\bigr)
       +B(e,X)CY-B(e,Y)CX+B(CX,Y)e,
\]
while \eqref{eq:J2} gives $[e,CX]=C[e,X]$.  Therefore
$\mathfrak t=\h\oplus\K e$ is a Lie subalgebra.  Applying $\alpha$
and using $d\alpha(u,v)=-\alpha([u,v])$ yields
\[
 d\alpha(CX,CY)=-B(CX,Y)=-\Omega_C(CX,CY),
 \qquad d\alpha(e,CX)=0.
\]
The restriction to $\ker\alpha=\h$ is nondegenerate, so $\alpha$
is contact. The identities $\alpha(e)=1$ and $\iota_e d\alpha=0$ identify its Reeb
element.  This includes $C=0$, when $\mathfrak t=\K e$.
\end{proof}

When $e$ belongs to the image of $C$, the characteristic subalgebra carries a locally conformally symplectic form. We compute its Lee form and distinguish the two-dimensional case, where the equation for $d\Omega_C$ does not determine that form.

\begin{proposition}\label{prop:lcs-type}
Let $(C,e)$ be a Jacobi structure with $e\in\Imop C$.  Then
$\h=\Imop C$ is a Lie subalgebra, and
\[
\theta=-\iota_e\Omega_C,
\qquad \theta(CX)=-B(e,X),
\qquad d_{\h}\theta=0,
\qquad d_{\h}\Omega_C=\theta\wedge\Omega_C.
\]
In particular, $\Omega_C$ is closed if $\theta=0$.  The converse holds
when $\dim\h\geq4$.  In dimension two, $\Omega_C$ is always closed.
\end{proposition}

\begin{proof}
Define $\h=\Imop C$. Since $e\in\h$, all terms in the bracket computed in the preceding proof belong to $\h$. Thus $\h$ is a Lie subalgebra. By \eqref{eq:Omega-C},
\[
 \theta(CX)=-\Omega_C(e,CX)=-B(e,X),\qquad \theta(e)=0.
\]
With $S=[CX,Y]+[X,CY]$, equation~\eqref{eq:J1} therefore gives
$d_{\h}\theta(CX,CY)=B(e,S)$.  Invariance and \eqref{eq:J2} imply
\[
 B(e,[CX,Y])=B([e,CX],Y)
 =B(C[e,X],Y)=-B(e,[X,CY]),
\]
thus $d_{\h}\theta=0$.

For $\xi\in\h$, one has $\Omega_C(\xi,CZ)=B(\xi,Z)$.  Thus
\eqref{eq:J1-scalar} gives
\[
\begin{aligned}
 d_{\h}\Omega_C(CX,CY,CZ)
 &=-\underset{X,Y,Z}{\mathrm{cyc}}\,B([CX,CY],Z)\\
 &=-\underset{X,Y,Z}{\mathrm{cyc}}\,B(e,X)B(CY,Z)\\
 &=(\theta\wedge\Omega_C)(CX,CY,CZ).
\end{aligned}
\]
Suppose $\dim\h=2n\geq4$. The symplectic isomorphism $\lambda\mapsto\lambda\wedge\Omega_C^{n-1}$ from $\h^*$ to $\Lambda^{2n-1}\h^*$ shows that $\theta\wedge\Omega_C=0$ forces $\theta=0$. In dimension two, every three-form vanishes. Thus $\Omega_C$ is closed even when the Lee form $-\iota_e\Omega_C$ specified by the Jacobi structure is nonzero. Unlike in dimensions at least four, the equation for $d\Omega_C$ alone does not determine the Lee form. If $\h=0$, then $e=0$, and the assertions follow directly.
\end{proof}

Over $\mathbb R$, let $G$ integrate $\g$, and let $T$ be the connected immersed subgroup with Lie algebra $\mathfrak t=\Imop C+\mathbb R e$. Equation~\eqref{exp4:eq-characteristic-distribution} gives $\mathcal D_g=(L_g)_*\mathfrak t$. Its maximal connected leaves are therefore the left cosets $gT$, which need not be embedded or closed. Extend the forms in the two propositions by left translation on $T$, and then to $gT$ by $L_g$. This gives the contact or locally conformally symplectic structure on each leaf in our conventions. In the contact case, the Reeb field is the restriction of $e^L$. In the other case, closedness of the Lee form gives $\theta=df$ locally, and $e^{-f}\Omega_C$ is symplectic. A global potential need not exist.

\subsection{Reconstruction when the characteristic subalgebra is nondegenerate}
\label{subsec:jacobi-nondegenerate}

To extend Theorem~\ref{thm:leaf-reconstruction} to Jacobi structures, we must prescribe the \emph{full characteristic subalgebra}
\begin{equation}\label{eq:jacobi-nd-hypothesis}
 \mathfrak t=\Imop C+\K e,
 \qquad B|_{\mathfrak t}\text{ nondegenerate}.
\end{equation}
The Kostant--Meyer data of Subsection~\ref{sec:leaf-reconstruction} reconstruct the ambient quadratic bracket without any additional obstruction. We first prove this for both types of characteristic subalgebra, and then express the data on the subalgebra through skew operators. All reconstruction statements below are over a field $\K$ of characteristic zero.

We begin by showing that the Jacobi pair is determined by its restriction to this subalgebra.

\begin{lemma}\label{lem:jacobi-nd-splitting}
Let $(C,e)$ be a Jacobi structure on $(\g,B)$ satisfying
\eqref{eq:jacobi-nd-hypothesis}.  Then $(\mathfrak t,B_{\mathfrak t})$,
where $B_{\mathfrak t}=B|_{\mathfrak t}$, is a quadratic Lie algebra,
and
\begin{equation}\label{eq:jacobi-nd-splitting}
 \g=\mathfrak t\mathbin{\perp}\mathfrak m,
 \qquad
 \mathfrak m=\mathfrak t^\perp=\Ker C\cap e^\perp,
 \qquad [\mathfrak t,\mathfrak m]\subseteq\mathfrak m.
\end{equation}
Relative to this decomposition,
\begin{equation}\label{eq:jacobi-nd-support}
 C=C_0\oplus0,\qquad C_0=C|_{\mathfrak t},\qquad e\in\mathfrak t,
 \qquad \Imop C_0+\K e=\mathfrak t.
\end{equation}
The pair $(C_0,e)$ is a Jacobi structure on
$(\mathfrak t,B_{\mathfrak t})$.
\end{lemma}

\begin{proof}
Propositions~\ref{prop:contact-type} and~\ref{prop:lcs-type} show that $\mathfrak t$ is a Lie subalgebra. The restriction of $B$ is therefore invariant, and it is nondegenerate by hypothesis. This gives the orthogonal direct sum. Lemma~\ref{lem:image-kernel-orthogonal} also gives
\[
 (\Imop C+\K e)^\perp=(\Imop C)^\perp\cap e^\perp
                       =\Ker C\cap e^\perp.
\]
For $x,y\in\mathfrak t$ and $u\in\mathfrak m$, invariance gives
$B([x,u],y)=-B(u,[x,y])=0$, proving the last inclusion in
\eqref{eq:jacobi-nd-splitting}.  Since $C\mathfrak m=0$ and
$\Imop C\subseteq\mathfrak t$, we obtain
\eqref{eq:jacobi-nd-support} and, in particular, $\Imop C_0=\Imop C$. By orthogonality, $r_C$ is the natural image of $r_{C_0}\in\Lambda^2\mathfrak t$ in $\Lambda^2\g$. Since inclusion of a Lie subalgebra preserves the Schouten bracket, the tensor equations for $(C,e)$ are precisely those for $(C_0,e)$ viewed in $\Lambda^\bullet\g$.
\end{proof}

For the extension theorem, replace $\h$ by $\mathfrak t$ in the Kostant data of Section~\ref{sec:leaf-reconstruction}. Thus $(\mathfrak m,B_{\mathfrak m})$ is a quadratic vector space, $\rho:\mathfrak t\to\mathfrak{so}(\mathfrak m,B_{\mathfrak m})$ is a representation, and $\gamma\in(\Lambda^3\mathfrak m^*)^{\mathfrak t}$. Define the alternating maps $\beta_\rho$ and $\varphi_\gamma$ by
\begin{equation}\label{eq:jacobi-nd-moment}
B_{\mathfrak t}(\beta_\rho(u,v),x)
=B_{\mathfrak m}(\rho(x)u,v),\qquad  B_{\mathfrak m}(\varphi_\gamma(u,v),w)=\gamma(u,v,w).
\end{equation}
We keep the normalization of $\mathcal K_{\rho,\gamma}$ in
\eqref{eq:kostant-quartic-shuffle}--\eqref{eq:kostant-quartic}, with
$B_{\mathfrak t}$ in place of $B_{\h}$.

The following result reconstructs a Jacobi structure from its characteristic subalgebra when the latter is nondegenerate.
\begin{theorem}
\label{thm:jacobi-nd-reconstruction}
Let $(\mathfrak t,B_{\mathfrak t})$ be a quadratic Lie algebra carrying
a Jacobi structure $(C_0,e_0)$ such that
$\Imop C_0+\K e_0=\mathfrak t$.  Choose
$(\mathfrak m,B_{\mathfrak m},\rho,\gamma)$ as above.  On
$\g=\mathfrak t\oplus\mathfrak m$, with
$B=B_{\mathfrak t}\perp B_{\mathfrak m}$, set
\begin{equation}\label{eq:jacobi-nd-bracket}
[x+u,y+v]_{\rho,\gamma}=[x,y]_{\mathfrak t}+\beta_\rho(u,v)+\rho(x)v-\rho(y)u+\varphi_\gamma(u,v).
\end{equation}
This is a quadratic Lie bracket if and only if
$\mathcal K_{\rho,\gamma}=0$.  Whenever this condition holds,
\begin{equation}\label{eq:jacobi-nd-extension}
 C(x+u)=C_0x,\qquad e=e_0
 \qquad(x\in\mathfrak t,\ u\in\mathfrak m)
\end{equation}
defines a Jacobi structure whose characteristic subalgebra is exactly
$\mathfrak t$ and whose restriction to it is $(C_0,e_0)$.

Conversely, every Jacobi structure satisfying
\eqref{eq:jacobi-nd-hypothesis} is obtained in this way.  Relative to
the canonical decomposition \eqref{eq:jacobi-nd-splitting}, all its
data are uniquely recovered by
\begin{equation}\label{eq:jacobi-nd-recovered-data}
\begin{aligned}
 B_{\mathfrak t}&=B|_{\mathfrak t},\qquad
 B_{\mathfrak m}=B|_{\mathfrak m},\qquad C_0=C|_{\mathfrak t},
 \qquad e_0=e,\\
 \rho(x)u&=[x,u],\qquad
 \gamma(u,v,w)=B\bigl(\operatorname{pr}_{\mathfrak m}[u,v],w\bigr).
\end{aligned}
\end{equation}
For a fixed reconstructed bracket and a fixed pair $(C_0,e_0)$ on the subalgebra, \eqref{eq:jacobi-nd-extension} is the unique Jacobi structure with that restriction and characteristic subalgebra $\mathfrak t$.
\end{theorem}

\begin{proof}
The assertion about the ambient bracket follows from Theorem~\ref{thm:kostant-quadratic-reconstruction}. We give the calculation to show why the obstruction does not involve $C_0$ or $e_0$. Invariance of $B$ follows from invariance on $\mathfrak t$, skewness of $\rho$, the first identity in \eqref{eq:jacobi-nd-moment}, and alternation of $\gamma$. The Jacobi identity on three elements of $\mathfrak t$ holds by assumption. On two elements of $\mathfrak t$ and one of $\mathfrak m$, it is the representation identity for $\rho$. On one element of $\mathfrak t$ and two of $\mathfrak m$, its components are the equivariance identities for $\beta_\rho$ and $\varphi_\gamma$, proved in \eqref{eq:moment-map-equivariance} and~\eqref{eq:varphi-equivariance}.

It remains to consider $u,v,w\in\mathfrak m$.  Write
$J(u,v,w)=\sum_{\mathrm{cyc}}[u,[v,w]]$.  Its two components are
\[
\begin{aligned}
 J(u,v,w)_{\mathfrak t}
   &=\sum_{\mathrm{cyc}}\beta_\rho(u,\varphi_\gamma(v,w)),\\
 J(u,v,w)_{\mathfrak m}
   &=\sum_{\mathrm{cyc}}
      \bigl(\varphi_\gamma(u,\varphi_\gamma(v,w))
                     -\rho(\beta_\rho(v,w))u\bigr).
\end{aligned}
\]
Pairing the first expression with $x\in\mathfrak t$ gives
\[
 B_{\mathfrak t}(J(u,v,w)_{\mathfrak t},x)
 =\gamma(\rho(x)u,v,w)+\gamma(u,\rho(x)v,w)
                          +\gamma(u,v,\rho(x)w)=0,
\]
since $\gamma\in (\Lambda^3m^*)^{\mathfrak t}$ by assumption. 
Pairing the second with $z\in\mathfrak m$, and using alternation of
$\gamma$, gives
\[
\begin{aligned}
 B_{\mathfrak m}(J(u,v,w)_{\mathfrak m},z)
 =&-\sum_{\mathrm{cyc}}
  \bigl(B_{\mathfrak t}(\beta_\rho(u,z),\beta_\rho(v,w))+B_{\mathfrak m}(\varphi_\gamma(u,z),\varphi_\gamma(v,w))\bigr)\\
 =&-\mathcal K_{\rho,\gamma}(u,v,w,z).
\end{aligned}
\]
Nondegeneracy gives the equivalence with the Kostant condition in both directions. We now check the Jacobi structure. The map $C=C_0\oplus0$ is $B$-skew, its bivector lies in $\Lambda^2\mathfrak t$, and $e=e_0\in\mathfrak t$. Therefore
\[
 [r_C,r_C]_{\g}=[r_{C_0},r_{C_0}]_{\mathfrak t}
                 =2e_0\wedge r_{C_0},
 \qquad
 [e,r_C]_{\g}=[e_0,r_{C_0}]_{\mathfrak t}=0.
\]
The equality $\Imop C+\K e=\mathfrak t$ follows from the corresponding equality on the prescribed leaf.

Conversely, Lemma~\ref{lem:jacobi-nd-splitting} gives the canonical orthogonal decomposition and the Jacobi pair on $\mathfrak t$. Invariance and the ambient Jacobi identity make $\rho(x)u=[x,u]$ an orthogonal representation. Pairing $[u,v]$ with $x\in\mathfrak t$ gives the first identity in \eqref{eq:jacobi-nd-moment}, so the $\mathfrak t$-component is $\beta_\rho(u,v)$. The $\mathfrak m$-component defines $\gamma$ through \eqref{eq:jacobi-nd-recovered-data}. Invariance of $B$ makes $\gamma$ alternating, and the Jacobi identity with one argument in $\mathfrak t$ makes it $\mathfrak t$-invariant. The calculation above gives $\mathcal K_{\rho,\gamma}=0$. The bracket and the orthogonal projections determine all these data uniquely. Finally, by Lemma~\ref{lem:jacobi-nd-splitting}, a skew map with characteristic subalgebra $\mathfrak t$ vanishes on $\mathfrak t^\perp$. Its prescribed restriction to $\mathfrak t$ therefore determines the whole pair.
\end{proof}

One can also see directly that the Jacobi equations on mixed pairs impose no further condition relating $(C,e)$ to $\rho$. By Theorem~\ref{thm:operator-jacobi}, for $x\in\mathfrak t$ and $u,v\in\mathfrak m$ we have
\[
 M_C(x,u)=-C([C_0x,u])=0,\qquad M_C(u,v)=0,
 \qquad [\ad_e,C]u=-C(\rho(e)u)=0,
\]
because $C\mathfrak m=0$ and $[\mathfrak t, \mathfrak m]\subset \mathfrak m.$ The right-hand side of \eqref{eq:J1} also vanishes on these pairs, since $C\mathfrak m=0$, $e\perp\mathfrak m$, and $C\mathfrak t\subseteq\mathfrak t$. In particular, there is no requirement that $\rho(e)$ vanish. The complement $\mathfrak m$ need not be an ideal or a Lie subalgebra, and the $\mathfrak t$-component of $[u,v]$ is $\beta_\rho(u,v)$.

We now describe the data on a quadratic leaf of locally conformally symplectic or contact type.
\subsubsection{Quadratic leaves that are locally conformally symplectic}
When $e\in\Imop C$, we have $\mathfrak t=\Imop C$ and $C_0$ is invertible. The leaf data consist of a nondegenerate two-form $\Omega$ and a specified closed Lee form $\theta$ satisfying $d\Omega=\theta\wedge\Omega$. We use the standard Lie algebra definition of a locally conformally symplectic structure \cite{AngellaBazzoniParton2020}, with the sign $\theta=-\iota_e\Omega$ from Proposition~\ref{prop:lcs-type}. The next lemma expresses these equations using the invariant metric. When the Lee form vanishes, it recovers the correspondence between symplectic forms and invertible skew derivations (see \cite{BajoBenayadiMedina2007}, Lemma~1.1).

\begin{lemma}\label{lem:jacobi-lcs-metric}
Let $(\mathfrak t,B_{\mathfrak t})$ be a quadratic Lie algebra with a nondegenerate two-form $\Omega$ and a one-form $\theta$. Define $D$ and $\ell$ by
\begin{equation}\label{eq:jacobi-lcs-metric-data}
 \Omega(x,y)=B_{\mathfrak t}(Dx,y),\qquad
 \theta(x)=B_{\mathfrak t}(\ell,x).
\end{equation}
Then $D$ is invertible and $B_{\mathfrak t}$-skew.  The equations
$d\theta=0$ and $d\Omega=\theta\wedge\Omega$ are equivalent respectively to
$\ell\in Z(\mathfrak t)$ and the identity
\begin{equation}\label{eq:jacobi-lcs-twisted-derivation}
[Dx,y]+[x,Dy]-D[x,y]=B_{\mathfrak t}(\ell,x)Dy-B_{\mathfrak t}(\ell,y)Dx+B_{\mathfrak t}(Dx,y)\ell,
       \qquad x,y\in\mathfrak t.
\end{equation}
\end{lemma}

\begin{proof}
Nondegeneracy and alternation of $\Omega$ make $D$ invertible and skew. Invariance of $B_{\mathfrak t}$ gives
\[
 d\theta(x,y)=-B_{\mathfrak t}(\ell,[x,y])
             =-B_{\mathfrak t}([\ell,x],y).
\]
Thus $d\theta=0$ is equivalent to centrality of $\ell$. To check the other equation, expand the Chevalley--Eilenberg differential and use invariance and skewness:
\begin{equation}\label{eq:jacobi-lcs-dOmega}
 d\Omega(x,y,z)
 =B_{\mathfrak t}([Dx,y]+[x,Dy]-D[x,y],z).
\end{equation}
On the other hand,
\[
\begin{aligned}
 (\theta\wedge\Omega)(x,y,z)
 ={}&B_{\mathfrak t}(\ell,x)B_{\mathfrak t}(Dy,z)
       -B_{\mathfrak t}(\ell,y)B_{\mathfrak t}(Dx,z)+B_{\mathfrak t}(Dx,y)B_{\mathfrak t}(\ell,z).
\end{aligned}
\]
Comparing the two expressions and using nondegeneracy of $B_{\mathfrak t}$ proves the claim.
\end{proof}
These data give the following reconstruction in the locally conformally symplectic case.
\begin{theorem}
\label{thm:jacobi-lcs-reconstruction}
Let $(\mathfrak t,B_{\mathfrak t})$ be a quadratic Lie algebra. Let $D$ be an invertible $B_{\mathfrak t}$-skew operator and let $\ell\in Z(\mathfrak t)$, with these data satisfying \eqref{eq:jacobi-lcs-twisted-derivation}. For every ambient reconstruction \eqref{eq:jacobi-nd-bracket} with $\mathcal K_{\rho,\gamma}=0$, the formulas
\begin{equation}\label{eq:jacobi-lcs-reconstruction}
 C(x+u)=-D^{-1}x,\qquad e=-D^{-1}\ell
 \qquad(x\in\mathfrak t,\ u\in\mathfrak m)
\end{equation}
define a Jacobi structure with
\begin{equation}\label{eq:jacobi-lcs-image-kernel}
 \Imop C=\mathfrak t,\qquad \Ker C=\mathfrak m,
 \qquad \Omega_C=\Omega,\qquad -\iota_e\Omega_C=\theta,
\end{equation}
where $\Omega$ and $\theta$ are given by
\eqref{eq:jacobi-lcs-metric-data}.

Conversely, every Jacobi structure for which $e\in\Imop C$ and
$B|_{\Imop C}$ is nondegenerate is obtained uniquely in this way,
relative to its canonical orthogonal decomposition. The data on its leaf are
\begin{equation}\label{eq:jacobi-lcs-converse-data}
 D=-(C|_{\mathfrak t})^{-1},\qquad \ell=-De,
 \qquad \mathfrak t=\Imop C,
\end{equation}
and its ambient data are \eqref{eq:jacobi-nd-recovered-data}.
In particular, $(\mathfrak t,B_{\mathfrak t},\Omega,\theta)$ and an admissible ambient bracket determine the Jacobi pair uniquely.
\end{theorem}

\begin{proof}
Define $C_0=-D^{-1}$ and $e=C_0\ell$. The inverse of an invertible skew map is skew, so $C_0$ is $B_{\mathfrak t}$-skew. Since $DC_0=-\id$,
\[
 \Omega(C_0x,C_0y)
 =-B_{\mathfrak t}(x,C_0y)=B_{\mathfrak t}(C_0x,y).
\]
Thus $C_0$ defines the form $\Omega$. Also, $De=-\ell$ gives $\iota_e\Omega=-\theta$ and $\theta(e)=0$. Lemma~\ref{lem:jacobi-lcs-metric} gives the two locally conformally symplectic equations. We check both Jacobi equations directly, without assuming the converse of Proposition~\ref{prop:lcs-type}. For $x,y,z\in\mathfrak t$, the identity for the form gives
\[
\begin{aligned}
 d\Omega(C_0x,C_0y,C_0z)
    &=-\sum_{\mathrm{cyc}}B_{\mathfrak t}([C_0x,C_0y],z),\\
 (\theta\wedge\Omega)(C_0x,C_0y,C_0z)
    &=-\sum_{\mathrm{cyc}}B_{\mathfrak t}(e,x)
                                      B_{\mathfrak t}(C_0y,z).
\end{aligned}
\]
Here $\theta(C_0x)=B_{\mathfrak t}(\ell,C_0x)
=-B_{\mathfrak t}(e,x)$.  Equality of these expressions is
\eqref{eq:J1-scalar} for $(C_0,e)$, and therefore gives the first Jacobi equation.
For the second, Cartan's formula on Chevalley--Eilenberg cochains gives
\[
 \mathcal L_e\Omega=d(\iota_e\Omega)+\iota_e(d\Omega)=-d\theta+\theta(e)\Omega-\theta\wedge\iota_e\Omega=0.
\]
Invariance of $B_{\mathfrak t}$ gives, more explicitly,
\[
 (\mathcal L_e\Omega)(x,y)
 =-\Omega([e,x],y)-\Omega(x,[e,y])
 =B_{\mathfrak t}((\ad_eD-D\ad_e)x,y).
\]
It follows that $[\ad_e,D]=0$ and therefore $[\ad_e,C_0]=0$. Thus the second Jacobi equation follows from the prescribed data and requires no additional hypothesis. Theorem~\ref{thm:jacobi-nd-reconstruction} extends this pair to the ambient algebra and gives \eqref{eq:jacobi-lcs-image-kernel}.

Conversely, Lemma~\ref{lem:image-kernel-orthogonal} makes $C_0=C|_{\mathfrak t}$ invertible and gives $\mathfrak m=\Ker C$. Proposition~\ref{prop:lcs-type} supplies $\Omega_C$ and the closed form $\theta=-\iota_e\Omega_C$. The unique map defined by $\Omega_C(x,y)=B_{\mathfrak t}(Dx,y)$ satisfies $DC_0=-\id$ by \eqref{eq:Omega-C}, and the metric dual of $\theta$ is $\ell=-De$. Lemma~\ref{lem:jacobi-lcs-metric} gives centrality of $\ell$ and \eqref{eq:jacobi-lcs-twisted-derivation}. Theorem~\ref{thm:jacobi-nd-reconstruction} gives uniqueness of the ambient data.
\end{proof}

For nonzero $\ell$, the map $D$ need not be a derivation, so the nilpotency conclusion in Remark~\ref{rem:JacobsonDerivation} does not apply. Moreover, $\ell$ is required to be central in $\mathfrak t$, not in the ambient algebra. There is no condition $\rho(\ell)=0$. To see the distinction, extend $\theta$ by zero on $\mathfrak m$ and denote the resulting form by $\widetilde\theta$. Then
\begin{equation}\label{eq:jacobi-lcs-extended-lee}
 d\widetilde\theta(u,v)
 =-B_{\mathfrak t}(\ell,\beta_\rho(u,v))
 =-B_{\mathfrak m}(\rho(\ell)u,v).
\end{equation}
The differential vanishes on all other pairs of summands. Thus the extension is closed on $\g$ exactly when $\rho(\ell)=0$. Closedness on the leaf alone does not imply closedness on the ambient algebra.

In dimension two, the Lee form must be kept as part of the prescribed data. Both sides of $d\Omega=\theta\wedge\Omega$ are zero three-forms, so this equation determines neither $\theta$ nor $e$. For example, on an abelian quadratic plane, any one-form can serve as a Lee form for a fixed nondegenerate $\Omega$. When $\dim\mathfrak t=2n\geq4$, the map $\lambda\mapsto\lambda\wedge\Omega$ is injective. Wedging with $\Omega^{n-2}$ reduces this assertion to the isomorphism $\lambda\mapsto\lambda\wedge\Omega^{n-1}$. In these dimensions, $\Omega$ determines $\theta$ and therefore $e$.

\subsubsection{Quadratic leaves of contact type}
\label{subsec:jacobi-contact-quadratic}

When $e\notin\Imop C$, the full characteristic subalgebra is
$\mathfrak t=\Imop C\oplus\K e$, and Proposition~\ref{prop:contact-type}
gives its contact form $\alpha$. We assume now that $B|_{\mathfrak t}$
is nondegenerate. The invariant metric then determines the Reeb element
and the operator from the metric dual of $\alpha$. In particular, the
restriction of the metric to $\ker\alpha$ is automatically nondegenerate.
We first describe these data and then apply
Theorem~\ref{thm:jacobi-nd-reconstruction}.

\begin{lemma}\label{lem:jacobi-contact-metric}
Let $(\mathfrak t,B_{\mathfrak t})$ be a quadratic Lie algebra of dimension
$2n+1$, and let $\alpha\in\mathfrak t^*$. Define $a\in\mathfrak t$ by
\begin{equation}\label{eq:jacobi-contact-metric-data}
 \alpha(x)=B_{\mathfrak t}(a,x),\qquad x\in\mathfrak t.
\end{equation}
Then $\alpha$ is a contact form if and only if
\begin{equation}\label{eq:jacobi-contact-metric-conditions}
 q:=B_{\mathfrak t}(a,a)\neq0,\qquad
 \Ker(\ad_a)=\K a.
\end{equation}
When these conditions hold, its Reeb element and kernel are
\begin{equation}\label{eq:jacobi-contact-reeb-splitting}
 e=q^{-1}a,\qquad \h=\ker\alpha=a^\perp,\qquad
 \mathfrak t=\K e\mathbin{\perp}\h.
\end{equation}
Here the orthogonal complement is taken in $\mathfrak t$. The map
\begin{equation}\label{eq:jacobi-contact-D}
 D:=\ad_a|_{\h}\colon\h\longrightarrow\h
\end{equation}
is invertible and skew with respect to $B_{\mathfrak t}|_{\h}$, and
\begin{equation}\label{eq:jacobi-contact-Omega}
 \Omega(x,y):=-d_{\mathfrak t}\alpha(x,y)
             =B_{\mathfrak t}(Dx,y),\qquad x,y\in\h,
\end{equation}
is nondegenerate.
\end{lemma}

\begin{proof}
Invariance gives
\begin{equation}\label{eq:jacobi-contact-dalpha}
 d_{\mathfrak t}\alpha(x,y)
 =-B_{\mathfrak t}(a,[x,y])
 =-B_{\mathfrak t}([a,x],y),\qquad x,y\in\mathfrak t.
\end{equation}
Thus the kernel of $d_{\mathfrak t}\alpha$ is $\Ker(\ad_a)$. If $\alpha$
is contact, this kernel is a line on which $\alpha$ does not vanish.
Since $a\neq0$ and $[a,a]=0$, the line is $\K a$, and
$\alpha(a)=B_{\mathfrak t}(a,a)\neq0$. This proves
\eqref{eq:jacobi-contact-metric-conditions}.

Conversely, these conditions give the orthogonal decomposition
$\mathfrak t=\K a\mathbin{\perp}a^\perp$ and make $\alpha$ nonzero on
the kernel of $d_{\mathfrak t}\alpha$. The restriction of
$d_{\mathfrak t}\alpha$ to $\h=\ker\alpha$ is therefore a
nondegenerate alternating form, so $\alpha$ is contact. The normalization $\alpha(e)=1$ gives $e=q^{-1}a$,
and \eqref{eq:jacobi-contact-dalpha} gives $\iota_e d_{\mathfrak t}\alpha=0$.
The map $\ad_a$ preserves $\h$, since
$B_{\mathfrak t}(a,[a,x])=0$, and its restriction has zero kernel.
It is consequently invertible on $\h$. Invariance makes it skew, and
\eqref{eq:jacobi-contact-dalpha} proves
\eqref{eq:jacobi-contact-Omega}.
\end{proof}

The space $\h$ in this lemma is generally not a Lie subalgebra. More
precisely, for $x,y\in\h$,
\begin{equation}\label{eq:jacobi-contact-bracket-component}
 \alpha([x,y])=\Omega(x,y),\qquad
 \operatorname{pr}_{\K e}[x,y]=\Omega(x,y)e.
\end{equation}
Thus it is not a subalgebra whenever $\h\neq0$. In particular, we cannot
apply the Poisson reconstruction to $\h$ alone. The operator $D$ is the
restriction of the inner derivation $\ad_a$ on $\mathfrak t$, but we do
not regard it as a derivation of a Lie algebra structure on $\h$.

These data give the contact counterpart of
Theorem~\ref{thm:jacobi-lcs-reconstruction}.

\begin{theorem}\label{thm:jacobi-c-reconstruction}
Let $(\mathfrak t,B_{\mathfrak t})$ be a quadratic Lie algebra, and choose
$a\in\mathfrak t$ satisfying
\eqref{eq:jacobi-contact-metric-conditions}. Define $\alpha,e,\h,D$ as
in Lemma~\ref{lem:jacobi-contact-metric}, and let
\begin{equation}\label{eq:jacobi-contact-projection}
 P x=x-\alpha(x)e,\qquad x\in\mathfrak t,
\end{equation}
be the orthogonal projection onto $\h$. For every ambient reconstruction
\eqref{eq:jacobi-nd-bracket} with $\mathcal K_{\rho,\gamma}=0$, the formulas
\begin{equation}\label{eq:jacobi-contact-reconstruction}
 C(x+u)=-D^{-1}P x,\qquad e=q^{-1}a
 \qquad(x\in\mathfrak t,\ u\in\mathfrak m)
\end{equation}
define a Jacobi structure of contact type. Its characteristic subalgebra
is $\mathfrak t$, and
\begin{equation}\label{eq:jacobi-contact-image-kernel}
 \Imop C=\h,\qquad
 \Ker C=\K e\mathbin{\perp}\mathfrak m,\qquad
 \Omega_C=-d_{\mathfrak t}\alpha|_{\h\times\h}.
\end{equation}
Its contact form is $\alpha$, with Reeb element $e$.

Conversely, every Jacobi structure for which $e\notin\Imop C$ and
$B|_{\mathfrak t}$ is nondegenerate, where
$\mathfrak t=\Imop C\oplus\K e$, is obtained uniquely in this way,
relative to its canonical orthogonal decomposition. In particular,
$B(e,e)\neq0$, $Ce=0$, and the data on its leaf are
\begin{equation}\label{eq:jacobi-contact-converse-data}
\begin{aligned}
 a&=\frac{e}{B(e,e)},\qquad q=\frac{1}{B(e,e)},\qquad
 \alpha(x)=\frac{B(e,x)}{B(e,e)}\quad(x\in\mathfrak t),\\
 \h&=\Imop C=e^\perp\cap\mathfrak t,\qquad
 D=\ad_a|_{\h}=-(C|_{\h})^{-1}.
\end{aligned}
\end{equation}
Its ambient data are \eqref{eq:jacobi-nd-recovered-data}. Thus a contact
form on a quadratic leaf and an admissible ambient bracket determine
the Jacobi pair uniquely.
\end{theorem}

\begin{proof}
Set $C_0=-D^{-1}P$ on $\mathfrak t$. The inverse of an invertible skew
map is skew, and $P$ is the projection in an orthogonal decomposition.
Therefore $C_0$ is skew, vanishes on $\K a$, and has image $\h$. Moreover,
\begin{equation}\label{eq:jacobi-contact-inverse-identities}
 \ad_a C_0=C_0\ad_a=-P.
\end{equation}
For $x,y\in\h$, these identities give
\[
 \Omega(C_0x,C_0y)
 =B_{\mathfrak t}(DC_0x,C_0y)
 =-B_{\mathfrak t}(x,C_0y)=B_{\mathfrak t}(C_0x,y).
\]
Thus the form defined by $C_0$ on its image is the prescribed $\Omega$. We verify both Jacobi equations directly. For $x,y\in\h$,  define
$S=[C_0x,y]+[x,C_0y]$. Invariance and
\eqref{eq:jacobi-contact-inverse-identities} give
\[
 \alpha(S)=B_{\mathfrak t}(\ad_a C_0x,y)
   +B_{\mathfrak t}(\ad_a x,C_0y)=-B_{\mathfrak t}(x,y)+B_{\mathfrak t}(x,y)=0.
\]
Consequently $S\in\h$. Since $\ad_a$ is a derivation, we get $\ad_a[C_0x,C_0y]=-S$ and  $\ad_a C_0S=-S$. Therefore $M_{C_0}(x,y)\in\Ker(\ad_a)=\K e$. Its coefficient is
\[
 \alpha(M_{C_0}(x,y))
 =\alpha([C_0x,C_0y])
 =-B_{\mathfrak t}(x,C_0y)=B_{\mathfrak t}(C_0x,y).
\]
This proves $M_{C_0}(x,y)=B_{\mathfrak t}(C_0x,y)e$ for all $x,y\in\h$, which is \eqref{eq:J1} on these pairs, since $e\perp\h$. On a mixed
pair $a,y$, with $y\in\h$, we obtain
\[
 M_{C_0}(a,y)=-C_0[a,C_0y]=C_0y.
\]
This is again \eqref{eq:J1}, because $B_{\mathfrak t}(e,a)=1$ and
$C_0a=0$. Bilinearity gives the first Jacobi equation on all of
$\mathfrak t$. Finally, \eqref{eq:jacobi-contact-inverse-identities}
and $e=q^{-1}a$ imply $[\ad_e,C_0]=0$. Thus $(C_0,e)$ is a Jacobi
structure. Theorem~\ref{thm:jacobi-nd-reconstruction} extends it to the
ambient algebra. Its image and kernel are those in
\eqref{eq:jacobi-contact-image-kernel}, and
$\mathfrak t=\h\oplus\K e$ is its characteristic subalgebra.

Conversely, Lemma~\ref{lem:jacobi-nd-splitting} gives the Jacobi pair
$(C_0,e)$ on $\mathfrak t$ and makes $C$ vanish on $\mathfrak m$.
Proposition~\ref{prop:contact-type} gives the unique contact form
$\alpha$ with $\alpha(e)=1$ and $\ker\alpha=\Imop C$.
Apply Lemma~\ref{lem:jacobi-contact-metric} to its metric dual $a$.
It gives $e=a/q$ and therefore
$B(e,e)=q^{-1}\neq0$. It also identifies
$\h=\Imop C$ with $e^\perp\cap\mathfrak t$, so $\h$ is nondegenerate.
Skewness now gives $\Ker C_0=\K e$ and makes $C_0|_{\h}$ invertible.
Since $\Omega_C=-d_{\mathfrak t}\alpha|_{\h\times\h}$, for $x,y\in\h$
we have
\[
 B_{\mathfrak t}(DC_0x,C_0y)
 =\Omega_C(C_0x,C_0y)
 =B_{\mathfrak t}(C_0x,y)=-B_{\mathfrak t}(x,C_0y).
\]
Nondegeneracy and surjectivity of $C_0|_{\h}$ give
$DC_0|_{\h}=-\id_{\h}$. This proves all the formulas in
\eqref{eq:jacobi-contact-converse-data} and recovers
\eqref{eq:jacobi-contact-reconstruction}. The metric dual of $\alpha$
is unique, and Theorem~\ref{thm:jacobi-nd-reconstruction} recovers the
ambient data uniquely.
\end{proof}

Equivalently, one can prescribe the Reeb element rather than the contact
form. The conditions on $e\in\mathfrak t$ and the resulting operator are
\begin{equation}\label{eq:jacobi-contact-reeb-data}
\begin{gathered}
 B_{\mathfrak t}(e,e)\neq0,\qquad \Ker(\ad_e)=\K e,\\
 C|_{\h}=-B_{\mathfrak t}(e,e)(\ad_e|_{\h})^{-1},\qquad
 C|_{\K e\oplus\mathfrak m}=0,\qquad
 \h=e^\perp\cap\mathfrak t.
\end{gathered}
\end{equation}
There is no additional condition that $\rho(e)$ vanish. In particular,
the contact form is prescribed on the characteristic subalgebra, not
on the whole ambient algebra. If $\widetilde\alpha$ denotes its
extension by zero on $\mathfrak m$, then
\begin{equation}\label{eq:jacobi-contact-extended-form}
 d\widetilde\alpha(u,v)
 =-B_{\mathfrak t}(a,\beta_\rho(u,v))
 =-B_{\mathfrak m}(\rho(a)u,v),\qquad u,v\in\mathfrak m.
\end{equation}
Thus the differential of this extension need not vanish on
$\mathfrak m\times\mathfrak m$.

Unlike the locally conformally symplectic case, the contact case imposes
a strong restriction on the quadratic leaf itself. Over $\mathbb R$,
the classification in dimension greater than one is due to Diatta
\cite{Diatta2008}, Theorem~1. We record the consequence for our data,
including its extension to an arbitrary field of characteristic zero.

\begin{proposition}\label{prop:jacobi-contact-dimension}
Over a field $\K$ of characteristic zero, a quadratic Lie algebra
admitting a contact form is either a line or a simple Lie algebra of
dimension three. In the latter case it becomes
$\mathfrak{sl}_2$ after extension to an algebraic closure, and every
$a$ with $B_{\mathfrak t}(a,a)\neq0$ satisfies
\eqref{eq:jacobi-contact-metric-conditions}.
Consequently, a Jacobi structure with a nondegenerate characteristic
subalgebra of contact type has $\rank C=0$ or $\rank C=2$.
Over $\mathbb R$, the leaf of dimension three is isomorphic to
$\mathfrak{sl}_2(\mathbb R)$ or $\mathfrak{su}(2)$.
\end{proposition}

\begin{proof}
Suppose $\dim\mathfrak t>1$. Over $\mathbb R$, the classification is the aforementioned 
Diatta's theorem. The argument proving
semisimplicity in his proof applies unchanged over $\mathbb C$.
By Lemma~\ref{lem:jacobi-contact-metric}, the centralizer of the metric
dual of a contact form has dimension one. Since the dimension of a
centralizer in a complex semisimple Lie algebra is at least its rank,
the complex algebra has rank one and is isomorphic to
$\mathfrak{sl}_2(\mathbb C)$.

For an arbitrary field $\K$ of characteristic zero, choose a subfield
$k\subseteq\K$, finitely generated over $\mathbb Q$, over which the
bracket, the invariant form, and the contact form are defined.
Embed $k$ into $\mathbb C$ and apply the complex classification.
It follows that the algebra has dimension three and a nondegenerate
Killing form. These properties hold over $k$ and remain true after
extension to $\K$ and to an algebraic closure. Thus the algebra is
simple and becomes $\mathfrak{sl}_2$ over an algebraic closure.

In a simple quadratic Lie algebra of dimension three, $\ad_a$ is
nonzero and skew for every $a\neq0$. Since an alternating form has
even rank, $\ad_a$ has rank two and kernel $\K a$. Consequently every
$a$ with $B_{\mathfrak t}(a,a)\neq0$ satisfies
\eqref{eq:jacobi-contact-metric-conditions}. The case of a line is
immediate, and $\dim\mathfrak t=\rank C+1$ gives the asserted
restriction on $\rank C$.
\end{proof}

The case of a line is included in the reconstruction: then $\h=0$,
$C=0$, and $\alpha(e)=1$ is the contact condition. The inverse in
\eqref{eq:jacobi-contact-reconstruction} is understood as the unique
endomorphism of the zero vector space. In particular, on a solvable
quadratic ambient algebra, every nondegenerate characteristic
subalgebra of contact type is a line. Any contact characteristic
subalgebra of dimension at least three in such an ambient algebra must
therefore be degenerate for the invariant form.

\begin{example}
\normalfont
Let $\K$ be  a field of characteristic zero and let $\mathfrak t=\mathfrak{sl}_2(\K)$ have basis $H,U,V$ with
\[
 [H,U]=2U,\qquad [H,V]=-2V,\qquad [U,V]=H,
\]
and take the invariant form determined by
$B_{\mathfrak t}(H,H)=2$, $B_{\mathfrak t}(U,V)=1$, with all other
basis pairings zero except those imposed by symmetry. For $a=H$, we
have $q=2$, $\h=\Span\{U,V\}$, and
\[
 e=\tfrac12H,\qquad C_0H=0,\qquad
 C_0U=-\tfrac12U,\qquad C_0V=\tfrac12V.
\]
Here $\alpha(H)=2$, $\alpha(U)=\alpha(V)=0$, and
$\Omega(U,V)=2$. The associated bivector is
$r_{C_0}=\tfrac12U\wedge V$, so
\[
 [r_{C_0},r_{C_0}]=\tfrac12H\wedge U\wedge V
                 =2e\wedge r_{C_0},\qquad [e,r_{C_0}]=0.
\]
This also verifies the signs in
\eqref{eq:jacobi-contact-reconstruction}. Then combining the classification result of Proposition \ref{prop:jacobi-contact-dimension} with the reconstruction Theorem \ref{thm:jacobi-c-reconstruction} we see that the reconstructed skew operator in every admissible ambient reconstruction extending this pair is zero on $\mathfrak m$.
\end{example}

\section{Induced brackets and cochains}
\label{subsec:shared-cochains}

Before applying generalized double extensions to Poisson and Jacobi structures, we introduce the cochain operators needed for reconstruction. They will give a cohomological interpretation of the Poisson equations and of the inner Jacobi branch when the structure preserves a marked isotropic ideal. We first construct the induced Lie bracket and its differential, and then establish the cochain identity used to describe the obstruction.

Let $\bLie$ be any Lie algebra over $\K$ and
let $d$ be its Chevalley--Eilenberg differential with trivial
coefficients.  For $A\in\End(\bLie)$,  define
\begin{align}
 b_A(x,y):=[x,y]_A&:=[Ax,y]+[x,Ay],
       \label{eq:descendent-A}\\
 M_A(x,y)&:=[Ax,Ay]-A b_A(x,y).
       \label{eq:jacobi-bA}
\end{align}
The operation $b_A$ need not be a Lie bracket.  We write $\widehat d_A$
for the alternating operator on scalar cochains given by
\begin{equation}\label{eq:d-A}
 (\widehat d_A\alpha)(x_0,\ldots,x_k)
 =\sum_{0\leq i<j\leq k}(-1)^{i+j}
 \alpha(b_A(x_i,x_j),x_0,\ldots,\widehat{x_i},\ldots,
                  \widehat{x_j},\ldots,x_k),
\end{equation}
where the remaining arguments keep their order. This operator vanishes on constants. When $b_A$ is a Lie bracket, write $\bLie_A$ for the resulting Lie algebra and $d_A=\widehat d_A$ for its differential.

\begin{lemma}\label{lem:jacobi-descendent}
Let $A\in\End(\bLie)$ and $f_E\in\bLie^*$ satisfy
\begin{equation}\label{eq:descendent-hypotheses}
 M_A(x,y)=f_E(x)Ay-f_E(y)Ax,\qquad f_E(b_A(x,y))=0.
\end{equation}
Then
\begin{equation}\label{eq:jacobi-descendent}
 [x,y]_{A,f_E}=b_A(x,y)+f_E(x)y-f_E(y)x
\end{equation}
is a Lie bracket, defining $\bLie_{A,f_E}$.  The map
$A:\bLie_{A,f_E}\to\bLie$ is a Lie algebra homomorphism, and $f_E$ is
a character of $\bLie_{A,f_E}$.  Its scalar differential satisfies
\begin{equation}\label{eq:jacobi-descendent-differential}
 d_{A,f_E}\alpha=\widehat d_A\alpha-k f_E\wedge\alpha
 \qquad (\alpha\in\Lambda^k\bLie^*).
\end{equation}
Consequently the module $\K_{f_E}$, defined by $x\cdot c=f_E(x)c$, has
Chevalley--Eilenberg differential
\begin{equation}\label{eq:jacobi-twisted-differential}
 \partial_{A,f_E}=d_{A,f_E}+f_E\wedge,
 \qquad
 \partial_{A,f_E}\alpha=\widehat d_A\alpha-(k-1)f_E\wedge\alpha.
\end{equation}
This includes $\partial_{A,f_E}c=c f_E$ in degree zero.  For $f_E=0$ these
constructions reduce to $\bLie_A$ and $d_A$.
\end{lemma}

\begin{proof}
To prove the Jacobi identity, we realize the new bracket through a representation. Let $\mathfrak l=\bLie\oplus\K$ be the direct product Lie algebra, with bracket $[(u,s),(v,t)]=([u,v],0)$, and let $V$ be the underlying vector space of $\bLie$. Define $\rho(u,s)x=[u,x]+sx$. This is a representation of $\mathfrak l$ on $V$, since $[\rho(u,s),\rho(v,t)]=\rho([u,v],0)$. Define $T:V\to\mathfrak l$ by $Tx=(Ax,f_E(x))$. Then
\[
 \rho(Tx)y-\rho(Ty)x
 =b_A(x,y)+f_E(x)y-f_E(y)x
 =[x,y]_{A,f_E}.
\]
The hypotheses \eqref{eq:descendent-hypotheses} give
\[
 A[x,y]_{A,f_E}=[Ax,Ay],
 \qquad
 f_E([x,y]_{A,f_E})=f_E(b_A(x,y))=0.
\]
Consequently
\[
 [Tx,Ty]
 =([Ax,Ay],0)
 =T\bigl(\rho(Tx)y-\rho(Ty)x\bigr),
\]
so $T$ is an $\mathcal O$-operator relative to $\rho$. By \cite{Kupershmidt1999}, Proposition~2.18, the operation $[\cdot,\cdot]_{A,f_E}$ is a Lie bracket and $T$ is a Lie algebra homomorphism. Composing with the two projections shows that $A:\bLie_{A,f_E}\to\bLie$ is a homomorphism and $f_E$ is a character.

We next compare the differentials. Define $s(x,y)=f_E(x)y-f_E(y)x$. Replacing $b_A$ by $b_A+s$ adds the following term on a $k$-form $\alpha$:
\[
 \sum_{i<j}(-1)^{i+j}
 \alpha(f_E(x_i)x_j-f_E(x_j)x_i,
           x_0,\ldots,\widehat{x_i},\ldots,
                    \widehat{x_j},\ldots,x_k).
\]
Fix an index $i$.  For each $j>i$, the term containing $f_E(x_i)$
has $x_j$ in the first slot.  Moving it to its original position in
the list with $x_i$ omitted requires $j-1$ transpositions.  Its sign
is therefore $(-1)^{i+j}(-1)^{j-1}=-(-1)^i$.
For $j<i$, the term containing $f_E(x_i)$ comes from the second
summand for the pair $(j,i)$.  Moving $x_j$ from the first slot to
position $j$ gives the sign
$-(-1)^{j+i}(-1)^j=-(-1)^i$ again.  There are $k$ choices of
$j\neq i$, so the total contribution of $f_E(x_i)$ is
\[
 -k(-1)^i f_E(x_i)
       \alpha(x_0,\ldots,\widehat{x_i},\ldots,x_k).
\]
Summing over $i$ gives $-k f_E\wedge\alpha$, which proves
\eqref{eq:jacobi-descendent-differential}.  For $k=0$, both scalar
differentials are zero, and the same formula holds.

Since scalar multiplications commute and $f_E([x,y]_{A,f_E})=0$, the action $x\cdot c=f_E(x)c$ defines a representation on $\K$. The action terms in the differential with these coefficients are
with these coefficients, the action terms are
\[
 \sum_i(-1)^i f_E(x_i)
       \alpha(x_0,\ldots,\widehat{x_i},\ldots,x_k)
 =(f_E\wedge\alpha)(x_0,\ldots,x_k).
\]
Adding them to $d_{A,f_E}$ proves
\eqref{eq:jacobi-twisted-differential}, including
$\partial_{A,f_E}c=c f_E$ in degree zero.  Finally, the product rule and
$f_E\wedge f_E=0$ give
\[
\begin{aligned}
 (d_{A,f_E}+f_E\wedge)^2\alpha
 &=d_{A,f_E}^2\alpha+d_{A,f_E}(f_E\wedge\alpha)
                         +f_E\wedge d_{A,f_E}\alpha\\
 &=(d_{A,f_E}f_E)\wedge\alpha=0,
\end{aligned}
\]
where $d_{A,f_E}f_E(x,y)=-f_E([x,y]_{A,f_E})=0$.
Thus this operator is a differential.
\end{proof}

\begin{remark}
\normalfont
The differential $d+\chi\wedge$ associated with a character $\chi$ is
the one used in generalized Lie bialgebra and Jacobi theory
(see \cite{IglesiasMarrero2001}, equation~(3.2) in Section~3.1, and
\cite{IglesiasMarrero2003}, Section~2.2 and Appendix~A).
\end{remark}
To relate the original differential to the induced one, we use the following transgression map. For $k\geq1$, define
\begin{equation}\label{eq:transgression-A}
 (\mathscr T_A\alpha)(x_1,\ldots,x_k)
 :=\sum_{i=1}^k
 \alpha(Ax_1,\ldots,Ax_{i-1},x_i,Ax_{i+1},\ldots,Ax_k).
\end{equation}
In degree zero, define $\mathscr T_A=0$. Summing over all positions makes the result alternating, so $\mathscr T_A$ preserves the degree of a cochain.

\begin{lemma}\label{lem:jacobi-transgression}
\label{lem:transgression-cochain-map}
Under \eqref{eq:descendent-hypotheses}, the transgression satisfies
\begin{equation}\label{eq:jacobi-transgression-map}
 \partial_{A,f_E}\mathscr T_A=\mathscr T_A d.
\end{equation}
In particular, if $M_A=0$, then
\begin{equation}\label{eq:transgression-cochain-map}
 d_A\mathscr T_A=\mathscr T_A d.
\end{equation}
\end{lemma}

\begin{proof}
We prove the identity by comparing the coefficient of a formal parameter in two expressions. Define $A_t=A+t\id$, where $t$ is a formal scalar parameter. All expressions are polynomial in $t$, so no limiting argument is needed. Write $s(x,y)=f_E(x)y-f_E(y)x$. Since $A[x,y]_{A,f_E}=[Ax,Ay]$ and $[x,y]_{A,f_E}=b_A(x,y)+s(x,y)$, we have
\[
\begin{aligned}
 A_t[x,y]_{A,f_E}
   &=[Ax,Ay]+t b_A(x,y)+t s(x,y),\\
 [A_tx,A_ty]
   &=[Ax,Ay]+t\bigl([Ax,y]+[x,Ay]\bigr)+t^2[x,y].
\end{aligned}
\]
Subtracting cancels the constant and $b_A$ terms and yields
\begin{equation}\label{eq:jacobi-At-defect}
 A_t[x,y]_{A,f_E}-[A_tx,A_ty]=t s(x,y)-t^2[x,y].
\end{equation}

For a $k$-form $\alpha$ with $k\geq1$, the coefficient of $t$ in
\[
 (A_t^*\alpha)(x_1,\ldots,x_k)
 =\alpha(Ax_1+tx_1,\ldots,Ax_k+tx_k)
\]
is precisely $\mathscr T_A\alpha$: to obtain a term of degree one,
choose $tx_i$ in exactly one slot and $Ax_j$ in every other slot.
Consequently
\[
 A_t^*\alpha=A^*\alpha+t\mathscr T_A\alpha+O(t^2),
\]
where $O(t^2)$ denotes a polynomial divisible by $t^2$.
The same statement for $d\alpha$, which has degree $k+1$, shows
that the coefficient of $t$ in
$d_{A,f_E}A_t^*\alpha-A_t^*d\alpha$ is
\[
 d_{A,f_E}\mathscr T_A\alpha-\mathscr T_A d\alpha.
\]
We can also compute this coefficient by expanding the two differentials. For $x_0,\ldots,x_k\in\bLie$, the difference is
\[
\begin{aligned}
 &\bigl(d_{A,f_E}A_t^*\alpha-A_t^*d\alpha\bigr)(x_0,\ldots,x_k)\\
 &\quad=\sum_{i<j}(-1)^{i+j}
 \alpha\bigl(A_t[x_i,x_j]_{A,f_E}-[A_tx_i,A_tx_j],
       A_tx_0,\ldots,\widehat{A_tx_i},\ldots,
                   \widehat{A_tx_j},\ldots,A_tx_k\bigr).
\end{aligned}
\]
The first argument in every summand has zero constant term by the
homomorphism identity for $A$.  Therefore its coefficient of $t$ is
obtained by taking $s(x_i,x_j)$ in that argument and the constant
term in every other argument.  By \eqref{eq:jacobi-At-defect}, the
result is
\begin{equation}\label{eq:jacobi-transgression-coefficient}
 \sum_{i<j}(-1)^{i+j}
 \alpha\bigl(f_E(x_i)x_j-f_E(x_j)x_i,
 Ax_0,\ldots,\widehat{Ax_i},\ldots,\widehat{Ax_j},\ldots,Ax_k\bigr).
\end{equation}

Fix $i$ and collect the terms containing $f_E(x_i)$.  For each
$j\neq i$, the argument $x_j$ is the unique one left unchanged
among the arguments other than $x_i$.  If $i<j$, moving it from the
first slot to its place in the ordered list with $x_i$ omitted
contributes $(-1)^{j-1}$, so the total sign is $-(-1)^i$.
If $j<i$, the minus sign in $-f_E(x_i)x_j$ and the $j$ transpositions
give $-(-1)^{j+i}(-1)^j=-(-1)^i$ as well.  Summing over all
$j\neq i$ therefore gives
\[
 -(-1)^i f_E(x_i)
 (\mathscr T_A\alpha)(x_0,\ldots,\widehat{x_i},\ldots,x_k).
\]
Summing over $i$ identifies \eqref{eq:jacobi-transgression-coefficient}
with $-f_E\wedge\mathscr T_A\alpha$.  Comparing the two coefficient
computations gives
\[
 d_{A,f_E}\mathscr T_A\alpha-\mathscr T_A d\alpha
 =-f_E\wedge\mathscr T_A\alpha.
\]
Move the last term to the left and use
$\partial_{A,f_E}=d_{A,f_E}+f_E\wedge$ to obtain
\eqref{eq:jacobi-transgression-map}.

For a constant $c$, both sides vanish because
$\mathscr T_Ac=0$ and $dc=0$.  This completes the proof in all
degrees.  When $f_E=0$, the differential is $d_A$, which gives
\eqref{eq:transgression-cochain-map}.
\end{proof}

The cochain identity also gives a map on cohomology.

\begin{remark}
	\normalfont
The transgression sends closed scalar cochains on $\bLie$ to closed cochains on $\bLie_{A,f_E}$ with values in $\K_{f_E}$. It also sends exact cochains to exact cochains, since $\alpha=d\beta$ implies $\mathscr T_A\alpha=\partial_{A,f_E}\mathscr T_A\beta$. It therefore induces maps
\[
H^k(\bLie,\K)\longrightarrow H^k(\bLie_{A,f_E},\K_{f_E}),
\qquad [\alpha]\longmapsto[\mathscr T_A\alpha].
\]
We use these maps below to express the obstruction as a class in degree three.
\end{remark}

\section{Poisson reconstruction along isotropic flags}
\label{sec:gd-reconstruction}

We now reconstruct Poisson structures that preserve an isotropic flag, using the generalized double extension from Section~\ref{sec:double-ext}. The construction passes from quadratic Poisson structures in dimension $n$ to those in dimension $n+2k$. We first express the Poisson equation in terms of the blocks of a skew operator, and then interpret the remaining scalar equation in cohomology. In the invertible case, this recovers and extends the quadratic symplectic constructions of \cite{BajoBenayadiMedina2007}. The same component calculation will also be used for Jacobi structures.

Let $(\tilde{\mathfrak g},T)$ be a generalized double extension with data $(\psi,\alpha,\Gamma)$ satisfying Theorem~\ref{thm:gd-bracket}. Its bracket and form are \eqref{eq:gd-bracket} and \eqref{eq:double-form}. We begin with the block form of a $T$-skew operator and determine when it preserves the flag $I\subseteq I^\perp\subseteq\tilde{\mathfrak g}$. For $V:\bLie\to\g$ and $U:I\to\g$, define
\[
(V^\dagger x)(y)=B(Vy,x),\qquad
f(U^\dagger x)=B(Uf,x).
\]
The superscript $*$ on an endomorphism of $\bLie$ denotes its ordinary
dual.

\begin{lemma}
	\label{lem:gd-skew}\label{prop:double-skew}
Relative to an adapted splitting, every $T$-skew operator has a unique expression
	\begin{equation}\label{eq:gd-six-blocks}
		\mathcal S=
		\begin{pmatrix}
			A&-U^\dagger&Q\\
			V&C&U\\
			\eta^\flat&-V^\dagger&-A^*
		\end{pmatrix},
	\end{equation}
	where $C\in\mathfrak o(\g,B)$, $\eta\in\Lambda^2\bLie^*$, and
	$Q:I\to\bLie$ satisfies $f(Qk)+k(Qf)=0$.
	The condition $\mathcal S(I)\subseteq I$ is equivalent to $U=Q=0$.
	It implies $\mathcal S(I^\perp)\subseteq I^\perp$, and the operator
	then has the form
	\begin{equation}\label{eq:gd-preserving-operator}
		\mathcal S=
		\begin{pmatrix}
			A&0&0\\
			V&C&0\\
			\eta^\flat&-V^\dagger&-A^*
		\end{pmatrix},\qquad
		\mathcal S(y+x+f)=Ay+Vy+Cx+\eta^\flat y-V^\dagger x-A^*f.
	\end{equation}
	The maps it induces on $\mathfrak d/I^\perp$ and $I^\perp/I$ are $A$ and $C$.  The blocks $V$ and $\eta$ depend on the splitting.
\end{lemma}

\begin{proof}
	Apply $T(\mathcal SX,Y)+T(X,\mathcal SY)=0$ to each of the six types
	of pairs from $\bLie,\g,I$.  The pair $(\g,\g)$ makes the block $C:\g\to\g$
	skew.  The pair $(\bLie,I)$ identifies the block $I\to I$
	as $-A^*$.  The pairs $(\bLie,\g)$ and $(\g,I)$ determine
	$-V^\dagger$ and $-U^\dagger$.  The remaining pairs make
	$\eta^\flat$ and $Q$ alternating.  These conditions also make the
	sum vanish for arbitrary $X,Y$, proving sufficiency and uniqueness.
	The image of $f\in I$ has $\bLie$-component $Qf$ and
	$\g$-component $Uf$, so
preserving $I$ is equivalent to $U=Q=0$. Finally, if $X\in I^\perp$ and $f\in I$, then $T(\mathcal SX,f)=-T(X,\mathcal Sf)=0$. Thus the operator preserves the whole flag and induces the stated quotient maps.
\end{proof}

In view of Theorem~\ref{thm:gd-recognition}, we now restrict attention to $T$-skew Poisson operators that preserve the marked isotropic ideal $I$.

\subsection{The Poisson defect}
\label{subsec:gd-defect}

Take the operator \eqref{eq:gd-preserving-operator} on $\widetilde{\mathfrak{g}}=\bLie\oplus\g\oplus I$, where $I=\bLie^*$. Its blocks are $A\in\End(\bLie)$, $C\in\mathfrak o(\g,B)$, $V\in\Hom(\bLie,\g)$, and $\eta\in\Lambda^2\bLie^*$. We will separate the Poisson equation into conditions on $A$ and $C$, conditions coupling the summands through $V$, and a scalar equation for $\eta$. The maps below collect the terms in these conditions. Throughout the component formulas, take $y_i\in\bLie$, $x_i\in\g$, and $f_i\in I=\bLie^*$, as in Lemma~\ref{lem:gd-defect}. Define
\begin{align}
 [y_1,y_2]_A&:=[Ay_1,y_2]_{\bLie}+[y_1,Ay_2]_{\bLie},\notag\\
 M_A(y_1,y_2)&:=[Ay_1,Ay_2]_{\bLie}-A[y_1,y_2]_A,
       \label{eq:gd-base-defect}\\
 K_{y_1}&:=[\psi_{Ay_1}+\ad_{Vy_1},C]-C\psi_{y_1}C,
       \label{eq:gd-K}\\
 L(y_1,y_2)&:=[Vy_1,Vy_2]_{\g}+(\psi_{Ay_1}-C\psi_{y_1})Vy_2\notag\\
       &\quad-(\psi_{Ay_2}-C\psi_{y_2})Vy_1-V[y_1,y_2]_A,
       \label{eq:gd-H}\\
 \alpha_A(y_1,y_2)&:=\alpha(Ay_1,y_2)+\alpha(y_1,Ay_2),
       \label{eq:gd-alpha-A}\\
 L^{\mathrm{gen}}(y_1,y_2)&:=L(y_1,y_2)+\alpha(Ay_1,Ay_2)-C\alpha_A(y_1,y_2).
       \label{eq:gd-Hgen}
\end{align}
Their domains and targets are
\[
 \begin{gathered}
 [\cdot,\cdot]_A,\ M_A:\Lambda^2\bLie\longrightarrow\bLie,\qquad
 K:\bLie\longrightarrow\End(\g),\quad y_1\longmapsto K_{y_1},\\
 L,\ \alpha_A,\ L^{\mathrm{gen}}:\Lambda^2\bLie\longrightarrow\g.
 \end{gathered}
\]
In \eqref{eq:gd-K}, the bracket is the commutator in $\End(\g)$ and $\ad_{Vy_1}$ is computed in $\g$. Thus $K_{y_1}x_2\in\g$ for $y_1\in\bLie$ and $x_2\in\g$. The defect $M_C:\Lambda^2\g\to\g$ uses the bracket of $\g$:
\[
 M_C(x_1,x_2)=[Cx_1,Cx_2]_{\g}
       -C\bigl([Cx_1,x_2]_{\g}+[x_1,Cx_2]_{\g}\bigr),
       \qquad x_1,x_2\in\g.
\]
The final scalar equation involves the three-forms $\kappa_V,\Xi_{A,V}\in\Lambda^3\bLie^*$. For $y_1,y_2,y_3\in\bLie$, define
\begin{align}
 \kappa_V(y_1,y_2,y_3)
   &=\sum_{\mathrm{cyc}(1,2,3)}B(\psi_{y_1}Vy_2,Vy_3),
       \label{eq:gd-kappa}\\
 \Xi_{A,V}(y_1,y_2,y_3)
   &=(\mathscr T_A\Gamma)(y_1,y_2,y_3)+\sum_{\mathrm{cyc}(1,2,3)}B(\alpha_A(y_1,y_2),Vy_3)
       +\kappa_V(y_1,y_2,y_3).
       \label{eq:gd-Xi}
\end{align}
All three terms in \eqref{eq:gd-Xi} are alternating. For $\kappa_V$, this follows from skewness of $\psi_{y_1}$, and for $\sum_{\mathrm{cyc}(1,2,3)}B(\alpha_A(y_1,y_2),Vy_3)$ it follows from alternation of $\alpha_A$. We use $d_A:\Lambda^2\bLie^*\to\Lambda^3\bLie^*$ for the cochain operator of Section~\ref{subsec:shared-cochains}, even before $[\cdot,\cdot]_A$ is known to be a Lie bracket. Thus $d_A\eta+\Xi_{A,V}\in\Lambda^3\bLie^*$ is a three-form.

The following calculation shows how these maps occur in the Poisson defect.

\begin{lemma}
\label{lem:gd-defect}
For $y_i\in\bLie$, $x_i\in\g$, and $f_i\in I=\bLie^*$, the Poisson defect $M_{\mathcal{S}}$ of $\mathcal{S}$ has components given by 
\begin{align}
 M_{\mathcal S}(y_1,y_2)
   &=M_A(y_1,y_2)+L^{\mathrm{gen}}(y_1,y_2)+d_A\eta(y_1,y_2,\cdot)+\Xi_{A,V}(y_1,y_2,\cdot),
       \label{eq:gd-defect-bb}\\
 M_{\mathcal S}(y_1,x_2)
   &=K_{y_1}x_2-B(L^{\mathrm{gen}}(y_1,\cdot),x_2),
       \label{eq:gd-defect-bg}\\
 M_{\mathcal S}(x_1,x_2)
   &=M_C(x_1,x_2)+ B(K_{(\cdot)}x_1,x_2),
       \label{eq:gd-defect-gg}\\
 M_{\mathcal S}(y_1,f_2)
   &=- f_2(M_A(y_1,\cdot)),
       \label{eq:gd-defect-bI}\\
 M_{\mathcal S}(x_1,f_2)&=0,\qquad M_{\mathcal S}(f_1,f_2)=0.
       \label{eq:gd-defect-zero}
\end{align}
Here $y_3\in\bLie$ is the argument of each displayed covector. In \eqref{eq:gd-defect-bb}, the terms $M_A$, $L^{\mathrm{gen}}$, and $(d_A\eta+\Xi_{A,V})(y_1,y_2,\cdot)$ belong to $\bLie$, $\g$, and $I$, respectively. In \eqref{eq:gd-defect-bg} and \eqref{eq:gd-defect-gg}, the terms $K_{y_1}x_2$ and $M_C(x_1,x_2)$ belong to $\g$, while the covectors belong to $I$. The defect in \eqref{eq:gd-defect-bI} lies entirely in $I$.  In particular, $\alpha$ contributes to $L^{\mathrm{gen}}$ and $\Xi_{A,V}$, while $\Gamma$ contributes only to $\Xi_{A,V}$. Neither affects $K$.
\end{lemma}

\begin{proof}
Write two elements of the adapted splitting as
\[
 X_1=y_1+x_1+f_1,\qquad X_2=y_2+x_2+f_2
       \quad\text{in }\bLie\oplus\g\oplus\bLie^*,
\]
where $y_i\in\bLie$, $x_i\in\g$, and $f_i\in I=\bLie^*$.
Their images are
\[
 \mathcal SX_i=Ay_i+(Vy_i+Cx_i)
       +(\eta^\flat y_i-V^\dagger x_i-A^*f_i),\qquad i=1,2,
\]
where $Ay_i$ is the $\bLie$-component, $Vy_i+Cx_i$ is the
$\g$-component, and $\eta^\flat y_i-V^\dagger x_i-A^*f_i$ is the
$I$-component.  Recall that
\[
 M_{\mathcal S}(X_1,X_2)
   =[\mathcal SX_1,\mathcal SX_2]
      -\mathcal S\bigl([\mathcal SX_1,X_2]+[X_1,\mathcal SX_2]\bigr).
\]
Bilinearity and skew-symmetry reduce the computation to six pairs of summands. For each pair, set the other components of $X_1,X_2$ equal to zero. We identify covectors in $I$ by evaluating them at an arbitrary $y_3\in\bLie$.

First take $X_1=y_1$ and $X_2=y_2$. The $\bLie$- and $\g$-components of $[\mathcal Sy_1,y_2]+[y_1,\mathcal Sy_2]$ are $[y_1,y_2]_A$ and $\alpha_A(y_1,y_2)+\psi_{y_1}Vy_2-\psi_{y_2}Vy_1$. The corresponding components of $[\mathcal Sy_1,\mathcal Sy_2]$ are $[Ay_1,Ay_2]$ and $\alpha(Ay_1,Ay_2)+\psi_{Ay_1}Vy_2-\psi_{Ay_2}Vy_1+[Vy_1,Vy_2]$. Subtracting the components of $\mathcal S([\mathcal Sy_1,y_2]+[y_1,\mathcal Sy_2])$ gives $M_A(y_1,y_2)$ and $L^{\mathrm{gen}}(y_1,y_2)$. For the $I$-component, evaluate at $y_3$ and group the terms by the data they involve. The terms containing $\eta$ are
\begin{align*}
 &\eta(y_1,[Ay_2,y_3])+\eta(y_1,[y_2,Ay_3])
       -\eta(y_2,[Ay_1,y_3])\\
 &\quad-\eta(y_2,[y_1,Ay_3])-\eta([y_1,y_2]_A,y_3)\\
 &\quad=-\eta([y_1,y_2]_A,y_3)+\eta([y_1,y_3]_A,y_2)
       -\eta([y_2,y_3]_A,y_1)=( d_A\eta)(y_1,y_2,y_3).
\end{align*}
The terms involving $\Gamma$ are
\[
 \Gamma(Ay_1,Ay_2,y_3)+\Gamma(Ay_1,y_2,Ay_3)
       +\Gamma(y_1,Ay_2,Ay_3)
       =(\mathscr T_A\Gamma)(y_1,y_2,y_3).
\]
The terms involving $\alpha$ are
\[
 -B(\alpha(Ay_1,y_3),Vy_2)+B(\alpha(Ay_2,y_3),Vy_1)
-B(\alpha(y_1,Ay_3),Vy_2)+B(\alpha(y_2,Ay_3),Vy_1)+\]
   \[    +B(\alpha_A(y_1,y_2),Vy_3)=\sum_{\mathrm{cyc}(1,2,3)}B(\alpha_A(y_1,y_2),Vy_3).\]

The terms involving $\psi$ and two copies of $V$ are $B(\psi_{y_3}Vy_1,Vy_2)
       +B(\psi_{y_1}Vy_2-\psi_{y_2}Vy_1,Vy_3)
       =\kappa_V(y_1,y_2,y_3)$. Thus the $I$-component is $(d_A\eta+\Xi_{A,V})(y_1,y_2,\cdot)$, proving
\eqref{eq:gd-defect-bb}.

Next take $X_1=y_1$ and $X_2=x_2$. The defect has zero $\bLie$-component, and its $\g$-component is
\begin{align*}
 &\psi_{Ay_1}Cx_2+[Vy_1,Cx_2]
      -C\bigl(\psi_{Ay_1}x_2+[Vy_1,x_2]+\psi_{y_1}Cx_2\bigr)\\
 &\quad=\bigl([\psi_{Ay_1}+\ad_{Vy_1},C]-C\psi_{y_1}C\bigr)x_2
       =K_{y_1}x_2.
\end{align*}
Skew-symmetry of $\mathcal S$ and invariance of $T$ make the pairing of the defect with a third argument alternating. Therefore
\[
\begin{aligned}
 T(M_{\mathcal S}(y_1,x_2),y_3)
   &=-T(M_{\mathcal S}(y_1,y_3),x_2)=-B(L^{\mathrm{gen}}(y_1,y_3),x_2).
\end{aligned}
\]
This identifies the $I$-component and proves \eqref{eq:gd-defect-bg}.

For $X_1=x_1$ and $X_2=x_2$, the $\bLie$-component is again zero, and the $\g$-component is
\[
 [Cx_1,Cx_2]-C\bigl([Cx_1,x_2]+[x_1,Cx_2]\bigr)=M_C(x_1,x_2).
\]
Alternation identifies the $I$-component through
\[
 T(M_{\mathcal S}(x_1,x_2),y_3)
   =T(M_{\mathcal S}(y_3,x_1),x_2)=B(K_{y_3}x_1,x_2).
\]
This proves \eqref{eq:gd-defect-gg}.

For $X_1=y_1$ and $X_2=f_2$, the defect takes values in $I$, since $[\widetilde{\mathfrak{g}},I]\subseteq I$ and $\mathcal S(I)\subseteq I$. Its value at $y_3$ is
$$M_{\mathcal S}(y_1,f_2)(y_3)=f_2(A[Ay_1,y_3])+f_2(A[y_1,Ay_3])-f_2([Ay_1,Ay_3])=-f_2(M_A(y_1,y_3)).$$
This proves \eqref{eq:gd-defect-bI}.

Finally, take $(X_1,X_2)=(x_1,f_2)$ or $(X_1,X_2)=(f_1,f_2)$. Since $[I,\g\oplus I]=0$, $\mathcal Sx_1\in\g\oplus I$, and $\mathcal Sf_i\in I$, every bracket in the defect vanishes. Thus $M_{\mathcal S}(x_1,f_2)=0$ and $M_{\mathcal S}(f_1,f_2)=0$, proving \eqref{eq:gd-defect-zero}. Reversing the arguments gives the remaining pairs by skew-symmetry.

These six formulas determine
$M_{\mathcal S}(X_1,X_2)$ for the two elements above by bilinearity.
\end{proof}

\needspace{5\baselineskip}
Setting these components equal to zero gives the reconstruction criterion.

\begin{theorem}
\label{thm:gd-poisson}
The operator \eqref{eq:gd-preserving-operator} on a generalized double
extension is Poisson if and only if
\begin{align}
 M_A&=0, &M_C&=0,\label{eq:gd-P12}\\
 K_y&=0\quad(y\in\bLie),&
 L^{\mathrm{gen}}(y,z)&=0\quad(y,z\in\bLie),\label{eq:gd-P34}\\
 d_A\eta+\Xi_{A,V}&=0.&&\label{eq:gd-P5}
\end{align}
Here $\bLie_A=(\bLie,[\cdot,\cdot]_A)$ is a Lie algebra and $d_A$ is its scalar
differential.  
\end{theorem}

\begin{proof}
If $M_{\mathcal S}=0$, the three components of \eqref{eq:gd-defect-bb} and the $\g$-components of \eqref{eq:gd-defect-bg} and \eqref{eq:gd-defect-gg} give the stated conditions. Since $M_A=0$, the operation $[\cdot,\cdot]_A$ is a Lie bracket and $d_A$ is its Chevalley--Eilenberg differential. Conversely, these conditions make every component in Lemma~\ref{lem:gd-defect} vanish. The operator is skew by Lemma~\ref{lem:gd-skew}, so it is Poisson.
\end{proof}

Combining this criterion with Theorem~\ref{thm:gd-recognition} gives a reconstruction whenever a specified nonzero isotropic ideal is preserved.

\begin{corollary}\label{cor:Sol_Nil_Poisson}
Let $(\mathfrak d,T)$ be a quadratic Lie algebra with a Poisson operator $\mathcal S$, and let $I\neq0$ be a marked isotropic ideal satisfying $\mathcal S(I)\subseteq I$. Then every adapted Witt splitting identifies $(\mathfrak d,T,\mathcal S)$ with the generalized double extension in Theorem~\ref{thm:gd-poisson}. Its quadratic quotient $I^\perp/I$ has dimension $\dim\mathfrak d-2\dim I$. This applies, in particular, to nilpotent and solvable quadratic Poisson algebras with such an ideal.
\end{corollary}

\begin{proof}
Theorem~\ref{thm:gd-recognition} reconstructs the quadratic bracket from the specified ideal. Lemma~\ref{lem:gd-skew} gives the operator in the adapted splitting, and Theorem~\ref{thm:gd-poisson} gives its compatibility equations. The dimension formula follows from $\dim I^\perp=\dim\mathfrak d-\dim I$.
\end{proof}
	
For a nonzero indecomposable quadratic Lie algebra that is neither simple nor one dimensional, a minimal nonzero proper ideal is isotropic, as explained after Corollary~\ref{cor:reconstruction_nilpotent_solvable}. The same reconstruction applies if the Poisson operator preserves the chosen ideal. We next specialize the equations to ordinary double extensions and $T^*_{\theta}$-extensions.
	\begin{corollary}\label{cor:ordi_Poisson}
Let $\mathcal S_{\mathrm{ord}}$ be a $T_q$-skew operator preserving $I$ on an ordinary double extension with form \eqref{eq:gd-ordinary-extra-form}. Transport it to the canonical adapted presentation by $\mathcal S=F_q^{-1}\mathcal S_{\mathrm{ord}}F_q$, where $F_q$ is given by \eqref{eq:gd-ordinary-isometry}, and write $\mathcal S$ as in \eqref{eq:gd-preserving-operator}.
Then $\mathcal S_{\mathrm{ord}}$ is Poisson if and only if the
equations in Theorem~\ref{thm:gd-poisson} hold with
\[
 L^{\mathrm{gen}}=L,\qquad
 \Xi_{A,V}=\mathscr T_A\Gamma_q+\kappa_V,\qquad
 \Gamma_q(y,z,w)=-\tfrac12q([y,z]_{\bLie},w),
\]
where $L$ and $\kappa_V$ are given by
\eqref{eq:gd-H} and \eqref{eq:gd-kappa}, respectively.
In particular, for $q=0$ one has $F_q=\id$ and
$\Xi_{A,V}=\kappa_V$.

For a $T^*_{\theta}$-extension, the operator \eqref{eq:gd-preserving-operator} is Poisson if and only if
\begin{align}
 M_A&=0, &\ d_A\eta+{\mathscr T}_A\Gamma&=0,\label{eq:tstar-ext}
\end{align}
where in this case $\Gamma(y_1, y_2, y_3)=\theta(y_1, y_2)(y_3). $
\end{corollary}

\subsection{Cohomological interpretation}
\label{subsec:gd-coho1}
The conditions on $A,C,V$ leave the equation $d_A\eta+\Xi_{A,V}=0$ to determine $\eta$. To interpret it cohomologically, identify a scalar cochain $\omega\in\Lambda^k\bLie^*$ with a multivector in $\Lambda^k I\subseteq\Lambda^k\widetilde{\g}$ through $I=\bLie^*$. We first show that the Schouten bracket with the bivector $\pi_{\mathcal S}$ associated with \eqref{eq:gd-preserving-operator} acts on this exterior algebra as $d_A$. Here $d_A$ is given by \eqref{eq:gd-scalar-differential} with $[\cdot,\cdot]_A$ in place of the original bracket, without yet assuming $d_A^2=0$. We then prove that \eqref{eq:gd-P12}--\eqref{eq:gd-P34} make $\Xi_{A,V}$ a cocycle, so it defines a class in $H^3(\bLie_A,\K)$. Throughout, we work in the fixed adapted presentation $\widetilde{\g}=\bLie\oplus\g\oplus I$, with $I=\bLie^*$.

\begin{proposition}
\label{prop:gd-bivector}\label{prop:double-bivector-decomp}
For a basis $(y_i)$ of $\bLie$, let $(f^i)$ be its dual basis in
$\bLie^*=I$.  Then the bivector $\pi_{\mathcal S}$ associated with the operator \eqref{eq:gd-preserving-operator} is given by 
\begin{equation}\label{eq:gd-bivector}
 \pi_{\mathcal S}
   =\pi_C+\sum_i f^i\wedge(Ay_i+Vy_i)
            +\frac12\sum_{i,j}\eta(y_i,y_j)f^i\wedge f^j.
\end{equation}
The expression is independent of the basis.  For every scalar cochain
$\omega\in\Lambda^k\bLie^*=\Lambda^k I$, one has
\begin{equation}\label{eq:gd-schouten-differential}
 [\pi_{\mathcal S},\omega]=d_A\omega.
\end{equation}
Here the Schouten bracket is taken in $\Lambda^\bullet\widetilde{\g}$,
with $\omega$ and $d_A\omega$ viewed as multivectors in
$\Lambda^k I$ and $\Lambda^{k+1}I$, respectively.
\end{proposition}

\begin{proof}

We recover $\pi_{\mathcal S}$ from \eqref{eq:rA-def} on $(\widetilde{\g},T)$, that is, from $\pi_{\mathcal S}(X^\flat,Y^\flat)=T(\mathcal SX,Y)$ with $X^{\flat}=T(X,\cdot).$ Choose a basis $(E_\alpha)$ of $\widetilde{\g}$ and its $T$-dual basis $(E^\alpha)$, so $T(E_\alpha,E^\beta)=\delta_\alpha^\beta$. Expanding the bivector in this basis gives
\[
 \pi_{\mathcal S}=\frac12\sum_{\alpha,\beta}
 T(\mathcal SE_\alpha,E_\beta)E^\alpha\wedge E^\beta.
\]
We apply this reconstruction to the adapted splitting
$\widetilde{\g}=\bLie\oplus\g\oplus I$.
Choose a basis $(u_a)$ of $\g$ with $B$-dual basis $(u^a)$,
so that $B(u_a,u^b)=\delta_a^b$.
By \eqref{eq:double-form}, the basis $(y_i,u_a,f^i)$ of
$\widetilde{\g}$ has $T$-dual basis $(f^i,u^a,y_i)$.
The block expression \eqref{eq:gd-preserving-operator} gives
\[
\begin{aligned}
 T(\mathcal S y_i,y_j)&=\eta(y_i,y_j),&
 T(\mathcal S y_i,u_a)&=B(Vy_i,u_a),\\
 T(\mathcal S y_i,f^j)&=f^j(Ay_i),&
 T(\mathcal Su_a,u_b)&=B(Cu_a,u_b),\\
 T(\mathcal S u_a,f^i)&=0,&
 T(\mathcal Sf^i,f^j)&=0.
\end{aligned}
\]
Substituting these pairings and the corresponding $T$-dual vectors
into the reconstruction formula, and using the identities for dual bases
\[
 \sum_a B(Vy_i,u_a)u^a=Vy_i,
 \qquad \sum_j f^j(Ay_i)y_j=Ay_i,
\]
we obtain
\[
 \pi_{\mathcal S}
 =\pi_C+\sum_i f^i\wedge(Ay_i+Vy_i)
       +\frac12\sum_{i,j}\eta(y_i,y_j)f^i\wedge f^j,
\]
This is \eqref{eq:gd-bivector}, with $\pi_C=\frac12\sum_{a,b}B(Cu_a,u_b)u^a\wedge u^b$. The expression is independent of the bases because the nondegenerate pairing $T$ determines a unique bivector from $T(\mathcal SX,Y)$.

We now compute its Schouten bracket with cochains in $\Lambda^\bullet I$. Every factor of $\pi_C$, $\sum_i f^i\wedge Vy_i$, and the last term of $\pi_{\mathcal S}$ lies in $\g\oplus I=I^\perp$. Since $[I^\perp,I]=0$, the derivation rule shows that only $\rho_A=\sum_i f^i\wedge Ay_i$ contributes to $[\pi_{\mathcal S},\omega]$ for $\omega\in\Lambda^\bullet I$. The result therefore lies in $\Lambda^\bullet I$. For $f\in\bLie^*=I$,
\[
 [\rho_A,f]=\sum_i f^i\wedge\ad_{Ay_i}^*f,
\]
because $[y,f]=\ad_y^*f$ by \eqref{eq:gd-bracket}.
Regard this element of $\Lambda^2 I=\Lambda^2\bLie^*$ as a scalar
two-cochain and evaluate it on $y,z\in\bLie$:
\begin{align*}
 [\rho_A,f](y,z)
 &=\sum_i\bigl(-f^i(y)f([Ay_i,z])
                    +f^i(z)f([Ay_i,y])\bigr)\\
 &=-f([Ay,z])+f([Az,y])=-f([y,z]_A)=(d_Af)(y,z).
\end{align*}
Thus \eqref{eq:gd-schouten-differential} holds in degree one.
The restriction $D:=[\pi_{\mathcal S},\cdot]$ to
$\Lambda^\bullet I=\Lambda^\bullet\bLie^*$ and $d_A$ are both
derivations of degree one on this exterior algebra, even when
$[\cdot,\cdot]_A$ is not a Lie bracket.  More explicitly, on a
decomposable cochain they both obey
\[
 D(f_1\wedge\cdots\wedge f_k)
 =\sum_{j=1}^k(-1)^{j-1}
    f_1\wedge\cdots\wedge Df_j\wedge\cdots\wedge f_k.
\]
They agree on each $f_j\in\bLie^*$ by the preceding calculation
and vanish on scalars.  Since decomposable cochains span every
exterior power, this proves \eqref{eq:gd-schouten-differential}
in every degree.
\end{proof}

We use this correspondence to prove that $\Xi_{A,V}$ is a three-cocycle on $\bLie_A$. The equation $d_A\eta+\Xi_{A,V}=0$ will then be solvable exactly when its cohomology class vanishes.

\begin{theorem}
\label{thm:gd-poisson-obstruction}
Fix $A,C,V$ satisfying \eqref{eq:gd-P12}--\eqref{eq:gd-P34}.  Then
$d_A\Xi_{A,V}=0$.  A two-form $\eta$ completing the data to a Poisson
structure exists exactly when
\begin{equation}\label{eq:gd-poisson-H3}
 [\Xi_{A,V}]=0\quad\text{in }H^3(\bLie_A,\K).
\end{equation}
If $\eta_0$ is one solution, all solutions in the fixed splitting are
\begin{equation}\label{eq:gd-poisson-affine}
 \eta_0+Z^2(\bLie_A,\K).
\end{equation}
For an ordinary double extension with invariant form $q$, the obstruction is $[\kappa_V+{\mathscr T}_A\Gamma_q]$. For a $T^*_{\theta}$-extension, it is $[{\mathscr T}_A\Gamma].$
\end{theorem}
\begin{proof}
Construct $\pi_{\mathcal S}$ with the prescribed $A,C,V$ and with $\eta=0$, without imposing the scalar equation. By Lemma~\ref{lem:gd-defect}, the only possibly nonzero component of the Poisson defect is $M_{\mathcal S}(y,z)=\Xi_{A,V}(y,z,\cdot)\in I=\bLie^*$. Define $P:=\tfrac12[\pi_{\mathcal{S}}, \pi_{\mathcal S}]$, write $X^\flat=T^\flat X=T(X,\cdot)$, and let $\pi:\widetilde{\g}\to\bLie=\widetilde{\g}/I^\perp$ be the quotient map. The Schouten--defect identity in Lemma~\ref{lem:operator-cybe-equivalence}, namely \eqref{eq:J-from-schouten} and the pairing computation that follows it, gives
\[
 P(X^\flat,Y^\flat,Z^\flat)=T\bigl(M_{\mathcal S}(X,Y),Z\bigr)=\Xi_{A,V}(\pi X,\pi Y,\pi Z),
\]
for all $X,Y,Z\in\widetilde{\g}$. The second equality uses bilinearity and the vanishing of the other components, since
\[M_{\mathcal S}(X,Y)=M_{\mathcal S}(\pi X, \pi Y)=\Xi_{A,V}(\pi X,\pi Y,\cdot)\in I=\bLie^*,\]
and $T(f,Z)=f(\pi Z)$ for $f\in I=\bLie^*$.
Thus the three-form corresponding to $P$ under $T$ is $(\Lambda^3T^\flat)P=\pi^*\Xi_{A,V}$. It vanishes whenever an argument lies in $I^\perp=\ker\pi$, so it belongs to $\Lambda^3\Ann(I^\perp)$, where $\Ann(I^\perp)\subseteq\widetilde{\g}^*$ is the annihilator of $I^\perp$. Since $(I^\perp)^\perp=I$, nondegeneracy gives $T^\flat(I)=\Ann(I^\perp)$. Applying the inverse metric identification yields $P\in\Lambda^3 I$. More precisely, $I=\bLie^*$ and the evaluation formula in Proposition~\ref{prop:gd-bivector} give
\[
 P=\tfrac12[\pi_{\mathcal S},\pi_{\mathcal S}]=\Xi_{A,V}
       \quad\text{in }\Lambda^3 I=\Lambda^3\bLie^*
       \subseteq\Lambda^3\widetilde{\g}.
\]
On the right-hand side, the scalar three-cochain is viewed as a
trivector in the summand $I=\bLie^*$.

The Schouten bracket is a graded Lie bracket of degree $-1$.
For three copies of the bivector $\pi_{\mathcal S}$, the graded Jacobi identity
therefore gives $[\pi_{\mathcal S},[\pi_{\mathcal S},\pi_{\mathcal S}]]=0$.  Applying
\eqref{eq:gd-schouten-differential} directly on
$\Lambda^\bullet I=\Lambda^\bullet\bLie^*$, we obtain
\[
 0=[\pi_{\mathcal S},\tfrac12 [\pi_{\mathcal S}, \pi_{\mathcal S}]]=[\pi_{\mathcal S},\Xi_{A,V}]=d_A\Xi_{A,V}
       \quad\text{in }\Lambda^4 I=\Lambda^4\bLie^*.
\]
The equation $M_A=0$ has already made $\bLie_A$ a Lie algebra,
so $d_A$ is its Chevalley--Eilenberg differential with scalar coefficients and
$\Xi_{A,V}$ is a $d_A$-cocycle.

By Theorem~\ref{thm:gd-poisson}, $\eta$ completes the operator to a Poisson structure exactly when $d_A\eta=-\Xi_{A,V}$. This is solvable if and only if the cocycle is exact, which is \eqref{eq:gd-poisson-H3}. If $\eta_0$ is one solution, any other satisfies $d_A(\eta-\eta_0)=0$. Conversely, $d_A(\eta_0+\xi)=-\Xi_{A,V}$ for every $\xi\in Z^2(\bLie_A,\K)$. The special cases follow by substituting the corresponding expressions for $\Xi_{A,V}$.
\end{proof}

The same bivector formula describes the Poisson structure induced on the quotient by $I$.
\begin{proposition}
\label{cor:double-lift}\label{rem:double-lifting}
For every Poisson operator in Theorem~\ref{thm:gd-poisson}, the quotient
map $q:\widetilde{\g}\to\widetilde{\g}/\bLie^*$ sends $\pi_{\mathcal S}$ to $\pi_C$.
The algebra $\g$ is an ideal of $\widetilde{\g}/\bLie^*$, and $\pi_C$ is a
Poisson bivector there.  Over $\mathbb R$, any group morphism integrating $q$
is a Poisson map for the left-invariant bivectors corresponding to these structures.
In the quotient group those bivectors are tangent to the left cosets
of the connected immersed subgroup integrating $\g$, and their
restriction to that subgroup is the structure defined by $C$.
\end{proposition}

\begin{proof}
For a linear map $q$, its exterior square acts by
$(\Lambda^2q)(X\wedge Y)=qX\wedge qY$.
Every term of \eqref{eq:gd-bivector} other than $\pi_C$ has a factor
$f^i\in \bLie^*$, and $qf^i=0$.  Therefore
$(\Lambda^2q)(\pi_{\mathcal S})=\pi_C$.
Modulo $\bLie^*$, the brackets with $x\in\g$ are $[y,x]=\psi_yx$ and $[w,x]=[w,x]_{\g}$, both in $\g$. Therefore, $\g$ is an ideal in the quotient.  Its own bracket is unchanged,
so the inclusion $\g\hookrightarrow\widetilde{\g}/\bLie^*$ commutes with
all brackets and therefore with Schouten brackets of bivectors.
Since $M_C=0$, we have $[\pi_C,\pi_C]=0$ in $\Lambda^3\g$ and
also in $\Lambda^3(\widetilde{\g}/\bLie^*)$.

For the group assertion, let $Q:D\to\overline D$ be a group
morphism integrating $q$, and write $L_\sigma$ for left translation.
The identity $Q\circ L_\sigma=L_{Q(\sigma)}\circ Q$ implies, after
differentiating at the identity, $T_\sigma Q\circ T_1L_\sigma=T_1L_{Q(\sigma)}\circ q$. Thus, taking exterior squares and applying the result to $\pi_{\mathcal S}$
gives precisely $(\Lambda^2T_\sigma Q)\bigl((\Lambda^2T_1L_\sigma)\pi_{\mathcal S}\bigr)
     =(\Lambda^2T_1L_{Q(\sigma)})\pi_C$. These are the values of the two left-invariant bivectors, so $Q$
is Poisson.  Let $G$ be the connected immersed subgroup of
$\overline D$ with Lie algebra $\g$.  At a point $\overline d$
of a left coset, its tangent space is
$T_1L_{\overline d}(\g)$.  Since $\pi_C\in\Lambda^2\g$,
its left translate belongs to the exterior square of this tangent
space.  This proves tangency to every such coset.  On $G$ itself,
left translation and the inclusion of $G$ in $\overline D$ give
exactly the left-invariant bivector defined by $C$.
None of these arguments require the chosen copy of $\g$ to be a
subalgebra before quotienting by $\bLie^*$.
\end{proof}

\subsection{Quadratic symplectic extensions}
\label{subsc-Poisson special cases}

When the Poisson operator is invertible, its inverse is a skew derivation and defines a symplectic form. We now express the reconstruction in terms of this inverse. Every quadratic symplectic Lie algebra with an isotropic ideal preserved by its associated invertible skew derivation admits the presentation in Proposition~\ref{prop:gd-invertible-bbm}. We then recover double extensions by a line and $T^*_{\theta}$-extensions in Proposition~\ref{prop:gd-bbm-specializations}.

The correspondence between symplectic forms and invertible skew derivations is classical (see \cite{BajoBenayadiMedina2007}, Lemma~1.1 and Remark~2). The extension criteria in \cite{BajoBenayadiMedina2007}, Proposition~2.2 and Lemma~4.1, treat $T^*_{\theta}$-extensions and double extensions by a line. Fischer's related equations assume that the algebra corresponding to $\g$ is abelian (see \cite{Fischer2019}, Definition~3.4 and Theorem~5.4). Here $\g\simeq I^\perp/I$ need not be abelian, and $\psi$ need not be a representation. We keep all the data $(\psi,\alpha,\Gamma)$ and all off-diagonal blocks. The next proposition gives the necessary and sufficient conditions in a fixed generalized double-extension presentation.

\begin{proposition}
\label{prop:gd-invertible-bbm}
Let $\mathcal S$ be the $T$-skew operator \eqref{eq:gd-preserving-operator} on a generalized double extension, preserving $I$. It is invertible if and only if $A$ and $C$ are invertible, with the $\g$-block omitted when $\g=0$. Suppose this holds and define
\begin{equation}\label{eq:gd-inverse-data}
\begin{aligned}
 P&=A^{-1},\qquad D=C^{-1},\qquad W=-DVP,\\
 \beta(y,z)&=\eta(Py,Pz)+B(DVPy,VPz),\qquad y,z\in\bLie.
\end{aligned}
\end{equation}
Then $\beta\in\Lambda^2\bLie^*$ and
\begin{equation}\label{eq:gd-inverse-operator}
 \widetilde D:=\mathcal S^{-1}
   =\begin{pmatrix}
       P&0&0\\
       W&D&0\\
       \beta^\flat&-W^\dagger&-P^*
     \end{pmatrix}.
\end{equation}
The operator $\mathcal S$ is Poisson if and only if $\widetilde D$ is a $T$-skew derivation, or equivalently if and only if
\begin{equation}\label{eq:gd-inverse-symplectic-form}
 \omega(X,Y):=T(\widetilde DX,Y)
\end{equation}
is a symplectic form on $\widetilde{\g}$.  In this case
$\omega=-\Omega_{\mathcal S}$ with the convention
\eqref{eq:leaf-form}. More explicitly, the derivation condition for $\widetilde D$ is equivalent to
\begin{align}
 P&\in\Der(\bLie),\qquad D\in\Dera(\g),
       \label{eq:gd-inverse-der-base}\\
 [D,\psi_y]&=\psi_{Py}+\ad_{Wy},
       \label{eq:gd-inverse-der-mixed}\\
 D\alpha(y,z)+W[y,z]_{\bLie}
   &=\alpha(Py,z)+\alpha(y,Pz)
       +\psi_yWz-\psi_zWy,
       \label{eq:gd-inverse-der-alpha}\\
 (d\beta)(y,z,w)+(\mathscr I_P\Gamma)(y,z,w)
   &+\sum_{\mathrm{cyc}(y,z,w)}B(\alpha(y,z),Ww)=0,
       \label{eq:gd-inverse-der-scalar}
\end{align}
for all $y,z,w\in\bLie$, where $d$ is the scalar differential
of the original algebra $\bLie$ and
\begin{equation}\label{eq:gd-inverse-LP}
 (\mathscr I_P\Gamma)(y,z,w)
   :=\Gamma(Py,z,w)+\Gamma(y,Pz,w)+\Gamma(y,z,Pw).
\end{equation}
Conversely, every quadratic symplectic Lie algebra with an isotropic ideal preserved by its associated invertible skew derivation admits such a presentation, with an invertible Poisson operator $\mathcal S$ preserving that ideal.
\end{proposition}

\begin{proof}
The triangular block form gives the criterion for invertibility. Solving for the three components of $\mathcal S^{-1}(y+x+f)$ in succession gives
\[
\begin{aligned}
 \mathcal S^{-1}(y+x+f)
   ={}&Py+(Dx-DVPy)\\
      &+P^*\eta^\flat Py+P^*V^\dagger DVPy
          -P^*V^\dagger Dx-P^*f.
\end{aligned}
\]
Since $D$ is $B$-skew, $W^\dagger=P^*V^\dagger D$ and
$B(DVPy,VPz)$ is alternating in $y,z$.  This proves
\eqref{eq:gd-inverse-data}--\eqref{eq:gd-inverse-operator}.
Proposition~\ref{prop:invertible-case}, applied to $(\widetilde{\g},T)$, gives the equivalence with the derivation condition. To relate it to the symplectic form, take a $T$-skew endomorphism $E$ and define $\omega_E(X,Y):=T(EX,Y)$. Invariance of $T$ and \eqref{eq:gd-scalar-differential} give
\begin{equation}\label{eq:gd-skew-derivation-differential}
 (d\omega_E)(X,Y,Z)
   =T\bigl([EX,Y]+[X,EY]-E[X,Y],Z\bigr).
\end{equation}
This proves the equivalence with closedness. Nondegeneracy follows
from invertibility.  Finally,
$\Omega_{\mathcal S}(X,Y)=T(X,\mathcal S^{-1}Y)
=-T(\mathcal S^{-1}X,Y)$, which gives the stated sign.

We now substitute \eqref{eq:gd-inverse-operator} in the derivation identity for \eqref{eq:gd-bracket}. The $\bLie$-component on $(y,z)$ and the $\g$-component on $(x,u)$ give \eqref{eq:gd-inverse-der-base}. The $\g$-components on $(y,x)$ and $(y,z)$ give \eqref{eq:gd-inverse-der-mixed} and \eqref{eq:gd-inverse-der-alpha}, respectively. Evaluating the $I$-component of $[\widetilde Dy,z]+[y,\widetilde Dz]-\widetilde D[y,z]$ at $w\in\bLie$ gives
\[
 (d\beta)(y,z,w)+(\mathscr I_P\Gamma)(y,z,w)
       +\sum_{\mathrm{cyc}(y,z,w)}B(\alpha(y,z),Ww).
\]
This gives \eqref{eq:gd-inverse-der-scalar}, including its sign.
Identity~\eqref{eq:gd-skew-derivation-differential} shows that pairing
the derivation defect with $T$ is alternating in three arguments.
It therefore determines the remaining $I$-components on $(y,x)$
and $(x,u)$ from the components already considered.  The pair
$(y,f)$ gives only the derivation identity for $P$, and the pairs
$(x,f)$ and $(f,k)$ vanish.  Thus there are no further conditions.
Setting $V=\eta=0$ gives the diagonal specialization.

Conversely, let $(\mathfrak d,T,\omega)$ be quadratic symplectic and define $E$ by $\omega(X,Y)=T(EX,Y)$. By \eqref{eq:gd-skew-derivation-differential}, $E$ is an invertible $T$-skew derivation. If an isotropic ideal $I$ is $E$-invariant, finite dimensionality gives $E(I)=I$, so $E^{-1}$ also preserves $I$. Choose an adapted Witt splitting by Theorem~\ref{thm:gd-recognition}, and use Lemma~\ref{lem:gd-skew} to write $E^{-1}$ as in \eqref{eq:gd-preserving-operator}. Its diagonal blocks are invertible, and Proposition~\ref{prop:invertible-case} makes it Poisson.
\end{proof}
For a diagonal operator, the equations simplify as follows.
\begin{remark}
If $\mathcal S=A\oplus C\oplus(-A^*)$, then $W=\beta=0$ and $\widetilde D=P\oplus D\oplus(-P^*)$. The equations become
\[
\begin{gathered}
	P\in\Der(\bLie),\qquad D\in\Dera(\g),\qquad
	[D,\psi_y]=\psi_{Py},\\
	D\alpha(y,z)=\alpha(Py,z)+\alpha(y,Pz),\qquad
	\mathscr I_P\Gamma=0.
\end{gathered}
\]
\end{remark}
We can now compare with the two extension constructions of Bajo--Benayadi--Medina \cite{BajoBenayadiMedina2007}. For a $T^*_{\theta}$-extension, pulling the scalar obstruction equation back along $P=A^{-1}$ gives the comparison directly.

\begin{proposition}
\label{prop:gd-bbm-specializations}
Let $\mathcal S$ be an invertible $T$-skew operator of the form \eqref{eq:gd-preserving-operator}, preserving $I$. Use the inverse data $P,D,W,\beta$ from \eqref{eq:gd-inverse-data}, and omit the $\g$-block when $\g=0$.
\begin{enumerate}
\item[\textnormal{(i)}]
For an ordinary double extension by a line, take $\bLie=\K e$, let $e^*\in I$ satisfy $e^*(e)=1$, and use the ordinary presentation with $q=0$ and $\psi_e=\delta\in\Dera(\g)$. Write $Ae=ae$ and $Ve=v$, with $a\ne0$, and define
\begin{equation}\label{eq:gd-bbm-parameters}
 \lambda=-a^{-1},\qquad c=a^{-1}Dv.
\end{equation}
Then $\mathcal S$ is Poisson exactly when $D$ is an invertible $B$-skew derivation satisfying
\begin{equation}\label{eq:gd-bbm-compatibility}
 [\delta,D]-\lambda\delta=\ad_c.
\end{equation}
Its inverse has the form
\begin{equation}\label{eq:gd-bbm-double-derivation}
 \widetilde D e=-\lambda e-c,\qquad
 \widetilde D x=Dx+B(c,x)e^*,\qquad
 \widetilde D e^*=\lambda e^*,\qquad x\in\g.
\end{equation}
Every quadratic symplectic Lie algebra whose associated invertible skew derivation preserves a central line admits such a presentation. Over an algebraically closed field of characteristic zero, every nonzero quadratic symplectic Lie algebra has such a line.

\item[\textnormal{(ii)}]
For a quadratic $T^*_{\theta}$-extension, take $\g=0$ and a cyclic coadjoint two-cocycle $\theta$, so $\Gamma(y,z,w)=\theta(y,z)(w)$. Define
\begin{equation}\label{eq:gd-bbm-tstar-data}
\begin{aligned}
 F(y,z)&:=-\eta(Py,Pz),\\
 \Theta_P(y,z,w)&:=\theta(Py,z)(w)+\theta(Pz,w)(y)
                         +\theta(Pw,y)(z).
\end{aligned}
\end{equation}
Then $\mathcal S$ is Poisson exactly when
\begin{equation}\label{eq:gd-bbm-tstar-coboundary}
 P\in\Der(\bLie)\text{ is invertible},\qquad dF=\Theta_P.
\end{equation}
More precisely, under the equivalent base conditions $M_A=0$ and
$P\in\Der(\bLie)$, the map $P:\bLie\to\bLie_A$ is a Lie algebra
isomorphism.  For scalar cochains, write
\[
 (P^*\xi)(y_1,\ldots,y_k):=\xi(Py_1,\ldots,Py_k),
 \qquad\text{so that }F=-P^*\eta.
\]
The pull-back of the scalar Poisson equation is
\begin{equation}\label{eq:gd-bbm-tstar-pullback}
 P^*(d_A\eta+\mathscr T_A\Gamma)=-dF+\Theta_P.
\end{equation}
Thus $d_A\eta+\mathscr T_A\Gamma=0$ becomes the equation $\Theta=dF$ of Bajo--Benayadi--Medina, with their derivation $D$ identified with $P$ and their $\Theta$ with $\Theta_P$. Conversely, given an invertible derivation $P$ and a two-form $F$ satisfying $dF=\Theta_P$, take $A=P^{-1}$ and $\eta(y,z)=-F(Ay,Az)$. For fixed $P$, the obstruction is $[\Theta_P]\in H^3(\bLie,\K)$, which corresponds under pullback to $[\mathscr T_A\Gamma]\in H^3(\bLie_A,\K)$. A diagonal lift requires $\Theta_P=0$ as a cochain.
\end{enumerate}
\end{proposition}
\begin{proof}
For \textnormal{(i)}, the base algebra is one-dimensional, so
$\alpha=\Gamma=\eta=\beta=0$ and $[\bLie,\bLie]=0$.
Consequently \eqref{eq:gd-inverse-der-alpha} and
\eqref{eq:gd-inverse-der-scalar} are automatic.  Since
$Pe=-\lambda e$ and $We=-c$, the mixed equation becomes $[D,\delta]=-\lambda\delta-\ad_c$, which is \eqref{eq:gd-bbm-compatibility}.  The inverse block
formula gives \eqref{eq:gd-bbm-double-derivation}. Solving
\eqref{eq:gd-bbm-parameters} for $a$ and $v$ gives the converse. For the reconstruction assertion, let $(\mathfrak d,T,\omega)$ be
quadratic symplectic, with associated skew derivation $E$ that is invertible,
and suppose that $I=\K e^*$ is a central line that is $E$-invariant.
Write $Ee^*=\lambda e^*$, where $\lambda\ne0$.
Skewness gives $2\lambda T(e^*,e^*)=0$, so $I$ is isotropic.
The converse in Proposition~\ref{prop:gd-invertible-bbm} gives an
adapted presentation with base $\bLie$ of dimension one.  In this
presentation $\alpha=\Gamma=0$, so it is an ordinary double extension
with $q=0$, and the preceding calculation gives exactly
\textnormal{(i)}.  If $\K$ is algebraically closed of characteristic zero and
$\mathfrak d\ne0$, then $\mathfrak d$ is nilpotent by
Remark~\ref{rem:JacobsonDerivation}.  Its center is nonzero and $E$-invariant and therefore contains an eigenline of $E$.
This proves the last assertion of \textnormal{(i)}.  Algebraic closure
is used only to obtain this line, not for either explicit extension
construction.

For \textnormal{(ii)}, first note that, for invertible $A$,
\[
\begin{aligned}
 M_A=0
 &\quad\Longleftrightarrow\quad
 P[y,z]_{\bLie}=[Py,z]_{\bLie}+[y,Pz]_{\bLie}\quad\Longleftrightarrow\quad P\in\Der(\bLie).
\end{aligned}
\]
Indeed, evaluate $M_A$ at $(Py,Pz)$ and apply $P$.
Under this condition,
\[
 [Py,Pz]_A=[y,Pz]_{\bLie}+[Py,z]_{\bLie}=P[y,z]_{\bLie},
\]
so $P:\bLie\to\bLie_A$ is a Lie algebra isomorphism. Naturality of the scalar differential gives $P^*d_A=dP^*$, and substitution yields $(P^*d_A\eta)(y,z,w)=-(dF)(y,z,w)$. Also, $AP=\id$ and cyclicity of $\theta$ give
\[
 (P^*\mathscr T_A\Gamma)(y,z,w)=\Theta_P(y,z,w)=(\mathscr I_P\Gamma)(y,z,w).
\]
Adding these identities proves \eqref{eq:gd-bbm-tstar-pullback}.
Since $P^*$ is invertible on cochains, the scalar equation in
\eqref{eq:tstar-ext} is therefore equivalent to $dF=\Theta_P$.
Together with the equivalence for $M_A=0$, this proves
\eqref{eq:gd-bbm-tstar-coboundary} directly from the Poisson criterion. Finally, $d\Gamma=0$ in a $T^*_{\theta}$-extension and the derivation
property of $P$ imply that $\Theta_P=\mathscr I_P\Gamma$ is a cocycle.
The cochain isomorphism $P^*$ consequently identifies the two
obstruction classes.  The inverse substitution
$\eta(y,z)=-F(Ay,Az)$ gives the stated converse, and $\eta=0$ is
equivalent to $F=0$, giving the diagonal assertion.
\end{proof}
We make the correspondence with the original formulas explicit.
\begin{remark}
For an ordinary double extension by a line, \eqref{eq:gd-bbm-parameters}, \eqref{eq:gd-bbm-compatibility}, and \eqref{eq:gd-bbm-double-derivation} agree with the formulas in \cite{BajoBenayadiMedina2007}, Lemma~4.1. Conversely, the data in that lemma give our operator through
	\[
	C=D^{-1},\qquad Ae=-\lambda^{-1}e,\qquad
	Ve=-\lambda^{-1}D^{-1}c,\qquad \eta=0.
	\]
The diagonal case is exactly $c=0$. Thus we recover the reconstruction in \cite{BajoBenayadiMedina2007}, Theorem~4.2. For the cohomological equation in the quadratic $T^*_{\theta}$ case, see \cite{BajoBenayadiMedina2007}, Proposition~2.2 and its proof.
\end{remark}

\section{Jacobi reconstruction along the isotropic flag}
\label{sec:jacobi-double}

We now extend the reconstruction of Poisson structures to Jacobi structures. Besides the skew operator, we must specify the Euler element. Its position in the isotropic flag leads to three cases, and the Poisson defect computed in the previous section gives the equations in each case.

Let $(\tilde{\mathfrak g},T)$ be the generalized double extension with data $(\psi,\alpha,\Gamma)$ satisfying Theorem~\ref{thm:gd-bracket}, bracket \eqref{eq:gd-bracket}, and form \eqref{eq:double-form}. Denote the Euler element by $E$. In an adapted splitting, write $E=y_E+x_E+f_E\in\bLie\oplus\g\oplus\bLie^*$, with $y_E\in\bLie$, $x_E\in\g$, and $f_E\in\bLie^*$. Define $x_E^\flat(x)=B(x_E,x)$. We distinguish the following positions of $E$ in the flag:
\[
\begin{aligned}
 \text{inner:}\quad &E\in \bLie^*
      &&\Longleftrightarrow\quad y_E=0,\ x_E=0,\\
 \text{middle:}\quad &E\in (\bLie^*)^\perp\setminus \bLie^*
      &&\Longleftrightarrow\quad y_E=0,\ x_E\neq0,\\
 \text{outer:}\quad &E\notin (\bLie^*)^\perp
      &&\Longleftrightarrow\quad y_E\neq0.
\end{aligned}
\]
In the \emph{inner} branch, $E$ lies in $\bLie^*$, so $E=f_E$. This includes $E=0$, which recovers the Poisson case. In the \emph{middle} branch, its image in $\widetilde{\g}/(\bLie^*)^\perp$ is zero, but its image in $(\bLie^*)^\perp/\bLie^*$ is nonzero. Here $E=x_E+f_E$, and $f_E$ need not vanish. In the \emph{outer} branch, the image of $E$ in $\widetilde{\g}/(\bLie^*)^\perp$ is nonzero. These cases are disjoint and exhaustive.
\subsection{The Jacobi system when the ideal is preserved}
\label{subsec:gd-jacobi-system}

Keep the operator \eqref{eq:gd-preserving-operator} and write $E=y_E+x_E+f_E$. To express the blocks of $\ad_E$ in the commutator equation, define
\begin{equation}\label{eq:gd-adjoint-data}
\begin{aligned}
   D_E&:=\ad_{y_E}, \quad      F_E=\psi_{y_E}+\ad_{x_E},\\
 W(y)&:=\alpha(y_E,y)-\psi_y x_E,\\
 q(y,z)&:=\Gamma(y_E,y,z)+B(x_E,\alpha(y,z))+f_E([y,z]).
\end{aligned}
\end{equation}
The form $q$ is alternating.  
We use the formula
$(\mathcal L_{y_E}\eta)(y,z)=-\eta([y_E,y],z)-\eta(y,[y_E,z])$.

\needspace{5\baselineskip}
The next theorem gives the full Jacobi system for an operator preserving $I$. It separates the first Jacobi identity into rank constraints and component equations, and then gives the commutator equation.

\begin{theorem}
\label{thm:gd-jacobi}
The pair $(\mathcal S,E)$ is Jacobi on a generalized double extension $\widetilde{\g}$ if and only if the following
conditions hold.  The rank constraints are
\begin{equation}\label{eq:gd-J0}
 \begin{gathered}
 f(y_E)Ay=f(Ay)y_E,\qquad f(y_E)Vy=f(Ay)x_E
       \quad(y\in\bLie,\ f\in \bLie^*), \quad y_E=0 \text{ or } C=0.
 \end{gathered}
\end{equation}
The remaining equations for the first Jacobi identity are
\begin{align}
 M_A(y,z)&=f_E(y)Az-f_E(z)Ay+\eta(y,z)y_E,
       \label{eq:gd-J1}\\
 M_C(x,w)&=B(x_E,x)Cw-B(x_E,w)Cx+B(Cx,w)x_E,
       \label{eq:gd-J2}\\
 K_y&=f_E(y)C+x_E\otimes(Vy)^\flat-Vy\otimes x_E^\flat,
       \label{eq:gd-J3}\\
 L^{\mathrm{gen}}(y,z)&=f_E(y)Vz-f_E(z)Vy+\eta(y,z)x_E,
       \label{eq:gd-J4}\\
 \widehat d_A\eta+\Xi_{A,V}&=f_E\wedge\eta.
       \label{eq:gd-J5}
\end{align}
The commutator equation is equivalent to
\begin{align}
 [D_E,A]&=0,\label{eq:gd-J6}\\
 [F_E,C]&=0,\label{eq:gd-J7}\\
 WA+F_EV-VD_E-CW&=0,\label{eq:gd-J8}\\
 (\mathcal L_{y_E}\eta+\mathscr T_Aq)(y,z)
       +B(Wy,Vz)-B(Wz,Vy)&=0.
       \label{eq:gd-J9}
\end{align}

\end{theorem}

\begin{proof}
We first compare the Poisson defect with the right-hand side of the Jacobi identity,
\begin{equation}\label{eq:gd-jacobi-operator-identity}
 M_{\mathcal S}(X,Z)
   =T(E,X)\mathcal SZ-T(E,Z)\mathcal SX+T(\mathcal SX,Z)E.
\end{equation}
which we denote by $R(X,Z)$. We compare the three components for each of the six pairs of summands. For $y,z\in\bLie$, they are
$$  R(y,z)_{\bLie}=f_E(y)Az-f_E(z)Ay+\eta(y,z)y_E,\qquad  R(y,z)_{\g}=f_E(y)Vz-f_E(z)Vy+\eta(y,z)x_E $$
$$ R(y,z)_{\bLie^*} =f_E(y)\eta^\flat z-f_E(z)\eta^\flat y+\eta(y,z)f_E.$$
%\begin{align*}
% R(y,z)_{\bLie}&=f_E(y)Az-f_E(z)Ay+\eta(y,z)y_E,\\
% R(y,z)_{\g}&=f_E(y)Vz-f_E(z)Vy+\eta(y,z)x_E,\\
% R(y,z)_{\bLie^*}&=f_E(y)\eta^\flat z-f_E(z)\eta^\flat y+\eta(y,z)f_E.
%\end{align*}
Comparison with \eqref{eq:gd-defect-bb} gives \eqref{eq:gd-J1} and \eqref{eq:gd-J4}. Evaluating the $\bLie^*$-component at $w$ gives
\[
 f_E(y)\eta(z,w)-f_E(z)\eta(y,w)+f_E(w)\eta(y,z)
                         =(f_E\wedge\eta)(y,z,w),
\]
which gives \eqref{eq:gd-J5}. For $y\in\bLie$, $x\in\g$, the components are
$$R(y,x)_{\bLie}=-B(x_E,x)Ay+B(Vy,x)y_E,\qquad  R(y,x)_{\g}=f_E(y)Cx-B(x_E,x)Vy+B(Vy,x)x_E$$
$$ R(y,x)_{\bLie^*}=-f_E(y)V^\dagger x-B(x_E,x)\eta^\flat y+B(Vy,x)f_E.$$
%\begin{align*}
% R(y,x)_{\bLie}&=-B(x_E,x)Ay+B(Vy,x)y_E,\\
% R(y,x)_{\g}&=f_E(y)Cx-B(x_E,x)Vy+B(Vy,x)x_E,\\
% R(y,x)_{\bLie^*}&=-f_E(y)V^\dagger x-B(x_E,x)\eta^\flat y+B(Vy,x)f_E.
%\end{align*}
The $\g$-component gives \eqref{eq:gd-J3}.
The $\bLie^*$-component, evaluated on $z$, is
\[
 -f_E(y)B(Vz,x)-B(x_E,x)\eta(y,z)+B(Vy,x)f_E(z),
\]
which equals $-B(L^{\mathrm{gen}}(y,z),x)$ by \eqref{eq:gd-J4}. The $\bLie$-component vanishes by the second constraint in \eqref{eq:gd-J0}. Indeed, pair that constraint with $x$ using $B$ and vary $f\in\bLie^*$ to obtain $B(Vy,x)y_E=B(x_E,x)Ay$. For $x,w\in\g$, we obtain
$$R(x,w)_{\bLie}=B(Cx,w)y_E,\qquad  R(x,w)_{\g}=B(x_E,x)Cw-B(x_E,w)Cx+B(Cx,w)x_E,$$
$$R(x,w)_{\bLie^*}=-B(x_E,x)V^\dagger w+B(x_E,w)V^\dagger x+B(Cx,w)f_E.$$
%\begin{align*}
% R(x,w)_{\bLie}&=B(Cx,w)y_E,\\
% R(x,w)_{\g}&=B(x_E,x)Cw-B(x_E,w)Cx+B(Cx,w)x_E,\\
% R(x,w)_{\bLie^*}&=-B(x_E,x)V^\dagger w+B(x_E,w)V^\dagger x+B(Cx,w)f_E.
%\end{align*}
The $\bLie$-component must vanish, which is precisely the last constraint in \eqref{eq:gd-J0}. The $\g$-component gives \eqref{eq:gd-J2}. Evaluating the $\bLie^*$-component at $y$ gives
\[
 -B(x_E,x)B(Vy,w)+B(x_E,w)B(Vy,x)+f_E(y)B(Cx,w)=B(K_yx,w),
\]
where the last equality is \eqref{eq:gd-J3} paired with $w$.
Thus it agrees with \eqref{eq:gd-defect-gg}. For $y\in\bLie$, $f\in \bLie^*$, the components are
$$ R(y,f)_{\bLie}=-f(y_E)Ay+f(Ay)y_E,\qquad R(y,f)_{\g}=-f(y_E)Vy+f(Ay)x_E, $$
$$R(y,f)_{\bLie^*}=-f_E(y)A^*f-f(y_E)\eta^\flat y+f(Ay)f_E.$$
%\begin{align*}
% R(y,f)_{\bLie}&=-f(y_E)Ay+f(Ay)y_E,\\
% R(y,f)_{\g}&=-f(y_E)Vy+f(Ay)x_E,\\
% R(y,f)_{\bLie^*}&=-f_E(y)A^*f-f(y_E)\eta^\flat y+f(Ay)f_E.
%\end{align*}
The $\bLie$- and $\g$-components of the defect are zero, so these give
the first two constraints in \eqref{eq:gd-J0}.  The $\bLie^*$-component
evaluated on $z$ is
\[
 -f_E(y)f(Az)-f(y_E)\eta(y,z)+f(Ay)f_E(z)=-f(M_A(y,z)),
\]
where \eqref{eq:gd-J1} gives the equality.  This is the component
in \eqref{eq:gd-defect-bI}. For $x\in\g$, $f\in \bLie^*$, the right side is
\[
 R(x,f)=-f(y_E)Cx-B(x_E,x)A^*f+f(y_E)V^\dagger x.
\]
Its $\g$-component vanishes by $y_E\otimes C=0$.
Its $\bLie^*$-component evaluated at $y$ is
$-B(x_E,x)f(Ay)+f(y_E)B(Vy,x)$, which vanishes by the second rank
constraint.  Finally,
\[
 R(f,k)=-f(y_E)A^*k+k(y_E)A^*f.
\]
Evaluation on $y$ gives $-f(y_E)k(Ay)+k(y_E)f(Ay)=0$ by the first
rank constraint.  This proves necessity and sufficiency for
the first Jacobi identity.  

For the second Jacobi identity, we first compute $\ad_E$. The bracket on the three summands gives
$$ [E,y]=[y_E,y]+\alpha(y_E,y)-\psi_y x_E+q^\flat y,\qquad  [E,x]=(\psi_{y_E}+\ad_{x_E})x-W^\dagger x,$$
$$[E,f]=\ad_{y_E}^*f=-D_E^*f, $$
%\begin{align*}
% [E,y]&=[y_E,y]+\alpha(y_E,y)-\psi_y x_E+q^\flat y,\\
% [E,x]&=(\psi_{y_E}+\ad_{x_E})x-W^\dagger x,\\
% [E,f]&=\ad_{y_E}^*f=-D_E^*f,
%\end{align*}
where $D_E^*$ denotes the ordinary dual. In $[E,y]$, the $\bLie^*$-components of $[y_E,y]$, $[x_E,y]$, and $[f_E,y]$, evaluated at $z$, are $\Gamma(y_E,y,z)$, $B(\alpha(y,z),x_E)$, and $f_E([y,z])$. Their sum is $q(y,z)$. In $[E,x]$, the $\bLie^*$-component is
\[
 -B(\alpha(y_E,y),x)+B(\psi_y x_E,x)=-B(Wy,x).
\]
Consequently
\begin{equation}\label{eq:gd-adE-matrix}
 \ad_E=
 \begin{pmatrix}
 D_E&0&0\\ W&F_E&0\\ q^\flat&-W^\dagger&-D_E^*
 \end{pmatrix}.
\end{equation}
To compute the commutator, temporarily write $D=D_E$ and $F=F_E$. Define $Z:=WA+FV-VD-CW$ and $R:=q^\flat A+A^*q^\flat-D^*\eta^\flat-\eta^\flat D-W^\dagger V+V^\dagger W$. Matrix multiplication gives
\[
[\ad_E,\mathcal S]=\ad_E\mathcal S-\mathcal S\ad_E=
\begin{pmatrix}
	[D,A]&0&0\\ Z&[F,C]&0\\ R&-Z^\dagger&-[D,A]^*
\end{pmatrix}.
\]
%\[
% \ad_E\mathcal S=
% \begin{pmatrix}
% DA&0&0\\
% WA+FV&FC&0\\
% q^\flat A-W^\dagger V-D^*\eta^\flat&-W^\dagger C+D^*V^\dagger&D^*A^*
% \end{pmatrix}
%\]
%and
%\[
% \mathcal S\ad_E=
% \begin{pmatrix}
% AD&0&0\\
% VD+CW&CF&0\\
% \eta^\flat D-V^\dagger W-A^*q^\flat&-V^\dagger F+A^*W^\dagger&A^*D^*
% \end{pmatrix}.
%\]
For instance, $F$ and $C$ are $B$-skew, so
$Z^\dagger=A^*W^\dagger-V^\dagger F-D^*V^\dagger+W^\dagger C$,
which verifies the block from $\g$ to $\bLie^*$, including its sign.
For $y,z\in\bLie$, the block from $\bLie$ to $\bLie^*$ evaluates to
\begin{align*}
 (Ry)(z)
 ={}&q(Ay,z)+q(y,Az)-\eta(y,[y_E,z])-\eta([y_E,y],z)\\
 &-B(Wz,Vy)+B(Vz,Wy)\\
 ={}&(\mathscr T_Aq+\mathcal L_{y_E}\eta)(y,z)
                         +B(Wy,Vz)-B(Wz,Vy).
\end{align*}
Its vanishing is therefore \eqref{eq:gd-J9}. The other independent blocks give \eqref{eq:gd-J6}--\eqref{eq:gd-J8}, and their dual blocks vanish at the same time.
\end{proof}

The same substitutions as in Corollary~\ref{cor:ordi_Poisson} give the Jacobi equations for ordinary double extensions and $T^\ast_\theta$-extensions. Likewise, the argument of Corollary~\ref{cor:Sol_Nil_Poisson} gives the corresponding specializations and reconstruction results for nilpotent or solvable quadratic Jacobi Lie algebras, and for those that are decomposable and neither simple nor one-dimensional. We leave these substitutions to the reader.

\subsection{Cohomological interpretation}\label{scs:cohomological-interpretation-Jacobi}
We now interpret \eqref{eq:gd-J5} cohomologically in the inner branch, where $E=f_E\in\bLie^*$. The argument follows the Poisson case, but uses the differential with coefficients determined by $f_E$. The equations independent of $\eta$ become
\begin{align}
 M_A(y,z)&=f_E(y)Az-f_E(z)Ay,& f_E(b_A(y,z))&=0,
       \label{eq:gd-V12}\\
 M_C&=0,& K_y&=f_E(y)C,
       \label{eq:gd-V34}\\
 L^{\mathrm{gen}}(y,z)&=f_E(y)Vz-f_E(z)Vy.
       \label{eq:gd-V5}
\end{align}
The remaining equation is
$\widehat d_A\eta-f_E\wedge\eta+\Xi_{A,V}=0$.
By Lemma~\ref{lem:jacobi-descendent}, the first line makes
\begin{equation}\label{eq:gd-inner-descendent}
 [y,z]_{A,f_E}=b_A(y,z)+f_E(y)z-f_E(z)y
\end{equation}
a Lie bracket with character $f_E$. Write $\bLie_{A,f_E}$ for this Lie algebra and $\K_{f_E}$ for its one-dimensional module with action $y\cdot c=f_E(y)c$. The differential from Lemma~\ref{lem:jacobi-descendent} is
\begin{equation}\label{eq:gd-twisted-differential}
 \partial_{A,f_E}=d_{A,f_E}+f_E\wedge,
 \qquad
 \partial_{A,f_E}\omega=\widehat d_A\omega-(k-1)f_E\wedge\omega
       \quad(\omega\in\Lambda^k\bLie^*).
\end{equation}
In particular, on two-forms it is $\widehat d_A-f_E\wedge$ and on
three-forms it is $\widehat d_A-2 f_E\wedge$.

\needspace{5\baselineskip}
The next theorem identifies the obstruction using this differential.

\begin{theorem}
\label{thm:gd-inner-obstruction}
Fix $A,C,V,f_E$ satisfying \eqref{eq:gd-V12}--\eqref{eq:gd-V5}.
Then $\partial_{A,f_E}\Xi_{A,V}=0$.
A two-form $\eta$ completes these data to a Jacobi structure if and
only if
\begin{equation}\label{eq:gd-inner-primitive}
 \partial_{A,f_E}\eta+\Xi_{A,V}=0.
\end{equation}
Equivalently, the obstruction is
\begin{equation}\label{eq:gd-inner-H3}
 [\Xi_{A,V}]=0\quad\text{in }H^3(\bLie_{A,f_E},\K_{f_E}).
\end{equation}
In a fixed splitting, the solution space is
$\eta_0+Z^2(\bLie_{A,f_E},\K_{f_E})$ whenever one solution $\eta_0$ exists.
\end{theorem}

\begin{proof}
Set $y_E=x_E=0$ in the full system. Then $D_E=F_E=W=0$ and $q(y,z)=f_E([y,z])=-d f_E(y,z)$. The commutator $[\ad_{f_E},\mathcal S]=\ad_{f_E}\circ\mathcal S-\mathcal S\circ\ad_{f_E}$ has only one possibly nonzero block, from $\bLie$ to $I$, namely $q^\flat A+A^*q^\flat$. Evaluating it gives
\[
 \bigl((q^\flat A+A^*q^\flat)y\bigr)(z)
 =q(Ay,z)+q(y,Az)=f_E(b_A(y,z)).
\]
Thus $[\ad_{f_E},\mathcal S]=0$ is the second equation in \eqref{eq:gd-V12} and does not involve $\eta$. The rank constraints hold automatically, and \eqref{eq:gd-J1}--\eqref{eq:gd-J4} give the remaining equations in \eqref{eq:gd-V12}--\eqref{eq:gd-V5}. The scalar condition is $\widehat d_A\eta-f_E\wedge\eta+\Xi_{A,V}=0$. Since $\widehat d_A-f_E\wedge=\partial_{A,f_E}$ on two-forms, this is \eqref{eq:gd-inner-primitive}. We have therefore proved the reconstruction criterion without yet using closedness of $\Xi_{A,V}$.

To prove closedness, take $\eta=0$ and define $r=\pi_{\mathcal S}$ without imposing the scalar equation. The commutator still vanishes. More explicitly, Theorem~\ref{thm:operator-jacobi}, with $X^\flat=T(X,\cdot)$, gives $[f_E,r](X^\flat,Y^\flat)=T\bigl([\ad_{f_E},\mathcal S]X,Y\bigr)=0$. Thus the second equation in \eqref{eq:operator-jacobi-algebraic} holds for this choice of $\eta$. The failure of the first Jacobi equation is measured by $\mathcal F=\tfrac12[r,r]-f_E\wedge r$. Define the corresponding operator defect, as in \eqref{eq:J1}, by
\[\mathcal D(X,Y)
 =:M_{\mathcal S}(X,Y)-T(f_E,X)\mathcal SY
       +T(f_E,Y)\mathcal SX-T(\mathcal SX,Y)f_E.\]
The pairing identity in Theorem~\ref{thm:operator-jacobi} gives $\mathcal F(X^\flat,Y^\flat,Z^\flat)=T\bigl(\mathcal D(X,Y),Z\bigr)$. The equations for $M_A$ and $M_C$ in \eqref{eq:gd-V12}--\eqref{eq:gd-V34} concern the induced operators on $\widetilde{\g}/I^\perp\simeq\bLie$ and $I^\perp/I\simeq\g$. We call them the quotient equations. The equations for $K_y$ and $L^{\mathrm{gen}}$ in \eqref{eq:gd-V34}--\eqref{eq:gd-V5} concern mixed components. Together they cancel all components of $\mathcal D$ except the scalar one. To check this, use $\widetilde{\g}=\bLie\oplus\g\oplus I$ and take $y,z,t\in\bLie$, $x,w,u\in\g$, and $\lambda\in I$. Lemma~\ref{lem:gd-defect} and the comparison in the proof of Theorem~\ref{thm:gd-jacobi} give
\begin{align*}
 T(\mathcal D(y,z),\lambda)
   &=\lambda\bigl(M_A(y,z)-f_E(y)Az+f_E(z)Ay\bigr)=0,\\
 T(\mathcal D(y,z),x)
   &=B\bigl(L^{\mathrm{gen}}(y,z)-f_E(y)Vz+f_E(z)Vy,x\bigr)=0,\\
 T(\mathcal D(y,x),w)
   &=B\bigl((K_y-f_E(y)C)x,w\bigr)=0,\\
 T(\mathcal D(x,w),u)
   &=B\bigl(M_C(x,w),u\bigr)=0,\\
 T(\mathcal D(y,z),t)
   &=\Xi_{A,V}(y,z,t).
\end{align*}
The first four expressions vanish by the quotient and mixed equations, while the last uses $\eta=0$. Every other type of triple vanishes by the block structure, up to permutation. Since the pairing is alternating and trilinear, this determines it on all arguments.

Let $p:\widetilde{\g}\to\bLie$ be the quotient projection.  We have
proved
\[\mathcal F(X^\flat,Y^\flat,Z^\flat)
   =\Xi_{A,V}(pX,pY,pZ),
 \quad (\Lambda^3T^\flat)\mathcal F=p^*\Xi_{A,V}.
\]
The corresponding three-form vanishes whenever an argument belongs to $\ker p=I^\perp$. Since $T^\flat(I)=\Ann(I^\perp)$, the metric identification used in Theorem~\ref{thm:gd-poisson-obstruction} gives $\mathcal F=\Xi_{A,V}$ in $\Lambda^3 I=\Lambda^3\bLie^*\subseteq\Lambda^3\widetilde{\g}$. Here the scalar three-cochain on $\bLie$ is viewed as a trivector in $I=\bLie^*$. We can now apply the Schouten derivation rule. Since $r$ has degree two and $f_E$ has degree one,
\[
 [r,f_E\wedge r]=[r,f_E]\wedge r-f_E\wedge[r,r].
\]
Graded antisymmetry gives $[r,f_E]=-[f_E,r]=0$.
The graded Jacobi identity gives $[r,[r,r]]=0$.  Combining these
two identities yields
\[
 [r,\mathcal F]
 =\tfrac12[r,[r,r]]-[r,f_E\wedge r]
 =f_E\wedge[r,r]
   =2 f_E\wedge(\mathcal F+f_E\wedge r)=2 f_E\wedge\mathcal F,
\]
where the last equality uses $f_E\wedge f_E=0.$ By Proposition~\ref{prop:gd-bivector}, the left-hand side corresponds to $\widehat d_A\Xi_{A,V}$. Thus $\widehat d_A\Xi_{A,V}-2 f_E\wedge\Xi_{A,V}=0$. On three-forms, the differential in \eqref{eq:gd-twisted-differential} is $\widehat d_A-2 f_E\wedge$, so $\partial_{A,f_E}\Xi_{A,V}=0$.

By Lemma~\ref{lem:jacobi-descendent}, the first line of \eqref{eq:gd-V12} makes \eqref{eq:gd-inner-descendent} a Lie bracket and $f_E$ a character. Therefore $\partial_{A,f_E}$ is a differential with values in $\K_{f_E}$. The class of $\Xi_{A,V}$ vanishes exactly when a two-form has differential $-\Xi_{A,V}$, proving \eqref{eq:gd-inner-H3}. If $\eta_0$ is one such form, then
\[
 \partial_{A,f_E}\eta+\Xi_{A,V}=0
 \quad\Longleftrightarrow\quad
 \partial_{A,f_E}(\eta-\eta_0)=0.
\]
Thus the solutions in the fixed splitting form the stated translate of $Z^2$. For $f_E=0$, the argument reduces to the Poisson defect and its ordinary differential. For an ordinary extension, with $q=0,$ one has $\Xi_{A,V}=\kappa_V$ throughout.
\end{proof}

\needspace{5\baselineskip}
When the three blocks $A,C,V$ vanish, the criterion takes a particularly simple form.

\begin{corollary}
\label{cor:gd-pure-lower}
For $A=C=V=0$, every $\eta\in\Lambda^2\bLie^*$ defines a Poisson
operator on every generalized double extension.
The pair with $E=f_E\in \bLie^*$ is Jacobi exactly when $f_E\wedge\eta=0$.
\end{corollary}

\section{Examples}
\label{sec:poisson-jacobi-examples}

We illustrate the reconstruction criteria on oscillator algebras and the five-dimensional algebra $A_{5,3}$, which underlies a Wess--Zumino--Witten model. We then examine Poisson structures defined by curvature operators of a bi-invariant metric.

Oscillator algebras are the only nonabelian solvable Lie algebras admitting an irreducible bi-invariant Lorentzian metric \cite{Medina1985}. Their relation to the harmonic oscillator and their uses in mathematical physics make them natural examples here \cite{GadeaOubina1999,Levichev1986,Levichev1990,MullerRicci1988,NappiWitten1993,Streater1967}. The algebra $A_{5,3}$ is the unique indecomposable real Lie algebra of dimension five admitting an invariant metric. Its Wess--Zumino--Witten model has gauged versions giving four-dimensional plane-wave string backgrounds \cite{Kehagias1994}. These examples connect the extension constructions with differential geometry and mathematical physics.

In the first two Jacobi examples, the Euler element $E$ belongs to the isotropic ideal $I$ used in the corresponding Poisson example. Over $\mathbb R$, left translation gives its Euler vector field. We do not assume that the operator preserves $I$. For oscillator algebras, preservation follows from the Jacobi equations. For $A_{5,3}$, we also obtain a family whose operators do not preserve $I$.

\subsection{The oscillator algebra}
\label{sec:oscillator-poisson}

Work over $\mathbb R$. Let $(\mathfrak v,B_{\mathfrak v})$ be a Euclidean vector space of dimension $2n$, with $n\geq1$, and let $\delta\in\mathfrak o(\mathfrak v,B_{\mathfrak v})$ be invertible. Choose an orthonormal basis $p_1,q_1,\ldots,p_n,q_n$ such that
\begin{equation}\label{eq:osc-frequencies}
 \delta p_i=\lambda_iq_i,\qquad
 \delta q_i=-\lambda_ip_i,\qquad
 \lambda_i>0,\quad 1\leq i\leq n.
\end{equation}
Define $\mathfrak{osc}_{\delta}=\mathbb Re\oplus \mathfrak v\oplus\mathbb Rf$ with Lie bracket
\begin{equation}\label{eq:osc-bracket}
 [e,x]=\delta x,\qquad
 [x,z]=B_{\mathfrak v}(\delta x,z)f,\qquad
 [f,\mathfrak{osc}_{\delta}]=0,\qquad x,z\in \mathfrak v.
\end{equation}
The quadratic form is
\begin{equation}\label{eq:osc-metric}
 \begin{gathered}
 T|_{\mathfrak v\times \mathfrak v}=B_{\mathfrak v},\qquad T(e,f)=1,\\
 T(e,e)=T(f,f)=0,\qquad T(e,\mathfrak v)=T(f,\mathfrak v)=0.
 \end{gathered}
\end{equation}
This is the ordinary double extension of the abelian quadratic algebra $(\mathfrak v,B_{\mathfrak v})$ by the one-dimensional algebra $\bLie=\mathbb Re$, with
\begin{equation}\label{eq:osc-extension-data}
 \psi_{te}=t\delta,\qquad \alpha=\Gamma=0,\qquad I=\mathbb Rf.
\end{equation}
In particular, $T$ has signature $(2n+1,1)$.  Invertibility of $\delta$
gives
\[
 Z(\mathfrak{osc}_{\delta})=I,\qquad
 I^\perp=\mathfrak v\oplus\mathbb Rf=[\mathfrak{osc}_{\delta},\mathfrak{osc}_{\delta}],
 \qquad [I^\perp,I^\perp]=I.
\]
Thus the algebra is solvable of derived length three, but not nilpotent, since $\ad_e|_{\mathfrak v}=\delta$ is invertible. The middle algebra $\mathfrak v$ is abelian in the extension data, although its image in the chosen splitting of $\mathfrak{osc}_{\delta}$ is not a Lie subalgebra. The alternating form
\begin{equation}\label{eq:osc-omega}
 \omega_\delta(x,z):=B_{\mathfrak v}(\delta x,z)
\end{equation}
is symplectic.  Thus $I^\perp$ is a Heisenberg algebra and
$\mathfrak{osc}_{\delta}$ is its semidirect extension by $e$.

We first show that the equations force preservation of $I$ in both cases considered below.

\begin{lemma}\label{lem:osc-preservation}
Every Jacobi structure $(\mathcal S,E)$ on
$(\mathfrak{osc}_{\delta},T)$ with $E\in I$ satisfies
$\mathcal S(I)\subseteq I$.
In particular, every Poisson structure preserves $I$.
\end{lemma}

\begin{proof}
Every $T$-skew operator has the form
\begin{equation}\label{eq:osc-general-skew}
 \mathcal S=
 \begin{pmatrix}
 a&-u^\flat&0\\
 v&C&u\\
 0&-v^\flat&-a
 \end{pmatrix},\qquad
 a\in\mathbb R,\quad u,v\in\mathfrak v,\quad
 C\in\mathfrak o(\mathfrak v,B_{\mathfrak v}),
\end{equation}
where $u^\flat(x)=B_{\mathfrak v}(u,x)$ and similarly for $v^\flat$.
Write $E=bf$.
For $x\in\mathfrak v$, the middle component of the Poisson defect is
\[
 \operatorname{pr}_{\mathfrak v}M_{\mathcal S}(x,f)
 =-B_{\mathfrak v}(u,x)\delta u
       +B_{\mathfrak v}(\delta u,x)u.
\]
On $(x,f)$, the right-hand side of \eqref{eq:gd-jacobi-operator-identity} lies in $I$, since $T(E,x)=T(E,f)=0$. Taking $x=u$ and using skewness of $\delta$ gives $B_{\mathfrak v}(u,u)\delta u=0$. Positive definiteness and invertibility of $\delta$ imply $u=0$. Thus $\mathcal S f=-af$, so $\mathcal S$ preserves $I$. Taking $E=0$ proves the Poisson assertion.
\end{proof}

\subsubsection{Poisson structures} The reconstruction criterion now describes all Poisson operators on this algebra
\begin{theorem}\label{thm:osc-classification}
The Poisson structures on $(\mathfrak{osc}_{\delta},T)$ are exactly the operators
\begin{equation}\label{eq:osc-all-S}
 \mathcal S_{a,v,C}=
 \begin{pmatrix}
 a&0&0\\
 v&C&0\\
 0&-v^\flat&-a
 \end{pmatrix},\qquad
 a\in\mathbb R,\quad v\in \mathfrak v,\quad C\in\mathfrak o(\mathfrak v,B_{\mathfrak v}),
\end{equation}
subject to
\begin{equation}\label{eq:osc-final-conditions}
  aC=0,\qquad C\delta C=0.
\end{equation}
There is no condition on $v$.  Their bivectors, in the convention
\eqref{eq:rA-def}, are
\begin{equation}\label{eq:osc-all-r}
 r_{\mathcal S}=r_C+f\wedge(ae+v),\qquad
 r_C(x^\flat,z^\flat)=B_{\mathfrak v}(Cx,z),\qquad r_C\in\Lambda^2\mathfrak v,
\end{equation}
where $x^\flat=B_{\mathfrak v}(x,\cdot)$.
\end{theorem}

\begin{proof}
By Lemma~\ref{lem:osc-preservation}, every Poisson operator has the form \eqref{eq:osc-all-S}. In Theorem~\ref{thm:gd-poisson}, take $A=a\,\id_{\bLie}$ and $V(e)=v$. Both $\bLie$ and the middle algebra $\mathfrak v$ are abelian, so $M_A=M_C=0$. Since $\Lambda^2\bLie=\Lambda^3\bLie=0$, the equations $L^{\mathrm{gen}}=0$ and $d_A\eta+\Xi_{A,V}=0$ hold automatically. The term $\ad_{Vy}$ in \eqref{eq:gd-K} is computed in the abelian middle algebra, so it also vanishes. The only remaining condition is
\begin{equation}\label{eq:osc-K}
 K_e=a[\delta,C]-C\delta C=0.
\end{equation}

For $a=0$, equation~\eqref{eq:osc-K} is $C\delta C=0$. Suppose $a\neq0$ and define $F=C/a$. Dividing by $a^2$ gives $[\delta,F]=F\delta F$, or equivalently
\begin{equation}\label{eq:osc-determinant-identity}
 (\id_{\mathfrak v}-F)\delta(\id_{\mathfrak v}+F)=\delta.
\end{equation}
In an orthonormal basis, $F^\top=-F$. Since $\det\delta\neq0$, taking determinants gives $\det(\id_{\mathfrak v}+F)^2=1$. The eigenvalues of a real skew-symmetric matrix occur in pairs $\pm i\mu_j$, with $\mu_j\in\mathbb R$, so
\[
 \det(\id_{\mathfrak v}+F)=\prod_{j=1}^n(1+\mu_j^2)\geq1.
\]
Equality forces every $\mu_j=0$. Since a real skew-symmetric matrix is normal and therefore diagonalizable over $\mathbb C$, we obtain $F=0$ and $C=0$. Thus \eqref{eq:osc-K} is equivalent to \eqref{eq:osc-final-conditions}. Conversely, these conditions satisfy every equation in Theorem~\ref{thm:gd-poisson}. The bivector formula follows from Proposition~\ref{prop:gd-bivector}.
\end{proof}

Consequently, the solution set is the union of the two families
\begin{equation}\label{eq:osc-two-families}
 \begin{aligned}
 &r=f\wedge(ae+v), &&a\in\mathbb R,\quad v\in \mathfrak v.\\
 &r=r_C+f\wedge v, &&v\in \mathfrak v,\quad C^\top=-C,\quad C\delta C=0.
 \end{aligned}
\end{equation}
The intersection is $\{f\wedge v:v\in\mathfrak v\}$. Thus $A=0$ is necessary when $C\neq0$, but not for a Poisson structure in general. For example, $f\wedge e$ has $A\neq0$. Equation~\eqref{eq:osc-K} is the operator form of the tensor criterion in \cite{BoucettaMedina2010}, Theorem~1.2 part~(2). The derivation here applies to all positive frequencies, including repeated or resonant ones.

To describe the equation $C\delta C=0$, it is enough to examine the image of $C$ (see the discussion after Theorem~1.2 in \cite{BoucettaMedina2010}).

\begin{proposition}\label{prop:osc-isotropic-image}
For $C\in\mathfrak o(\mathfrak v,B_{\mathfrak v})$,
\begin{equation}\label{eq:osc-isotropic-image}
 C\delta C=0
 \quad\Longleftrightarrow\quad
 \Imop C\text{ is }\omega_\delta\text{-isotropic}.
\end{equation}
Every nonzero solution is obtained by choosing an even-dimensional $\omega_\delta$-isotropic subspace $W\subset\mathfrak v$, taking an invertible $B_{\mathfrak v}|_W$-skew map $C_W:W\to W$, and extending it by zero on $W^\perp$. In particular, the possible ranks are
\begin{equation}\label{eq:osc-C-ranks}
 \rank C=2k,\qquad 0\leq k\leq\left\lfloor\frac n2\right\rfloor.
\end{equation}
\end{proposition}

\begin{proof}
For $x,z\in \mathfrak v$, skewness gives
\[
 B_{\mathfrak v}(C\delta Cx,z)=-B_{\mathfrak v}(\delta Cx,Cz)
                 =-\omega_\delta(Cx,Cz).
\]
Nondegeneracy of $B_{\mathfrak v}$ proves \eqref{eq:osc-isotropic-image}.
If $W=\Imop C$, then skewness also gives $\Ker C=W^\perp$.
Positive definiteness implies $W\cap W^\perp=0$, so $C|_W$ is invertible.
Conversely, a skew map on $W$ extended by zero on $W^\perp$ is skew on
$\mathfrak v$, and isotropy of its image gives $C\delta C=0$.  The rank is even
because the form $B_{\mathfrak v}(C\cdot,\cdot)$ is alternating, and at most $n$
because an isotropic subspace of a symplectic space of dimension $2n$
has dimension at most $n$.  Every even rank in the stated range occurs.
\end{proof}

Equivalently, choose an $\omega_\delta$-Lagrangian subspace $L\subset\mathfrak v$ and $\rho\in\Lambda^2L$. Then $r=\rho+f\wedge v$ is Poisson for every $v\in\mathfrak v$, and every member of the $a=0$ family is obtained this way. The vector $v$ need not lie in $L$. Indeed, $[L,L]=0$ and $[v,L]\subseteq\mathbb Rf$, so wedging with $f$ makes the mixed Schouten bracket vanish.

\subsubsection{Jacobi structures}
\label{subsec:osc-jacobi}

Let the Euler element belong to $I=\mathbb Rf$ and write $E=bf$, with $b\in\mathbb R$. Since $f$ is central, $[E,r]=0$ for every bivector $r$. We therefore need only the first Jacobi equation.

\begin{theorem}\label{thm:osc-jacobi-equation}
All Jacobi structures on $(\mathfrak{osc}_{\delta},T)$ with $E\in I$
are the pairs $(\mathcal S_{a,v,C},bf)$, where
$\mathcal S_{a,v,C}$ is given by \eqref{eq:osc-all-S} and
\begin{equation}\label{eq:osc-jacobi-equation}
 a[\delta,C]-C\delta C=bC.
\end{equation}
The vector $v\in\mathfrak v$ is arbitrary, and the bivector is
\begin{equation}\label{eq:osc-jacobi-r}
 r=r_C+f\wedge(ae+v).
\end{equation}
\end{theorem}

\begin{proof}
Lemma~\ref{lem:osc-preservation} gives preservation of $I$ without an additional assumption. In \eqref{eq:gd-V12}--\eqref{eq:gd-V5}, take $A=a\,\id_{\bLie}$, $V(e)=v$, and $f_E(e)=b$. Both quotient algebras are abelian, the adjoint term in $K_e$ vanishes, and $\dim\bLie=1$. The sole remaining equation is $K_e=bC$, namely \eqref{eq:osc-jacobi-equation}. Centrality of $E$ gives the commutator equation, and \eqref{eq:osc-all-r} gives the bivector in \eqref{eq:osc-jacobi-r}.
\end{proof}
For $a\neq0$, the equation also forces $C$ to commute with $\delta$, which gives the following description of its solutions.
\begin{proposition}\label{prop:osc-jacobi-commutation}
Every solution of \eqref{eq:osc-jacobi-equation} with $a\neq0$ satisfies
\begin{equation}\label{eq:osc-jacobi-commutation}
 [\delta,C]=0.
\end{equation}
If $b=0$, then $C=0$.  If $b\neq0$, these solutions are exactly
\begin{equation}\label{eq:osc-jacobi-general-nonzero-a}
 C=-b\,\delta^{-1}P_W,
\end{equation}
where $W\subseteq\mathfrak v$ is a $\delta$-invariant subspace and $P_W$ is its $B_{\mathfrak v}$-orthogonal projection. Taking $W=0$ gives $C=0$.
\end{proposition}

\begin{proof}
To prove commutation, use the positive definite inner product $\langle P,Q\rangle_{\mathrm{HS}}:=\operatorname{tr}(P^{\top}Q)$ on $\End(\mathfrak v)$, where the transpose is taken with respect to $B_{\mathfrak v}$. Define $D=[\delta,C]$. The operators $C$, $D$, and $C\delta C$ are skew, and cyclicity of the trace gives
\begin{align*}
 \langle C,D\rangle_{\mathrm{HS}}
   &=-\operatorname{tr}(C\delta C-C^2\delta)=0,\\
 \langle C\delta C,D\rangle_{\mathrm{HS}}
   &=-\operatorname{tr}(C\delta C\delta C)
       +\operatorname{tr}(C\delta C^2\delta)=0.
\end{align*}
For the second equality, both traces on the right are
$\operatorname{tr}(\delta C\delta C^2)$.
Taking the inner product of \eqref{eq:osc-jacobi-equation} with $D$ gives $a\,\|[\delta,C]\|_{\mathrm{HS}}^2=0$. Since $a\neq0$, positive definiteness proves \eqref{eq:osc-jacobi-commutation}. The equation then reduces to $C\delta C=-bC$. If $b=0$, commutation and invertibility of $\delta$ give $C^2=0$. Skewness yields $\|Cx\|^2=-B_{\mathfrak v}(C^2x,x)=0$, so $C=0$.
Suppose $b\neq0$ and define $W=\Imop C$.
Then $\Ker C=W^\perp$ and $C|_W$ is invertible.
Commutation makes $W$ and $W^\perp$ invariant under $\delta$.
Restricting $C\delta C=-bC$ to $W$ gives
$C|_W=-b(\delta|_W)^{-1}$, proving
\eqref{eq:osc-jacobi-general-nonzero-a}.
Conversely, the orthogonal projection onto a $\delta$-invariant subspace commutes with $\delta$. Thus \eqref{eq:osc-jacobi-general-nonzero-a} defines a skew map satisfying the equation.
\end{proof}

\subsection{The algebra \texorpdfstring{$A_{5,3}$}{A5,3}}
\label{sec:a53-poisson}

Let $\mathfrak d=A_{5,3}=\Span\{y_1,y_2,x,f_1,f_2\}$ be the five-dimensional Lie algebra with nonzero brackets
\begin{equation}\label{eq:a53-brackets}
 [y_1,y_2]=x,\qquad [y_1,x]=-f_2,\qquad [y_2,x]=f_1.
\end{equation}
This is the presentation in \cite{Kehagias1994}, equation~(16), after the change of basis
\[
 J_1=-f_1,\qquad J_2=f_2,\qquad J_3=x,\qquad
 J_4=y_1,\qquad J_5=y_2.
\]
Equip $\mathfrak d$ with the form determined by symmetry and the equations
\begin{equation}\label{eq:a53-metric}
 T(y_i,f_j)=\delta_{ij},\qquad T(x,x)=1,
 \qquad\text{all other basis pairings zero.}
\end{equation}
The form is nondegenerate and invariant. Indeed, $T([\cdot,\cdot],\cdot)$ is alternating, and its only nonzero value on basis vectors, up to permutation, is $T([y_1,y_2],x)=1$. Over $\mathbb R$, the signature is $(3,2)$. The form $T$ and the Cartan three-form $\mathcal H(X,Y,Z)=T([X,Y],Z)$ give the metric and three-form of the Wess--Zumino--Witten model on the integrating nilpotent group. This is the $A_{5,3}$ WZW model of \cite{Kehagias1994}, Section~3, with central charge $c=5$.
The center and its orthogonal complement are
\begin{equation}\label{eq:a53-flag}
 \begin{gathered}
 I:=Z(\mathfrak d)=\Span\{f_1,f_2\},\qquad
 I^\perp=\K x\oplus I=[\mathfrak d,\mathfrak d],\\
 [\mathfrak d,[\mathfrak d,\mathfrak d]]=I,\qquad
 [\mathfrak d,I]=0.
 \end{gathered}
\end{equation}
Thus $\mathfrak d$ is three-step nilpotent, $I$ is an
isotropic ideal of dimension two, and $I^\perp$ is abelian.

The basis displayed above gives an adapted Witt splitting
$\mathfrak d=\bLie\oplus\g\oplus I$, where
\begin{equation}\label{eq:a53-quotients}
 \bLie=\mathfrak d/I^\perp\cong\Span\{y_1,y_2\},\qquad
 (\g,B)=(I^\perp/I,B)\cong(\K x,B),\qquad B(x,x)=1.
\end{equation}
Both quotient algebras are abelian, and $f_1,f_2$ are the basis of
$I=\bLie^*$ dual to $y_1,y_2$.  Their chosen lifts $y_1,y_2$ do
not span a Lie subalgebra of $\mathfrak d$.

\begin{proposition}\label{prop:a53-generalized}
Relative to the marked ideal $I=Z(\mathfrak d)$, the quadratic algebra
\eqref{eq:a53-brackets}--\eqref{eq:a53-metric} is the generalized
double extension with data
\begin{equation}\label{eq:a53-data}
 \psi=0,\qquad \alpha(y_1,y_2)=x,\qquad \Gamma=0.
\end{equation}
There is no ordinary double-extension presentation relative to the same ideal, even when the additional invariant form $q$ in \eqref{eq:gd-ordinary-extra-form} is allowed.
\end{proposition}

\begin{proof}
Since $\g$ is abelian, $\ad_{\alpha(y,z)}=0$, so \eqref{eq:gd-curvature} holds with $\psi=0$. Every term in \eqref{eq:gd-bianchi} vanishes. Equation~\eqref{eq:gd-gamma} also holds, since it is an equality of four-forms on the two-dimensional algebra $\bLie$. The bracket \eqref{eq:gd-bracket} gives
\[
 [y_1,y_2]=\alpha(y_1,y_2)=x,\quad
 [y_1,x]=-B(\alpha(y_1,\cdot),x)=-f_2,\quad
 [y_2,x]=-B(\alpha(y_2,\cdot),x)=f_1.
\]
All other brackets vanish, and \eqref{eq:double-form} agrees with \eqref{eq:a53-metric}. To rule out an ordinary presentation, consider any linear section of $\mathfrak d\to\bLie$. Its lifts have the form
\[
 \ell_1=y_1+s_1x+z_1,\qquad
 \ell_2=y_2+s_2x+z_2,\qquad s_i\in\K,\quad z_i\in I.
\]
Their bracket is $[\ell_1,\ell_2]=x-s_1f_1-s_2f_2\ne0$. Since $\bLie$ is abelian, no such section is a Lie algebra homomorphism. Theorem~\ref{thm:gd-recognition} therefore excludes an ordinary double extension relative to $I$.

The change-of-splitting formulas in Subsection~\ref{subsec:gd-splittings} give another proof. Since $\psi=0$ and both $\bLie$ and $\g$ are abelian, \eqref{eq:gd-splitting-alpha} gives $\alpha'=\alpha\ne0$ for every adapted splitting keeping $I=Z(\mathfrak d)$. Also, $\Gamma'=0$ because $\dim\bLie=2$. Thus no ordinary double extension relative to this ideal is possible, even with the additional form $q$.
\end{proof}

This does not exclude an ordinary double extension relative to a smaller isotropic ideal.

\subsubsection{Poisson structures} The Poisson equations force the operator to preserve the center.

\begin{lemma}\label{lem:a53-preservation}
Every Poisson structure on $(\mathfrak d,T)$ preserves
$I=Z(\mathfrak d)$, and consequently also $I^\perp$.
\end{lemma}

\begin{proof}
Begin with an arbitrary $T$-skew operator in the six-block form \eqref{eq:gd-six-blocks}, without assuming that it preserves $I$. Its block $C$ vanishes because $\dim\g=1$. Write
\[
 Qf_1=-\kappa y_2,\qquad Qf_2=\kappa y_1,\qquad
 Uf_1=p x,\qquad Uf_2=q x,
 \qquad \kappa,p,q\in\K.
\]
Let $\operatorname{pr}_{\g}$ denote projection onto the middle
summand $\K x$.  Centrality of $I$ and \eqref{eq:a53-brackets} give
\[
 \operatorname{pr}_{\g}M_{\mathcal S}(f_1,f_2)=\kappa^2x.
\]
The Poisson condition gives $\kappa=0$ and therefore $Q=0$. With this substitution, the defect also satisfies $\operatorname{pr}_{\g}M_{\mathcal S}(y_2,f_1)=-p^2x$ and $\operatorname{pr}_{\g}M_{\mathcal S}(y_1,f_2)=q^2x$. Thus $p=q=0$ and $U=0$. Lemma~\ref{lem:gd-skew} gives preservation of the whole flag. The argument holds over every field $\K$ of characteristic zero.
\end{proof}

We can now describe all Poisson structures on the marked algebra.

\begin{theorem}\label{thm:a53-classification}
In the ordered basis $(y_1,y_2,x,f_1,f_2)$, all Poisson structures
on $(A_{5,3},T)$ are exactly the operators
\begin{equation}\label{eq:a53-all-S}
 \mathcal S_{a,b,c,d,u,v,t}=
 \begin{pmatrix}
  a&b&0&0&0\\
  c&d&0&0&0\\
  u&v&0&0&0\\
  0&-t&-u&-a&-c\\
  t&0&-v&-b&-d
 \end{pmatrix},\qquad
 a,b,c,d,u,v,t\in\K,\qquad ad-bc=0.
\end{equation}
Equivalently, the skew tensor solutions of the classical Yang--Baxter equation on $A_{5,3}$ are
\begin{equation}\label{eq:a53-all-r}
\begin{aligned}
 r_{\mathcal S}={}&f_1\wedge(ay_1+cy_2+ux)\\
                 &+f_2\wedge(by_1+dy_2+vx)+t\,f_1\wedge f_2,
                 \qquad ad-bc=0.
 \end{aligned}
\end{equation}
\end{theorem}

\begin{proof}
By Lemma~\ref{lem:a53-preservation}, every Poisson operator has
form \eqref{eq:gd-preserving-operator}.  Its blocks are
\[
 A=\begin{pmatrix}a&b\\c&d\end{pmatrix},\qquad
 Vy_1=ux,\quad Vy_2=vx,\qquad C=0,\qquad
 \eta(y_1,y_2)=t,
\]
which gives \eqref{eq:a53-all-S} before imposing the Poisson
condition.  Apply Theorem~\ref{thm:gd-poisson} with the data
\eqref{eq:a53-data}.  Since $\bLie$ and $\g$ are abelian and
$\psi=C=0$, one has
\[
 M_A=M_C=0,\qquad K=0,\qquad L=0.
\]
The remaining equation in the middle summand is
\begin{equation}\label{eq:a53-curvature-condition}
 L^{\mathrm{gen}}(y_1,y_2)
   =\alpha(Ay_1,Ay_2)=(ad-bc)x=0.
\end{equation}
The scalar equation $d_A\eta+\Xi_{A,V}=0$ is automatic because
$\Lambda^3\bLie^*=0$.  This proves both necessity and sufficiency
of the determinant condition on $A$, with $V$ and $\eta$ arbitrary.
\end{proof}

\subsubsection{Jacobi structures}
\label{subsec:a53-jacobi}

Assume that
\begin{equation}\label{eq:a53-E}
 E=\sigma f_1+\tau f_2\in I=Z(\mathfrak d),
 \qquad \sigma,\tau\in\K.
\end{equation}
Thus $[E,r]=0$ for every bivector $r$.
For $E=0$, Theorem~\ref{thm:a53-classification} already gives all
Jacobi pairs.  When $E\neq0$, preservation of $I$ no longer follows
from the Jacobi equations.

\begin{theorem}\label{thm:a53-jacobi-classification}
Fix $E=\sigma f_1+\tau f_2\neq0$ and define
\begin{equation}\label{eq:a53-w}
 w:=-\tau y_1+\sigma y_2.
\end{equation}
All Jacobi bivectors on $(A_{5,3},T)$ with Euler element $E$ are the
following two families:
\begin{align}
 r&=E\wedge Z,\qquad Z\in\mathfrak d,
       \label{eq:a53-jacobi-preserving}\\
 r&=w\wedge x+\zeta\wedge w+\chi\wedge x+t\,f_1\wedge f_2,
       \qquad \zeta,\chi\in I,\quad t\in\K.
       \label{eq:a53-jacobi-nonpreserving}
\end{align}
The operators corresponding to \eqref{eq:a53-jacobi-preserving}
preserve $I$, and their bivectors are Poisson.
The operators corresponding to \eqref{eq:a53-jacobi-nonpreserving}
do not preserve $I$, and their bivectors are not Poisson.
\end{theorem}

\begin{proof}
Start with an arbitrary $T$-skew operator in the form \eqref{eq:gd-six-blocks}. Its middle block vanishes because $\dim\g=1$. As in Lemma~\ref{lem:a53-preservation}, write
\[
 Qf_1=-\kappa y_2,\qquad Qf_2=\kappa y_1,\qquad
 Uf_1=px,\qquad Uf_2=qx.
\]
The $\g$-component of $M_{\mathcal S}(f_1,f_2)$ is $\kappa^2x$, while the right-hand side of \eqref{eq:gd-jacobi-operator-identity} lies in $I$. Thus $\kappa=0$, and every solution has the matrix form
\begin{equation}\label{eq:a53-jacobi-matrix}
 \mathcal S=
 \begin{pmatrix}
  a&b&-p&0&0\\
  c&d&-q&0&0\\
  u&v&0&p&q\\
  0&-t&-u&-a&-c\\
  t&0&-v&-b&-d
 \end{pmatrix}.
\end{equation}
Evaluation of \eqref{eq:gd-jacobi-operator-identity} on the basis pairs
gives exactly
\begin{align}
 p(p+\tau)&=0,&q(p+\tau)&=0,
       \label{eq:a53-jacobi-U1}\\
 p(q-\sigma)&=0,&q(q-\sigma)&=0,
       \label{eq:a53-jacobi-U2}\\
 (a+d)p+\tau a-\sigma b&=0,&
 (a+d)q+\tau c-\sigma d&=0,
       \label{eq:a53-jacobi-A}\\
 ad-bc+(p+\tau)u+(q-\sigma)v&=0.
       \label{eq:a53-jacobi-V}
\end{align}
The first four equations are the $\bLie$-components on $(y_1,x)$ and $(y_2,x)$. The last three come from the $\bLie$- and $\g$-components on $(y_1,y_2)$. Every other component is zero or repeats one of these equations. Since $E$ is central, the second Jacobi identity imposes no further condition. We now solve the system according to whether $(p,q)$ vanishes.

If $p=q=0$, equations \eqref{eq:a53-jacobi-A}--\eqref{eq:a53-jacobi-V}
become
\begin{equation}\label{eq:a53-jacobi-preserving-conditions}
 \tau a=\sigma b,\qquad \tau c=\sigma d,\qquad
 ad-bc=\sigma v-\tau u.
\end{equation}
The first two equations say that $Aw=0$ for
$A=\left(\begin{smallmatrix}a&b\\c&d\end{smallmatrix}\right)$.
Since $w\neq0$, they force $\det A=0$.
Regard $\varepsilon:=\sigma f_1+\tau f_2$ as a covector on $\bLie$.
Then $\Ker\varepsilon=\K w$, so
\[
 A=z\otimes\varepsilon\quad\text{for some }z\in\bLie,
 \qquad (u,v)=h(\sigma,\tau)\quad\text{for some }h\in\K.
\]
Formula~\eqref{eq:a53-all-r} gives
\[
 r=E\wedge(z+hx)+t\,f_1\wedge f_2.
\]
Because $E\neq0$ and $\dim I=2$, every element of $\Lambda^2I$
is $E\wedge z_I$ for some $z_I\in I$.
This proves \eqref{eq:a53-jacobi-preserving}, including its converse.
Its operator preserves $I$ because $p=q=0$.
Since $E$ is central, $[E\wedge Z,E\wedge Z]=0$ and
$E\wedge(E\wedge Z)=0$, so this family is Poisson and Jacobi.

Now suppose $(p,q)\neq(0,0)$. Equations~\eqref{eq:a53-jacobi-U1}--\eqref{eq:a53-jacobi-U2} give
\begin{equation}\label{eq:a53-jacobi-U-nonzero}
 p=-\tau,\qquad q=\sigma.
\end{equation}
The remaining equations reduce to
\[
 \sigma a+\tau c=0,\qquad
 \sigma b+\tau d=0,\qquad ad-bc=0.
\]
The first two equations mean that $\Imop A\subseteq\Ker\varepsilon=\K w$, which also gives the determinant equation. There is no restriction on $u,v,t$. Thus $A=w\otimes\ell$ for some $\ell\in\bLie^*$. Define $\zeta:=\ell(y_1)f_1+\ell(y_2)f_2$ and $\chi:=uf_1+vf_2$. The bivector corresponding to \eqref{eq:a53-jacobi-matrix} is
\begin{align*}
 r={}&(py_1+qy_2)\wedge x
      +f_1\wedge(ay_1+cy_2+ux)\\
    &+f_2\wedge(by_1+dy_2+vx)+t\,f_1\wedge f_2,
\end{align*}
which becomes \eqref{eq:a53-jacobi-nonpreserving}.

For a direct verification of this family, note that
$[w,x]=E$ by \eqref{eq:a53-brackets}.
Define
\[
 X_0=w+\chi,\qquad Y_0=x-\zeta,\qquad
 \beta=t\,f_1\wedge f_2+\chi\wedge\zeta\in\Lambda^2I.
\]
Then $r=X_0\wedge Y_0+\beta$,
$[X_0,Y_0]=E$, and $\beta$ is central for the Schouten bracket.
Since $\Lambda^3I=0$, one has $E\wedge\beta=0$, and therefore
\[
 [r,r]=2E\wedge X_0\wedge Y_0=2E\wedge r.
\]
The component $E\wedge w\wedge x$ of $E\wedge r$ is nonzero,
so this bivector is not Poisson.
Finally, \eqref{eq:a53-jacobi-U-nonzero} gives $U\neq0$, which proves
that its operator does not preserve $I$.
\end{proof}

When $E\neq0$, only \eqref{eq:a53-jacobi-preserving} has operators preserving the marked ideal. Thus Theorem~\ref{thm:gd-jacobi} applies to this family, but not to \eqref{eq:a53-jacobi-nonpreserving}. For an example of the latter, take $E=f_1$, $w=y_2$, and
\[
 r=y_2\wedge x.
\]
Then $[r,r]=2f_1\wedge r$ and $[f_1,r]=0$.
The corresponding operator satisfies $\mathcal S(f_2)=x\notin I$.

\subsection{Poisson structures from curvature}
\label{sec:inner-poisson}

Over $\K=\mathbb R$, the form $B$ defines a bi-invariant pseudo-Riemannian metric on the simply connected Lie group $G$ integrating $\g$. On left-invariant vector fields, its Levi--Civita connection is $\nabla_XY=\tfrac12[X,Y]$, and its curvature is
\[
R(X,Y)=-\frac14\ad_{[X,Y]}.
\]
The metric is locally symmetric, so $\nabla R=0$ (see \cite{Milnor1976}). Every curvature operator is therefore of the form $-\tfrac14\ad_c$, with $c=[X,Y]$. This reduces the question of when a curvature operator is Poisson to the same question for an inner derivation. We use the criteria in Section~\ref{sec:poisson-preliminaries} to answer it. The Jacobi case can be treated in the same way and is left to the reader. For $c\in\g$, define
\begin{equation}\label{eq:inner-Dc}
	D_c:=\ad_c,
	\qquad
	\h_c:=\Imop D_c=[c,\g].
\end{equation}
By invariance of $B$, $D_c$ is a $B$-skew derivation, though it need not be invertible. For $\lambda\in\K$, define $A_{\lambda,c}:=\lambda D_c$. If $\lambda\neq0$, then $D_c$ and $A_{\lambda,c}$ have the same image, and one satisfies the CYBE exactly when the other does.

%We first record formulas for $r_{A_{\lambda,c}}$ and its components.
%
%\begin{proposition}\label{prop:inner-tensor-formula}
%	Let $c\in\g$ and $\lambda\in\K$.  For $\alpha,\beta\in\g^*$, write
%	$\alpha^\sharp:=B^\sharp\alpha$ and $\beta^\sharp:=B^\sharp\beta$.  Then
%	\begin{equation}\label{eq:inner-r-formula}
%		r_{A_{\lambda,c}}(\alpha,\beta)
%		=\lambda B\bigl([c,\alpha^\sharp],\beta^\sharp\bigr)
%		=\lambda B\bigl(c,[\alpha^\sharp,\beta^\sharp]\bigr),
%	\end{equation}
%	and $r_{A_{\lambda,c}}^\sharp=A_{\lambda,c}\circ B^\sharp$.  In a basis
%	$\{e_i\}$ of $\g$, denote the dual basis by $\{e^i\}$ and define the constants $c_{ij}^{k}$ by
%	$[e_i,e_j]=c_{ij}^{k}e_k$.  Let $B_{ij}$ be the Gram matrix, and let $(B^{ij})$ be its inverse.
%	Write $c=c^ke_k$ and
%	\[
%	r_{A_{\lambda,c}}=\frac12r^{ab}_{\lambda,c}e_a\wedge e_b.
%	\]
%	Then
%	\begin{equation}\label{eq:inner-components}
%		(A_{\lambda,c})^{b}{}_{p}=\lambda c^kc_{kp}^{b},
%		\qquad
%		r^{ab}_{\lambda,c}=\lambda c^kc_{kp}^{b}B^{pa}.
%	\end{equation}
%	If $c=[X,Y]$, then $c^k=X^iY^jc_{ij}^{k}$.  Substitution in
%	\eqref{eq:inner-components} gives the corresponding formulas when $c$ is a commutator.
%\end{proposition}
%
%\begin{proof}
%	The first equality is \eqref{eq:rA-def}, and the second follows from invariance of $B$.
%	The identity for the sharp map is \eqref{eq:rsharp-is-A}.  Finally,
%	$B^\sharp(e^a)=B^{pa}e_p$, so
%	\[
%	r_{A_{\lambda,c}}^\sharp(e^a)
%	=\lambda c^kc_{kp}^{b}B^{pa}e_b,
%	\]
%	which gives \eqref{eq:inner-components} with the convention for coefficients used in the
%	statement.
%\end{proof}

The following cyclic identity relates the CYBE for $A_{\lambda,c}$ to the square of $D_c$. It follows by direct computation and holds for every $B$-skew derivation.

\begin{lemma}\label{lem:derivation-cyclic-identity}
	Let $D\in\Der(\g)$ be $B$-skew.  Then
	\begin{equation}\label{eq:derivation-cyclic-identity}
		\underset{X,Y,Z}{\mathrm{cyc}}\,
		B([D^2X,Y],Z)
		=-2\underset{X,Y,Z}{\mathrm{cyc}}\,
		B([DX,DY],Z).
	\end{equation}
\end{lemma}

%\begin{proof}
%	Let $\eta(X,Y,Z):=B([X,Y],Z)$ be the Cartan $3$-form.  The derivation property and
%	$B$-skewness of $D$ give
%	\[
%	\eta(DX,Y,Z)+\eta(X,DY,Z)+\eta(X,Y,DZ)=0.
%	\]
%	Apply the same infinitesimal action once more.  The terms in which $D$ occurs twice in one
%	argument have cyclic sum
%	$\underset{X,Y,Z}{\mathrm{cyc}}B([D^2X,Y],Z)$.  The mixed terms occur
%	with coefficient $2$.  Using the cyclic invariance of the alternating form
%	$\eta$, their sum is
%	$\underset{X,Y,Z}{\mathrm{cyc}}B([DX,DY],Z)$.  The resulting identity is exactly
%	\eqref{eq:derivation-cyclic-identity}.
%\end{proof}

Combining this identity with the leaf criterion gives several equivalent tests for an inner derivation to be Poisson.

\begin{theorem}\label{thm:inner-leaf-criterion}
	Let $(\g,B)$ be a quadratic Lie algebra, let $c\in\g$, and let $\lambda\in\K^\times$.  Define
	$D_c=\ad_c$ and $A_{\lambda,c}=\lambda D_c$.  The following conditions are equivalent:
	\begin{enumerate}
		\item $A_{\lambda,c}$ is a Poisson structure.
		\item
		\begin{equation}\label{eq:inner-operator-cybe}
			[[c,X],[c,Y]]=[c,[c,[X,Y]]],
			\qquad X,Y\in\g.
		\end{equation}
		\item
		\begin{equation}\label{eq:inner-image-cyclic}
			\underset{X,Y,Z}{\mathrm{cyc}}\,
			B([[c,X],[c,Y]],Z)=0,
			\qquad X,Y,Z\in\g.
		\end{equation}
		\item
		\begin{equation}\label{eq:inner-square-cyclic}
			\underset{X,Y,Z}{\mathrm{cyc}}\,
			B([D_c^2X,Y],Z)=0,
			\qquad X,Y,Z\in\g.
		\end{equation}
		\item The subspace $\h_c=[c,\g]$ is a Lie subalgebra, and the form
		\begin{equation}\label{eq:inner-leaf-form}
			\Omega_{\lambda,c}(A_{\lambda,c}X,A_{\lambda,c}Y)
			:=B(A_{\lambda,c}X,Y)
			=\lambda B(c,[X,Y])
		\end{equation}
		is closed on $\h_c$.
	\end{enumerate}
	When these conditions hold, $(\h_c,\Omega_{\lambda,c})$ is the characteristic leaf of $A_{\lambda,c}$.
\end{theorem}

\begin{proof}
	Since $D_c$ is a derivation,
	\[
	[D_cX,Y]+[X,D_cY]=D_c[X,Y].
	\]
	Consequently
	\begin{equation}\label{eq:inner-defect}
		M_{A_{\lambda,c}}(X,Y)
		=\lambda^2\bigl([D_cX,D_cY]-D_c^2[X,Y]\bigr).
	\end{equation}
Thus the CYBE is equivalent to \eqref{eq:inner-operator-cybe}, proving the equivalence of (1) and (2). By the calculation in Lemma~\ref{lem:operator-cybe-equivalence}, vanishing of the defect is equivalent to \eqref{eq:inner-image-cyclic}, so (1) and (3) are equivalent. Lemma~\ref{lem:derivation-cyclic-identity} gives the equivalence of (3) and (4). Finally, Proposition~\ref{prop:leaf-criterion}, applied to $A_{\lambda,c}$, gives the equivalence of (1) and (5), together with the formula and nondegeneracy of the form on the leaf.
\end{proof}

Either cyclic identity implies that $\h_c$ is a subalgebra, as required in Theorem~\ref{thm:inner-leaf-criterion}(5). Although this is not apparent from \eqref{eq:inner-image-cyclic} or \eqref{eq:inner-square-cyclic} alone, nondegeneracy of $B$ makes each scalar identity equivalent to the CYBE.

In particular, the theorem gives two sufficient conditions.

\begin{corollary}\label{cor:inner-sufficient}
The operator $A_{\lambda,c}$ is a Poisson structure under either of the following hypotheses:
	\begin{enumerate}
		\item $\h_c=[c,\g]$ is abelian.
		\item $D_c^2\g\subseteq Z(\g)$.  In particular, $D_c^2=0$ is sufficient.
	\end{enumerate}
\end{corollary}

To apply these criteria to nilpotent algebras, write the lower central series of $\mathfrak{g}$ as
\begin{equation}\label{eq:inner-lower-central}
	\g_1:=\g,
	\qquad
	\g_{j+1}:=[\g,\g_j],
	\qquad j\geq1,
\end{equation}
where $[\g_i,\g_j]\subseteq\g_{i+j}$. The algebra is $N$-\emph{step nilpotent} if $\g_{N+1}=0$ and $\g_N\neq0$.

\begin{corollary}\label{prop:inner-nilpotent-depth}
Let $(\g,B)$ be a quadratic Lie algebra nilpotent of step $N$, and let $c\in\g_s$. Then
	\begin{equation}\label{eq:inner-depth-image}
		\h_c=[c,\g]\subseteq\g_{s+1},
		\qquad
		[\h_c,\h_c]\subseteq\g_{2s+2}.
	\end{equation}
	In particular, if $2s+2\geq N+1$, then $\h_c$ is abelian.  Every nonzero scalar multiple of $D_c$ is then a
	Poisson structure and its leaf is
	$(\h_c,\Omega_{\lambda,c})$.
\end{corollary}

%\begin{proof}
%	The two containments in \eqref{eq:inner-depth-image} follow from the bracket inclusion for
%	the lower central series.  Condition \eqref{eq:inner-depth-condition} makes
%	$\g_{2s+2}$ vanish.
%	Corollaries~\ref{cor:inner-sufficient} and~\ref{cor:inner-nijenhuis-leaf} then give the
%	conclusion.
%\end{proof}

When $c\in\g_2$, in particular when $c$ is a commutator arising from curvature, its position in the lower central series gives the following consequences.

\begin{corollary}\label{cor:inner-low-step-hierarchy}
Let $(\g,B)$ be a quadratic Lie algebra nilpotent of step $N$, and let $c\in\g_2$. In particular, one may take $c=[X,Y]$. Define $D_c=\ad_c$ and $\h_c=[c,\g]$. Then
	\begin{enumerate}
		\item $\h_c\subseteq\g_3$.
		\item If $N\leq3$, then $\h_c\subseteq Z(\g)$ and
		$[D_c,\ad_X]=0$ for every $X\in\g$.
		\item If $N\leq4$, then $D_c^2=0$.  Thus the pair
			$(\pi_{\lambda D_c},\lambda D_c)$ is Poisson for every $\lambda\in\K$.
		\item If $N\leq5$, then $\h_c$ is abelian.  Therefore every $\lambda D_c$ with
		$\lambda\in\K^\times$ is a Poisson structure with
		characteristic leaf $(\h_c,\Omega_{\lambda,c})$.
	\end{enumerate}
	For $N\geq6$, the first possible obstruction is that
	$(\h_c,\Omega_{\lambda,c})$ fails to be symplectic, as characterized in
	Theorem~\ref{thm:inner-leaf-criterion}.
\end{corollary}

%\begin{proof}
%	The first claim is \eqref{eq:inner-depth-image} with $s=2$.  If $N\leq3$, then
%	$\g_3\subseteq Z(\g)$.  Moreover
%	$[\ad_X,D_c]=\ad_{[X,c]}=0$ because $[X,c]\in\g_3$.  If $N\leq4$, then
%	\[
%	D_c^2\g=[c,[c,\g]]\subseteq[\g_2,\g_3]\subseteq\g_5=0,
%	\]
%	and Corollary~\ref{cor:inner-nijenhuis-leaf} applies.  If $N\leq5$, then
%	$[\h_c,\h_c]\subseteq[\g_3,\g_3]\subseteq\g_6=0$, and the last assertion follows from
%	Corollaries~\ref{cor:inner-sufficient} and~\ref{cor:inner-nijenhuis-leaf}.
%\end{proof}
%
%The three-step nilpotent example in the final section
%applies these results to a quadratic Lie algebra of nilpotency step three.

For solvable algebras, use the derived series $\g^{(0)}:=\g$ and $\g^{(j+1)}:=[\g^{(j)},\g^{(j)}]$. The algebra $\g$ is solvable when this series terminates. In the metabelian case, the next corollary gives a Poisson structure for each inner derivation in the stated family.

%We call
%$(\g,B)$ a \emph{solvable quadratic Lie algebra} .  Solvability by
%itself does not provide an invariant metric.
%
%The following remark illustrates this obstruction with the affine algebra of dimension two.
%
%\begin{remark}\label{rem:affine-not-quadratic-inner}
%	Let $\mathfrak{aff}(\K)=\Span\{H,E\}$ with $[H,E]=E$.  If $B$ were invariant and symmetric,
%	then
%	\[
%	B(E,H)=B([H,E],H)=B(H,[E,H])=-B(H,E),
%	\]
%	so $B(E,H)=0$ in characteristic zero.  In addition,
%	\[
%	B(E,E)=B([H,E],E)=B(H,[E,E])=0.
%	\]
%	Thus $E$ is orthogonal to the whole algebra and $B$ is degenerate.  Therefore
%	$\mathfrak{aff}(\K)$ is solvable but not quadratic.
%\end{remark}

\begin{corollary}\label{cor:inner-metabelian}
Let $(\g,B)$ be a metabelian quadratic Lie algebra, so $[\g',\g']=0$ for $\g':=[\g,\g]$. For every $c\in\g'$ and $\lambda\in\K$, the map $A_{\lambda,c}=\lambda\ad_c$ is a Poisson structure. If it is nonzero, its characteristic leaf is $([c,\g],\Omega_{\lambda,c})$. The subalgebra $[c,\g]$ is abelian and lies in $\g'$.
\end{corollary}

%\begin{proof}
%	Because $\g'$ is an ideal, $[c,\g]\subseteq\g'$.  The metabelian hypothesis makes this image
%	abelian.  Apply Corollaries~\ref{cor:inner-sufficient}
%	and~\ref{cor:inner-nijenhuis-leaf}.
%\end{proof}

The oscillator algebra gives a family of solvable examples obtained by ordinary double extension. In the notation of Subsection~\ref{sec:oscillator-poisson}, the quadratic Lie algebra $\mathfrak{osc}_{\delta}=\mathbb Re\oplus\mathfrak v\oplus\mathbb Rf$ is solvable and satisfies
	\[
	\mathfrak{osc}_{\delta}^{(1)}\subseteq \mathfrak v\oplus\mathbb R f,
	\qquad
	\mathfrak{osc}_{\delta}^{(2)}\subseteq\mathbb R f,
	\qquad
	\mathfrak{osc}_{\delta}^{(3)}=0.
	\]
It is nilpotent if and only if $\delta$ is nilpotent. For $c=v+\mu f\in\mathfrak v\oplus\mathbb R f$ and $x=ae+w+bf\in\mathfrak{osc}_{\delta}$, we have $[c,x]=-a\delta v+B_{\mathfrak v}(\delta v,w)f$. If $\delta v=0$, then $\ad_c=0$. Otherwise, $[c,\mathfrak{osc}_{\delta}]=\Span\{\delta v,f\}$ is a two-dimensional abelian subalgebra. For $A_{\lambda,c}=\lambda\ad_c$ with $\lambda\neq0$,
	\begin{equation*}
		r_{A_{\lambda,c}}=\lambda \delta v\wedge f,
		\qquad
		\Omega_{A_{\lambda,c}}(\delta v,f)=\lambda^{-1},
		\qquad
		A_{\lambda,c}^2(e)=-\lambda^2B_{\mathfrak v}(\delta v,\delta v)f.
	\end{equation*}
	Thus every such map is a Poisson structure.  It has square zero precisely when $\delta v$ is isotropic.

\begin{remark}
	\normalfont
This example also shows why solvability alone does not replace the criterion in Theorem~\ref{thm:inner-leaf-criterion}. What matters is the image $[c,\g]$ and the form on it. Solvability provides useful families in which that image is small or abelian.
\end{remark}

Over $\mathbb R$, a nonabelian solvable example necessarily has an indefinite quadratic form.

\begin{remark}\label{rem:positive-definite-solvable}
	\normalfont
If a real quadratic Lie algebra is solvable and nonabelian, its invariant form must be indefinite. In particular, this holds for nonabelian nilpotent algebras. To see this, suppose $B$ is positive definite and $\g$ is solvable. Then $\ad(\g)$ is a solvable subalgebra of the compact Lie algebra of $B$-skew endomorphisms, so it is abelian. Thus $[\g,\g]\subseteq Z(\g)$. For $z=[x,y]\in[\g,\g]$, invariance and centrality give $B(z,z)=B([x,y],z)=B(x,[y,z])=0$. Positive definiteness forces $z=0$, so $\g$ is abelian. In the nilpotent case, one can also note that each $\ad_x$ is both nilpotent and skew-adjoint, and therefore zero.
\end{remark}

We conclude with the Lie--Poisson interpretation of an inner tensor over $\K=\mathbb R$. Use $B$ to identify $\g$ with $\g^*$ and transport the Lie--Poisson tensor to the vector space $\g$. With the convention $\Pi_{\mathrm{LP}}(\mu)(X,Y)=\mu([X,Y])$ for $\mu\in\g^*$ and $X,Y\in\g$, the resulting linear tensor is
\begin{equation}\label{eq:transported-lie-poisson-inner}
	\widetilde\Pi_{\mathrm{LP}}(z)(\alpha,\beta)
	:=B\bigl(z,[\alpha^\sharp,\beta^\sharp]\bigr),
	\qquad z\in\g,\quad \alpha,\beta\in\g^*.
\end{equation}

The next proposition identifies its value at a fixed point and compares the leaf form with the form on the corresponding coadjoint orbit.

\begin{proposition}\label{prop:inner-frozen-lie-poisson}
	Let $c\in\g$, $\lambda\in\mathbb R^\times$, and $A_{\lambda,c}=\lambda\ad_c$.  Then
	\begin{equation}\label{eq:inner-frozen-value}
		r_{A_{\lambda,c}}=\widetilde\Pi_{\mathrm{LP}}(\lambda c).
	\end{equation}
For $\mu_{\lambda,c}:=\lambda B(c,\cdot)\in\g^*$, we have $-\ad_X^*\mu_{\lambda,c}=(A_{\lambda,c}X)^\flat$. Thus the identification given by $B$ sends the tangent space of the coadjoint orbit through $\mu_{\lambda,c}$ onto $\h_c=[c,\g]$. The Kirillov--Kostant--Souriau form satisfies
	\begin{equation}\label{eq:inner-KKS-form}
		\begin{split}
			\omega_{\mathrm{KKS},\mu_{\lambda,c}}
			\bigl((A_{\lambda,c}X)^\flat,(A_{\lambda,c}Y)^\flat\bigr)
			&=\mu_{\lambda,c}([X,Y])\\
			&=B(A_{\lambda,c}X,Y)\\
			&=\Omega_{A_{\lambda,c}}(A_{\lambda,c}X,A_{\lambda,c}Y).
		\end{split}
	\end{equation}

\end{proposition}

The word ``frozen'' distinguishes this construction from the full Lie--Poisson tensor. That tensor varies linearly over $\g^*$, whereas $r_{A_{\lambda,c}}$ is its value at $\mu_{\lambda,c}$, used as the value at the identity of a left-invariant bivector on a group. Even if $\h_c$ is abelian and its own Lie--Poisson tensor is zero, the frozen ambient tensor may be nonzero. When $A_{\lambda,c}$ is Poisson, the KKS value in \eqref{eq:inner-KKS-form} is also a symplectic form in the Chevalley--Eilenberg complex of $\h_c$. This additional closedness is precisely the condition in Theorem~\ref{thm:inner-leaf-criterion}.

A left-invariant bivector $\Pi$ on a Lie group $G$ is multiplicative if and only if it vanishes. Applying this observation to the inner tensors gives the following statement.

\begin{proposition}\label{prop:inner-nonmultiplicative}
For	$A_{\lambda,c}=\lambda\ad_c$ with $\lambda\neq0$, the left-invariant bivector	$\pi_{A_{\lambda,c}}$ is multiplicative if and only if $c\in Z(\g)$.  
\end{proposition}

Therefore, every nonzero Poisson tensor in this family is nonmultiplicative.

\end{document}